\RequirePackage{etoolbox} 
\csdef{input@path}{{style/}{graphics/}}
\documentclass[aos,MSNbibl,seceqn,dvips]{arximspdf}
\usepackage{graphicx}
\doi{10.1214/14-AOS1294}% Updated by VTEXPTS2LaTeX.exe, 29.01.2015
\volume{43}
\issue{2}
\pubyear{2015}
\firstpage{675}
\lastpage{712}
\docsubty{FLA}

\makeatletter
\newcommand{\lleft}{\left}
\newcommand{\rrvert}{\vert}
\newcommand{\rright}{\right}
\newcommand{\rrVert}{\Vert}
\newcommand{\llvert}{\vert}
\newcommand{\llVert}{\Vert}
\newtheorem{theorem}{Theorem}[section]
\newtheorem{lemma}{Lemma}[section]
\newproclaim{remark}{Remark}[section]
\newtheorem{proposition}{Proposition}[section]
\newproclaim{assumption}{Assumption}[section]
\newproclaim{Outline}{Outline}
\newproclaim{TSI}{Transformation to serial independence}
\newproclaim{TT}{Transformation of $\bC_\tau$}
\newproclaim{AT}{Approximating equations for the Stieltjes kernel $K_\tau$}
\newcommand{\bA}{\mathbf{A}}
\newcommand{\bC}{\mathbf{C}}
\newcommand{\bD}{\mathbf{D}}
\newcommand{\bF}{\mathbf{F}}
\newcommand{\bH}{\mathbf{H}}
\newcommand{\bI}{\mathbf{I}}
\newcommand{\bL}{\mathbf{L}}
\newcommand{\bP}{\mathbf{P}}
\newcommand{\bQ}{\mathbf{Q}}
\newcommand{\bR}{\mathbf{R}}
\newcommand{\bS}{\mathbf{S}}
\newcommand{\bU}{\mathbf{U}}
\newcommand{\bV}{\mathbf{V}}
\newcommand{\bW}{\mathbf{W}}
\newcommand{\bX}{\mathbf{X}}
\newcommand{\bZ}{\mathbf{Z}}
\newcommand{\bLambda}{\bolds{\Lambda}}
\newcommand{\bpsi}{\bolds{\psi}}
\newcommand{\oX}{\bar{X}}
\newcommand{\obC}{\bar{\bC}}
\newcommand{\obR}{\bar{\bR}}
\newcommand{\obV}{\bar{\bV}}
\newcommand{\cH}{\mathcal{H}}
\newcommand{\kK}{\mathfrak{K}}
\newcommand{\tK}{\tilde{K}}
\newcommand{\tX}{\tilde{X}}
\newcommand{\tZ}{\tilde{Z}}
\newcommand{\tcH}{\tilde{\cH}}
\newcommand{\ts}{\tilde{s}}
\newcommand{\tbC}{\tilde{\mathbf{C}}}
\newcommand{\tbL}{\tilde{\mathbf{L}}}
\newcommand{\tbR}{\tilde{\mathbf{R}}}
\newcommand{\tbX}{\tilde{\mathbf{X}}}
\newcommand{\tbZ}{\tilde{\mathbf{Z}}}
\newcommand{\bblambda}{\bar{\lambda}}
\makeatother

\begin{document}
\begin{frontmatter}

\title{On the Mar\v{c}enko--Pastur law for linear time series}
\runtitle{Mar\v{c}enko--Pastur law for time series}

\begin{aug}
% Corresponding author: Debashis Paul - debashis@wald.ucdavis.edu% Updated by VTEXPTS2LaTeX.exe, 30.01.2015 14:33
%Updated by VTEXPTS2LaTeX.exe, 29.01.2015 10:22
\author[A]{\fnms{Haoyang}~\snm{Liu}\thanksref{M1}\ead[label=e1]{liuhy@berkeley.edu}},
\author[B]{\fnms{Alexander}~\snm{Aue}\thanksref{M2,T2}\ead[label=e2]{aaue@ucdavis.edu}}
\and
\author[B]{\fnms{Debashis}~\snm{Paul}\corref{}\thanksref{M2,T3}\ead[label=e3]{debpaul@ucdavis.edu}}
\runauthor{H. Liu, A. Aue and D. Paul}
\affiliation{University of California, Berkeley\thanksmark{M1} and
University of California, Davis\thanksmark{M2}}
\address[A]{H. Liu\\
Haas School of Business\\
University of California\\
Berkeley, California 94720\\
USA\\
\printead{e1}}
\address[B]{A. Aue\\
D. Paul\\
Department of Statistics\\
University of California\\
Davis, California 95656\\
USA\\
\printead{e2}\\
\phantom{E-mail: }\printead*{e3}}
\end{aug}
\thankstext{T2}{Supported in part by NSF Grants DMS-12-09226,
DMS-13-05858 and DMS-14-07530.}
\thankstext{T3}{Supported in part by NSF Grants DMR-10-35468,
DMS-11-06690 and DMS-14-07530.}

% HISTORY:
%
\received{\smonth{10} \syear{2013}}% Updated by VTEXPTS2LaTeX.exe,
%29.01.2015 10:22
%
\revised{\smonth{11} \syear{2014}}% Updated by VTEXPTS2LaTeX.exe,
%29.01.2015 10:22

% ABSTRACT
%
\begin{abstract}
This paper is concerned with extensions of the classical Mar\v{c}enko--Pastur
law to time series. Specifically, $p$-dimensional linear processes are
considered
which are built from innovation vectors with independent, identically
distributed
(real- or complex-valued) entries possessing zero mean, unit variance
and finite fourth moments. The
coefficient matrices of the linear process are assumed to be simultaneously
diagonalizable. In this setting, the limiting behavior of the empirical spectral
distribution of both sample covariance and symmetrized sample autocovariance
matrices is determined in the high-dimensional setting $p/n\to c\in
(0,\infty)$
for which dimension $p$ and sample size $n$ diverge to infinity at the
same rate.
The results extend existing contributions available in the literature
for the
covariance case and are one of the first of their kind for the
autocovariance case.
\end{abstract}

% KEYWORDS
\end{frontmatter}

%s1 #&#
\section{Introduction}\label{secintro}

%- contemporaneous dependence (present values of the observed vectors
%may be correlated,
%but observations taken at different points in time are independent)
%References:

%- dependence induced by lag structure (past influences present values)
%References:

%- all work has been on covariance matrix $\mathbf{S}$ (cite and
%describe)

%- but prediction relies heavily on understanding the autocovariance
%structure

One of the exciting developments in statistics during the last
decade has been the development of the theory and methodologies for dealing
with high-dimensional data. The term \textit{high dimension} is primarily
interpreted as meaning that the dimensionality of the observed multivariate
data is comparable to the available number of replicates or subjects on which
the measurements on the different variables are taken. This is often expressed
in the asymptotic framework as $p/n \to c > 0$, where $p$ denotes the dimension
of the observation vectors (forming a triangular array) and $n$ the
sample size. Much of this development centered on understanding the
behavior of
the sample covariance matrix and especially its eigenvalues and eigenvectors,
due to their role in dimension reduction, in estimation of population
covariances and as building block of numerous inferential procedures for
multivariate data. Comprehensive reviews of this topic can be found in
Johnstone~\cite{j07} and Paul and Aue~\cite{pa13}.

One most notable high-dimensional phenomena associated with sample covariance
matrices is that the sample eigenvalues do not converge to their population
counterparts if dimension and sample sizes remain comparable even as
the sample
size increases. A formal way to express this phenomenon is through the
use of
the \textit{empirical spectral distribution} (ESD), that is, the empirical
distribution of the eigenvalues of the sample covariance matrix. The celebrated
work of Mar\v{c}enko and Pastur \cite{mp67} shows that if one studies
a triangular array of random
vectors $X_1,\ldots,X_n$, whose components form independent,
identically distributed (i.i.d.) random variables with zero mean, unit variance
and finite fourth moment, then as $p,n \to\infty$ such that $p/n \to
c \in
(0,\infty)$, the ESD of $\mathbf{S} = n^{-1} \sum_{j=1}^n X_jX_j^T$ converges
almost surely to a nonrandom probability distribution known as the
Mar\v{c}enko--Pastur distribution. Since this highly influential
discovery a
large body of literature under the banner of random matrix theory (RMT) has
been developed to explore the properties of the eigenvalues and
eigenvectors of
large random matrices. One may refer to Anderson et~al. \cite{agz09},
%%Guionnet and Zeitouni~\cite{agz09},
Bai and Silverstein~\cite{bs10} and Tao~\cite{t12} to study various
aspects of this
literature.

Many important classes of high-dimensional data, particularly
those arising in signal processing, economics and finance, have the feature
that in addition to the dimensional correlation, the observations are
correlated in time. Classical models for time series often assume a stationary
correlation structure and use spectral analysis methods or methods
built on
the behavior of the sample autocovariance matrices for inference and
prediction purposes.
%to infer about such type of data as well as for prediction of future
%observations.
In spite of this, to our knowledge, no work exists that analyzes the behavior
of the sample autocovariance matrices of a time series from a random
matrix perspective, even though Jin et~al. \cite{jwbnh14} have dealt
recently covered autocovariance matrices in the independent case.
A~striking observation is that, in the high-dimensional scenario, the
distribution of the eigenvalues of the symmetrized sample
autocovariance of a
given lag order tends to stabilize to a nondegenerate distribution
even in the
setting where the observations are i.i.d. This raises questions about the
applicability of sample autocovariance matrices as diagnostic tools for
determining the nature of temporal dependence in high-dimensional settings.
Thus a detailed study of the phenomena associated with the behavior of
the ESD
of the sample autocovariance matrices when the observations have both
dimensional and temporal correlation is of importance to gain a better
understanding of the ways in which the dimensionality affects the
inference for
high-dimensional time series.

All the existing work on high-dimensional time series dealing
with the limiting behavior of the ESD focuses on the sample covariance matrix
of the data when $X_1,\ldots,X_n$ are $p$-dimensional observations
recorded in
time and $p,n \to\infty$ such that $p/n \to c \in(0,\infty)$. This includes
the works of Jin et~al. \cite{jwml09}, who assume the process
$(X_t\dvtx t\in\mathbb{Z})$ has
i.i.d. rows with each row following a causal ARMA process. Pfaffel and
Schlemm~\cite{ps12} and Yao~\cite{y12} extend this framework to the setting
where the rows are arbitrary i.i.d. stationary processes with short-range
dependence. Zhang~\cite{z06}, Paul and Silverstein~\cite{ps09} and El
Karoui~\cite{e09}, under slightly different assumptions, consider the
limiting behavior of the ESD of the sample covariance when the data matrices
are of the form $\mathbf{A}^{1/2} \mathbf{Z}\mathbf{B}^{1/2}$ where
$\mathbf{A}$ and $\mathbf{B}$ are positive semidefinite matrices, and
$\mathbf{Z}$ has i.i.d. entries with zero mean, unit variance and
finite fourth
moment. This model is known as the separable covariance model, since the
covariance of the data matrix is the Kronecker product of $\mathbf{A}$~and~$\mathbf{B}$.
If the rows indicate spatial coordinates and columns indicate
time instances, then this model implies that spatial (dimensional) and temporal
dependencies in the data are independent of each other. The work of
this paper is also partly
related to the results of Hachem et~al. \cite{hln05}, %Loubaton and
%Nazim~\cite{hln05}
who prove the existence of the limiting ESD for sample covariance of
data matrices that are rectangular slices from a bistationary Gaussian
process on
$\mathbb{Z}^2$.

In this paper, the focus is on a class of time series known as
linear processes [or $\mathrm{MA}(\infty)$ processes]. %given by the
%representation
%X_t =\sum_{\ell=0}^\infty\bA_\ell Z_{t-\ell}, \qquad t \in\mathbb{Z},
%where $(\bA_\ell\dvtx\ell\in\mathbb{N})$ are $p\times p$ matrices, $
%identity matrix, and $(Z_t\dvtx t\in\mathbb{Z})$ are $p$-dimensional
%random
%vectors (innovations) with i.i.d. entries $Z_{tj}$ (real- or
%complex-valued)
%with zero mean, unit variance and finite fourth moment. It is assumed
%that the
%matrices $\bA_\ell$ are symmetric (in the real case) or Hermitian (in
%the
%complex case) and are simultaneously diagonalizable. Moreover, the
%stability
%requirement $\sum_{\ell=1}^\infty\ell\left\Vert\bA_\ell\right\Vert<
%imposed,
%where $\left\Vert\cdot\right\Vert$ denotes operator norm.
The assumptions to be imposed in Section~\ref{secmain} imply that, up
to an
unknown rotation, the coordinates of the linear process, say
$(X_t\dvtx t\in\mathbb{Z})$,
are uncorrelated stationary linear processes with short range dependence.
Extending the work of Jin et~al. \cite{jwbnh14} to the time series
case, the goal is to relate
the behavior of the ESD of the lag-$\tau$ symmetrized sample autocovariances,
defined as $\mathbf{C}_\tau= (2n)^{-1} \sum_{t=1}^{n-\tau} (X_t
X_{t+\tau}^* +
X_{t+\tau}X_t^*)$, with $^*$ denoting complex conjugation, to that of the
spectra of the coefficient matrices of the linear process %$(\bA_\ell
when $p,n \to\infty$ such that $p/n \to c \in(0,\infty)$. This
requires assuming
certain stability conditions on the joint distribution of the
eigenvalues of
the coefficient matrices which are described later. The class of models under
study here includes the class of causal autoregressive moving average (ARMA)
processes of finite orders satisfying the requirement that the coefficient
matrices are simultaneously diagonalizable and the joint empirical distribution
of their eigenvalues (when diagonalized in the common orthogonal or unitary
basis), converges to a finite-dimensional distribution. The results are
expressed in terms of the \textit{Stieltjes transform} of the ESD of
the sample
autocovariances. %The Stieltjes transform of the ESD of a Hermitian
%matrix can
%be expressed as the trace of the resolvent of the matrix normalized by
%the dimension.
Specifically, it is shown that the ESD of the symmetrized sample
autocovariance matrix of any lag order converges to a nonrandom probability
distribution on the real line whose Stieltjes transform can be
expressed in
terms a unique \textit{Stieltjes kernel}. The definition of the Stieltjes
kernel involves integration with respect to the limiting joint empirical
distribution of the eigenvalues of the coefficient matrices as well as the
spectral density functions of the one-dimensional processes that
correspond to
the coordinates of the process $(X_t\dvtx t\in\mathbb{Z})$, after
rotation in the common unitary or
orthogonal matrix that simultaneously diagonalizes the coefficient matrices.
%$(\bA_\ell\dvtx\ell\in\mathbb{N})$.
Thus this result neatly ties the
dimensional correlation, captured by the eigenvalues of the coefficient
matrices, with the temporal correlation, captured by the spectral density
of the coordinate processes.

The main contributions of this paper are the following: (i) A framework is
provided for analyzing the behavior of symmetrized autocovariance
matrices of
linear processes; (ii) for linear processes satisfying appropriate
regularity conditions,
a~concrete description of the limiting Stieltjes transform is
given in terms of the limiting joint ESD of the coefficient matrices
and the
spectral density of the coordinate processes after a rotation of the
coordinates of the observation. Extensions to these main results are
(iii) the
characterization of the behavior of the ESD of autocovariances of
linear filters
applied to the observed process; (iv) the description of the ESDs of a
class of tapered
estimates of the spectral density operator of the observed process that
can be
used to analyze the long-run variance and spectral coherence of the process.
These contributions surpass the work in the existing literature dealing with
high-dimensionality effects for time series in two different ways.
First, the
class of time series models that are analyzed in detail encompasses the setting
of stationary i.i.d. rows studied by Jin et~al. \cite{jwml09},
Pfaffel and
Schlemm~\cite{ps12} and Yao~\cite{y12}, as well as the setting of separable
covariance structure studied by Zhang~\cite{z06}, Paul and
Silverstein~\cite{ps09} and El Karoui~\cite{e09}. The proofs of the main
results also require more involved arguments. They are partly related
to the
constructions in Hachem et~al. \cite{hln05}, but additional technical arguments
are needed to go beyond Gaussanity. The results are also related to the work
of Hachem et~al. \cite{hln07}, who studies limiting spectral
distributions of covariance matrices for data with a given variance profile.
The connection is through the fact that after an approximation of lag operators
by circulant shift matrices, and appropriate row and column rotations,
the data matrix in our setting can be equivalently expressed as a
matrix with independent entries and
with a variance profile related to the spectral densities of the
different coordinates of the
time series. Second, the framework allows for a
unified analysis of the ESD of symmetrized autocovariance matrices of
all lag
orders as well as that of the tapered spectral density operator. None
of the
existing works deals with the behavior of autocovariances for time series
(note again that Jin et~al. \cite{jwbnh14} treat the i.i.d. case),
and this analysis
requires a nontrivial variation of the arguments used for dealing with the
Stieltjes transform of the sample covariance matrix. Moreover, even
though we
stick to the setting where the coefficient matrices are Hermitian and
simultaneously diagonalizable, the main steps in the derivation,
especially the
construction of a ``deterministic equivalent'' of the resolvent of the
symmetrized autocovariance matrix, is very general and can be applied to
linear processes with structures that go beyond the settings studied
in this paper, for example, when the simultaneous diagonalizability of the
coefficient matrices is replaced by a form of simultaneous block
diagonalizability,
even though the latter is not pursued in this paper due to lack of
clear statistical
motivation. The existence and uniqueness of the limits of the resulting
equations and their
solutions is the key to establishing the existence of liming ESDs of
the autocovariances. This step
requires certain regularity conditions on the coefficient matrices and
is not
pursued beyond the setting described in Section~\ref{secmain}. A number
of potential applications, for example, to problems in signal processing,
and dynamic and static factor models, are discussed in Section~\ref
{secexandapp}.

The remaining sections of the paper are organized as follows.
Extensions of
the main results in Section~\ref{secmain} are discussed in
Section~\ref{secex}.
The outcomes of a small simulation study are reported in Section~\ref
{secsimulations},
while the proofs of the main results are provided in Sections~\ref
{secstructure}--\ref{secpreal}.
Several technical lemmas are collected in the online Supplemental
Material (SM) \cite{lap-sm}.

%s2 #&#
\section{Main results}
\label{secmain}

Let $\mathbb{Z}$ denote the set of integers. A sequence of random vectors
$(X_t\dvtx t\in\mathbb{Z})$ with values in $\mathbb{C}^p$ is called
a linear
process or moving average process of order infinity, abbreviated by the acronym
MA($\infty$), if it has the representation
%
%e2.1 #&#
\begin{equation}
\label{eqmainfty} X_t=\sum_{\ell=0}^\infty
\bA_\ell Z_{t-\ell}, \qquad t\in\mathbb{Z},
\end{equation}
where $(Z_t\dvtx t\in\mathbb{Z})$ denotes a sequence of independent,
identically
distributed \mbox{$p$-}dimensional random vectors whose entries are
independent and satisfy $\mathbb{E}[Z_{jt}]=0$,
$\mathbb{E}[\llvert Z_{jt}\rrvert^2]=1$ and $\mathbb{E}[\llvert
Z_{jt }\rrvert^4]<\infty$,
where $Z_{jt}$ denotes the
$j$th coordinate of $Z_t$. In the complex-valued case this is meant
as $\mathbb{E}[\operatorname{Re}(Z_{jt})^2]=\mathbb{E}[\operatorname
{Im}(Z_{jt})^2]=1/2$. It is also assumed
that real and imaginary parts are independent. Let further $\bA_0=\bI$,
the identity matrix. To ensure finite fourth moments for $(X_t\dvtx
t\in\mathbb{Z})$
and a sufficiently fast decaying weak dependence structure,
Assumption~\ref{asmainfty}
below lists several additional conditions imposed on the coefficient
matrices $\bA_\ell$.

The results presented in this paper are concerned with the behavior of
the symmetrized
lag-$\tau$ sample autocovariances
\[
\bC_\tau=\frac{1}{2n}\sum_{t=1}^{n-\tau}
\bigl(X_tX_{t+\tau
}^*+X_{t+\tau}X_t^*
\bigr), \qquad\tau\in\mathbb{N}_0, %
\]
assuming observations for $X_1,\ldots,X_n$ are available. For $\tau
=0$, this definition
gives the covariance matrix $\bS=\bC_0$ discussed in the
\hyperref[secintro]{Introduction}. Note that in
order to make predictions in the linear process setting, it is
imperative to
understand the second-order dynamics which are captured in the
population autocovariance
matrices $\bolds{\Gamma}_\tau=\mathbb{E}[X_{t+\tau}X_t^*]$, $\tau
\in\mathbb{N}_0$, as all
of the popular prediction algorithms such as the Durbin--Levinson and
innovations
algorithms are starting from there; see, for example, L\"utkepohl \cite{l06}.
The set-up in (\ref{eqmainfty}) provides a (strictly) stationary
process and consequently
the definition of $\bolds{\Gamma}_\tau$ does not depend on the
value of $t$. The main
goal of this paper is to analyze the behavior of the matrices $\bC
_\tau$, which can be
viewed as a special sample counterpart to the corresponding $\bolds
{\Gamma}_\tau$, in
the high-dimensional setting for which $p=p(n)$ is a function of the
sample size such that
%
%e2.2 #&#
\begin{equation}
\label{eqpn} \lim_{n\to\infty}\frac{p}n= c\in(0,\infty),
\end{equation}
thereby extending the above mentioned Mar\v{c}enko--Pastur-type
results to more
general time series models and to autocovariance matrices. We can
weaken requirement (\ref{eqpn})
to ``$p/n$ bounded away from zero and infinity,'' in which case, the
asymptotic results
hold for subsequences $(p_k,n_k)$ satisfying $p_k/n_k$ converging to a
positive constant
$c_k$, provided that the structural assumptions on the model continue
to hold.
Let then $\hat
F_\tau$ denote the empirical spectral distribution (ESD) of $\bC_\tau
$ given by
\[
\hat F_\tau(\sigma) =\frac{1}p\sum
_{j=1}^p\mathbb{I}_{\{\sigma_j\leq\sigma\}},
\]
where $\sigma_1,\ldots,\sigma_p$ are the eigenvalues of
$\bC_\tau$. The proof techniques for establishing large-sample
results about
$\hat F_\tau$ are based on exploiting convergence properties of Stieltjes
transforms, which continue to play an important role in verifying theoretical
results in RMT; see, for example, Paul and Aue \cite{pa13} for a recent
summary. The Stieltjes transform of a distribution function $F$ on the real
line is the function
\[
s_F\dvtx\mathbb{C}^+\to\mathbb{C}^+,\qquad z\mapsto s_F(z)=
\int\frac
{1}{\sigma-z}\,dF(\sigma),
\]
where $\mathbb{C}^+=\{x+iy\dvtx x\in\mathbb{R},y>0\}$ denotes the upper
complex half plane. It can be shown that $s_F$ is analytic on $\mathbb
{C}^+$ and
that the distribution function $F$ can be reconstructed from $s_F$
using an
inversion formula; see \cite{pa13}. In order to make statements about
$\hat
F_\tau$, the following additional assumptions on the coefficient matrices
$\bA_\ell$ are needed. Let $\mathbb{N}$ and $\mathbb{N}_0$ denote
the positive
and nonnegative integers, respectively.

%as2.1 #&#
\begin{assumption}
\label{asmainfty}
(a) The matrices $(\bA_\ell\dvtx\ell\in\mathbb{N}_0)$ are
simultaneously diagonalizable random Hermitian matrices, independent
of\/
$(Z_t\dvtx t\in\mathbb{Z})$ and satisfying $\llVert\bA_\ell\rrVert
\leq
\bblambda_{\bA_\ell}$
for all $\ell\in\mathbb{N}_0$ and large $p$ with
\[
\sum_{\ell=0}^\infty\bblambda_{\bA_\ell}
\leq\bblambda_\bA<\infty\quad\mbox{and}\quad\sum
_{\ell=0}^\infty\ell\bblambda_{\bA_\ell}\leq
\bblambda_\bA^\prime<\infty.
\]
Note that one can set $\bblambda_{\bA_0}=1$.

(b) There are continuous functions
$f_\ell\dvtx\mathbb{R}^m\to\mathbb{R}$, $\ell\in\mathbb{N}_0$,
%$(f_\ell\dvtx\ell\in\mathbb{N}_0)$,
%a random distribution $F_p^\bA$ assigning
%(subject to multiplicities) mass $1/p$ to each of $p$ points, say,
such that, for every~$p$, there is a set of points
$\lambda_1,\ldots,\lambda_p\in\mathbb{R}^m$, not necessarily
distinct, and a
unitary $p\times p$ matrix $\bU$ such that
\[
\bU^*\bA_\ell\bU=\operatorname{diag}\bigl(f_\ell(
\lambda_1),\ldots,f_\ell(\lambda_p)\bigr),
\qquad\ell\in\mathbb{N},
\]
and $f_0(\lambda)=1$. [Note that the functions $f_\ell$ are allowed
to depend on $p=p(n)$ as long as they converge to continuous functions as
$n\to\infty$ uniformly.]

(c) With probability one, $F_p^\bA$, the ESD of
$\{\lambda_1,\ldots,\lambda_p\}$, converges weakly to a nonrandom probability
distribution function $F^\bA$ on $\mathbb{R}^m$ as $p\to\infty$.
\end{assumption}

Let $\bA=[\bA_0\dvtx\bA_1\dvtx\cdots]$ denote the matrix collecting
the coefficient
matrices of the linear process $(X_t\dvtx t\in\mathbb{Z})$. Define the
transfer functions
%
%e2.3 #&#
\begin{equation}
\label{eqtransferpsi} \psi(\lambda,\nu)=\sum_{\ell=0}^\infty
e^{i\ell\nu}f_\ell(\lambda) \quad\mbox{and}\quad\psi(\bA,\nu)=
\sum_{\ell=0}^\infty e^{i\ell\nu}
\bA_\ell,
\end{equation}
as well as the power transfer functions
\[
h(\lambda,\nu)=\bigl\llvert\psi(\lambda,\nu)\bigr\rrvert^2\quad\mbox{and}\quad\cH(\bA,\nu)=\psi(\bA,\nu)\psi(\bA,\nu)^*.
\]
Note\vspace*{1pt} that the contribution of the temporal dependence of the underlying time
series on the asymptotic behavior of $\hat F_\tau$ is quantified through
$h(\lambda,\nu)$. Specifically, $h(\lambda_j,\nu)$ with $\lambda
_j$ as in part
(b) of Assumption~\ref{asmainfty} is (up to normalization) the spectral
density of the $j$th coordinate of the process rotated with the help of the
unitary matrix $\bU$. With these definitions, the main results of this
paper can be stated as follows.

%th2.1 #&#
\begin{theorem}
\label{thesdmainfty} If\vspace*{1pt} a complex-valued linear process $(X_t\dvtx
t\in\mathbb{Z})$
with independent, identically\vspace*{1pt} distributed $Z_{jt}$, $\mathbb
{E}[Z_{jt}]=0$, $\mathbb{E}[\operatorname{Re}(Z_{jt})^2]=\mathbb
{E}[\operatorname{Im}(Z_{jt})^2]=1/2$,
$\operatorname{Re}(Z_{jt})$ and $\operatorname{Im}(Z_{jt})$ are independent,
and $\mathbb{E}[\llvert Z_{jt}\rrvert^4]<\infty$, satisfies
Assumption~\ref{asmainfty}, then, with probability one and in the
high-dimensional setting (\ref{eqpn}), $\hat F_\tau$~converges to a nonrandom
probability distribution $F_\tau$ with Stieltjes transform $s_\tau$ determined
by the equation
%
%e2.4 #&#
\begin{equation}
\label{eqstau} s_\tau(z)=\int\biggl[\frac{1}{2\pi}\int
_0^{2\pi}\frac{\cos(\tau
\nu)h(\lambda,\nu)}{1+c\cos(\tau\nu)K_\tau(z,\nu)}\,d\nu-z
\biggr]^{-1}\,dF^\bA(\lambda),
\end{equation}
where $K_\tau\dvtx\mathbb{C}^+\times[0,2\pi]\to\mathbb{C}^+$ is a
Stieltjes kernel; that is, $K_\tau(\cdot,\nu)$ is the Stieltjes
transform of a
measure with total mass $m_\nu=\int h(\lambda,\nu)\,dF^\bA(\lambda)$
for every
fixed $\nu\in[0,2\pi]$, whenever $m_\nu>0$. Moreover, $K_\tau$ is
the unique solution of
%
%e2.5 #&#
\begin{eqnarray}\label{eqktau}
&& K_\tau(z,\nu)
\nonumber\\[-8pt]\\[-12pt]\nonumber
&&\qquad =\int\biggl[\frac{1}{2\pi}\int
_0^{2\pi} \frac{\cos(\tau\nu^\prime)h(\lambda,\nu^\prime)}{1+c\cos(\tau
\nu^\prime)K_\tau(z,\nu^\prime)}\,d\nu^\prime-z
\biggr]^{-1}h(\lambda,\nu) \,dF^\bA(\lambda),\hspace*{-20pt}
\end{eqnarray}
subject to the restriction that $K_\tau$ is a Stieltjes kernel.
Otherwise, if
$m_\nu=0$, then $K_\tau(z,\nu)$ is identically zero on $\mathbb
{C}^+$ and so
still satisfies (\ref{eqktau}).
\end{theorem}

%th2.2 #&#
\begin{theorem}
\label{coesdmainfty} If a real-valued\vspace*{3pt} linear process $(X_t\dvtx
t\in\mathbb{Z})$ with independent, identically distributed
real-valued\vspace*{3pt} $Z_{jt}$, $\mathbb{E}[Z_{jt}]=0$, $\mathbb
{E}[Z_{jt}^2]=1$ and $\mathbb{E}[Z_{jt}^4]<\infty$, satisfies
Assumption~\ref{asmainfty} with real symmetric coefficient matrices
$(\bA_\ell\dvtx\ell\in\mathbb{N}_0)$, then the result of Theorem~\ref{thesdmainfty} is retained.
\end{theorem}

%re2.1 #&#
\begin{remark}
When each coefficient matrix $\bA_\ell$ is a multiple of the identity
matrix, that is, $\mathbf{A}_\ell=
\alpha_\ell I_p$ where $(\alpha_\ell\dvtx\ell\in\mathbb{N})$ is
a sequence of real numbers
satisfying $\sum_{\ell=1}^\infty\ell\llvert\alpha_\ell\rrvert<
\infty$, the
result of Theorem~\ref{thesdmainfty}
reduces to the results obtained in Pfaffel and Schlemm~\cite{ps12} and
Yao~\cite{y12}.
\end{remark}

%re2.2 #&#
\begin{remark}
One can relax the assumption of simultaneous diagonalizability of the
coefficient matrices of the linear process to certain forms of
near-simulta\-ne\-ous diagonalizability, so that the conclusions of Theorem~\ref{thesdmainfty} continue to hold for linear processes where the
MA coefficients are
Toeplitz matrices whose entries decay away from the diagonal at an
appropriate rate.
Specifically, if $\mathbf{A}_\ell$ is the Toeplitz matrix with $j$th row
equaling $(a_{k-j,\ell}\dvtx0\leq k<p)$, for the bi-infinite
sequence $(a_{k,\ell})$
satisfying the condition
\[
\sum_{\ell=0}^\infty\sum
_{k=-\infty}^\infty\llvert k\rrvert^\beta(\ell+1)
\llvert a_{k,\ell}\rrvert< \infty%
\]
for some $\beta> 0$, which in particular implies Assumption~\ref{asmainfty}(a)
by the Gershgorin theorem, then the existence of the limiting ESD of symmetrized
autocovariance matrices can be proved.
%
%This aspect relates to the work of
%Hachem et~al. \cite{hln05}, who study the convergence of the empirical
%distribution of the sample covariance matrix of rectangular slices of
%bistationary Gaussian random fields. The data in their setting can be
%seen as
%realizations from a Gaussian linear process where the MA coefficients
%are
%Toeplitz matrices whose entries decay away from the diagonal. Since
%such
%Toeplitz matrices can be approximated by circulant matrices, which are
%diagonalized in the Fourier basis, \cite{hln05} do not need to assume
%Hermitianity (or symmetry).
%
For brevity, instead of giving a thorough technical argument, we only
provide the main idea of proof. First, the $\mathrm{MA}(\infty)$
series is approximated by an $\mathrm{MA}(q_p)$ series with $q_p=O(p^{1/3})$,
using arguments along the
line of Section~\ref{secintlin}. Second, banding with bandwidth
$k_p$ is applied to the coefficient matrices $\mathbf{A}_\ell$. It
can be shown through an application
of norm inequality, that the limiting spectral behavior is unchanged under
the banding so long as $k_p\to\infty$ under~(\ref{eqpn}). Third,
circulant matrices are constructed from the banded Toeplitz matrices by
periodization. The resulting matrices are therefore simultaneously
diagonalizable, and the eigenvalues of the $\ell$th approximate
coefficient matrix approximate the transfer function of the sequence
$(a_{k,\ell}\dvtx k\in\mathbb{Z})$. The limiting spectral behavior
is seen to be unchanged
after the use of the rank inequality so long as $k_p/p\to0$ under
(\ref{eqpn}).
The rest of the derivations follow the arguments
in the proof of Theorem~\ref{thesdmainfty}. While this particular result
is related to the work of Hachem et~al. \cite{hln05}, who study the
convergence of the empirical
distribution of the sample covariance matrix of rectangular slices of
bistationary Gaussian random fields, \cite{hln05}
does not cover the transition to non-Gaussian processes or the spectral
behavior of
sample autocovariance matrices.
\end{remark}

Several extensions of Theorem~\ref{thesdmainfty} are discussed in
Section~\ref{secex} below. The proof steps needed in order to verify
the main result
are outlined in Section~\ref{secstructure}, and the details are in
Sections~\ref{secgauss}--\ref{secpreal}. %Appendix~\ref{apptech}
The online SM \cite{lap-sm} contains additional technical lemmas.

%s3 #&#
\section{Examples and applications}\label{secexandapp}

%In this section, an ARMA(1,1) example is provided that serves as a way
%to both
%motivate the conditions imposed through Assumption~\ref{asmainfty}
%and illustrate the
%result of Theorem~\ref{thesdmainfty}. In the subsequent sections,
%existing results for the row-wise independent case are summarized, and
%potential applications
%of the proposed methodology are discussed.

%s3.1 #&#
\subsection{An $\mathrm{ARMA}(1,1)$ example}

In this section, let $(X_t\dvtx t\in\mathbb{Z})$ be the causal
ARMA($1,1$) process given by the stochastic difference equations
\[
\Phi(L)X_t=\Theta(L)Z_t, \qquad t\in\mathbb{Z},
\]
where $\Phi(L)=I-\Phi_1L$ and $\Theta(L)=I+\Theta_1L$ are,
respectively, the matrix-valued autoregressive and moving average
polynomials in the lag operator $L$ for which it is assumed that
$\llVert
\Phi_1\rrVert\leq\bar{\phi}<1$ and $\llVert\Theta_1\rrVert\leq\bar
\theta
<\infty$. Moreover, $(Z_t\dvtx t\in\mathbb{Z})\sim\mathrm{IID}(0,I)$ with entries possessing finite fourth moments. Under these
conditions $(X_t\dvtx t\in\mathbb{Z})$ admits the MA($\infty$)
representation
\[
X_t=\bA(L)Z_t,\qquad t\in\mathbb{Z},
\]
with $\bA(L)=\sum_{\ell=0}^\infty\bA_\ell L^\ell=\Phi
^{-1}(L)\Theta(L)$. Assume now further that $\Phi_1$ and $\Theta_1$
are simultaneously diagonalizable. Then $\Phi_1=\bU\bLambda_\Phi\bU
^*$ and $\Theta_1=\bU\bLambda_\Theta\bU^*$, where $\bLambda_\Phi
=\operatorname{diag}(\phi_1,\ldots,\phi_p)$ and $\bLambda_\Theta
=\operatorname{diag}(\theta_1,\ldots,\theta_p)$ such\vspace*{1pt} that $\max
_j\llvert\phi
_j\rrvert\leq\bar\phi$ and $\max_j\llvert\theta_j\rrvert\leq\bar
\theta$. With
regard to Assumption~\ref{asmainfty}, let $\lambda_j=(\phi_j,\theta
_j)^\prime\in\mathbb{R}^2$. Part (c) of the assumption then requires
almost sure weak convergence of the ESD of $\{\lambda_1,\ldots,\lambda
_p\}$ to a nonrandom probability distribution function on
$\mathbb{R}^2$. Moreover, using that for each coordinate,
\[
\frac{1+\theta_jL}{1-\phi_jL} =(1+\theta_jL)\sum_{\ell=0}^\infty(
\phi_jL)^\ell=1+(\theta_j+\phi_j)
\sum_{\ell=1}^\infty\phi_j^{\ell-1}L^\ell,
\]
it follows that $\bA_\ell=\bU\operatorname{diag}(f_\ell(\lambda
_1),\ldots,f_\ell(\lambda_p))\bU^*$ with $f_0(\lambda_j)=1$ and
$f_\ell(\lambda_j)=(\theta_j+\phi_j)\phi_j^{\ell-1}$ for $\ell
\in\mathbb{N}$. This illustrates part (b) of Assumption~\ref{asmainfty}. The summability conditions stated in part (a) are clearly
satisfied. Generalization to arbitrary causal ARMA models follows in a
similar fashion.

%s3.2 #&#
\subsection{Time series with independent rows}\label{secind}

In this section, the situation of time series with independent rows is
considered. Our results describe the
limiting ESD of the symmetrized sample autocovariances in the setting
where the $j$th row of the time
series, denoted $\{\xi_{jt}\dvtx t \in\mathbb{Z}\}$, is given by
%
%e3.1 #&#
\begin{equation}
\label{eqidiot} \xi_{jt}=\sum_{\ell=0}^\infty
f_\ell(\lambda_j)Z_{j,t-\ell},
\end{equation}
where: (i) the $Z_{jt}$'s are independent,
identically distributed real- or complex-valued random variables
with mean zero, unit variance and finite fourth moments; (ii) the
$f_\ell$'s are continuous functions
from $\mathbb{R}^m \to\mathbb{R}$ satisfying $f_0(\lambda) \equiv
1$ and the summability condition $\sup_{\lambda\in\operatorname{supp}(F^{\mathbf{A}})}
\sum_{\ell=0}^\infty(\ell+1)\llvert f_\ell(\lambda)\rrvert< \infty$; (iii)
the $\lambda_j$'s are i.i.d. realizations
from an $m$-dimensional probability distribution denoted by $F^\mathbf
{A}$. If the supremum in condition (ii) is taken over $\mathbb{R}^m$,
condition (iii) can be weakened to require that the empirical
distribution of the $\lambda_j$'s converges almost surely to a
nonrandom distribution $F^\mathbf{A}$.
%%such that $\left\Vert f_\ell(\lambda)\right\Vert\leq\bar\lambda_

Let $h(\lambda,\nu) = \llvert\sum_{\ell=0}^\infty e^{i\ell\nu}f_\ell
(\lambda)\rrvert^2$.
Then the empirical distributions of the eigenvalues of the lag-$\tau$
symmetrized autocovariance
matrices converge almost surely to a nonrandom probability distribution
$F_\tau$
with Stieltjes transform $s_\tau$ determined by equations (\ref
{eqstau}) and (\ref{eqktau}),
where $K_\tau(z,\nu)$ is as in Theorem~\ref{thesdmainfty}. This
is in the spirit
of the works of
%Hachem et~al. \cite{hln05},
Jin et~al. \cite{jwml09}, Pfaffel and Schlemm \cite{ps12}
and Yao \cite{y12}, who studied the sample covariance case with
$f_\ell(\lambda)\equiv\bar{f}_\ell$ for all $\ell\geq0$,
Hachem et~al. \cite{hln05}, who considered the sample covariance case
for stationary Gaussian fields
and Jin et~al. \cite{jwbnh14}, who studied the symmetrized sample
autocovariance case with $f_1(\lambda)\equiv1$
and $f_\ell(\lambda)\equiv0$ for $\ell\geq1$ (i.e., when the $\xi
_{jt}$'s are i.i.d. with zero mean and unit variance).

%s3.3 #&#
\subsection{Signal processing and diagnostic checks}

The results derived here can be useful in dealing with a number of important
statistical questions. Signal detection in a noisy background is one of the
most important problems in signal processing and communications theory. Often
the observations are taken in time, and the standard assumption is that the
noise is i.i.d. in time, referred to as white noise. However, in
spatio-temporal signal processing, it is quite apt to formulate the
noise as
``colored'' or correlated in time, as well as in the spatial dimension. The
proposed model for the time series is a good prototype for such a noise
structure. Thus the problem of detecting a low-dimensional signal
embedded in
\mbox{high-}dimensional noise, for example, through a factor model framework,
can be
effectively addressed by making use of the behavior of the ESDs of
autocovariances of the noise. Another potential application of the
results is
in building diagnostic tools for high-dimensional time series. By
focusing on
the ESDs of the autocovariances for various lag orders, or that of a tapered
estimate of the spectral density operator, one can infer about the
nature of
dependence, provided the model assumptions hold.
%The latter is often a reasonable assumption for modeling panel data
%commonly arising in
%econometrics.
The proposed model also provides a broad class of alternatives for the
hypothesis of independence of observations in settings where those observations
are measured in time. Finally, in practical applications, it is of interest
whether the spectrum of the coefficient matrices of the linear process
can be
estimated from the data. The equations for the limiting Stieltjes
kernel and
its relation to the Stieltjes transform of the autocovariance matrices provide
a tool for attacking this problem. This aspect has been explored in the
Ph.D. thesis of the first author \cite{l13} and the methodology will be reported
elsewhere.

%s3.4 #&#
\subsection{Dynamic factor models}\label{secexandappdfm}

Forni and Lippi \cite{fl99} describe a class of time series models
that captures
the subject specific variations of microeconomic activities. This class of
models, referred to as Dynamic Factor Models (DFM), has proved
immensely popular in the
econometrics community and beyond. DFMs have, for example, been used for
describing the stock returns in \cite{ner92}, forecasting national
accounts in \cite{abr10},
modeling portfolio allocation in \cite{aw10} and modeling
psychological development in \cite{m94},
as well as in many other applications. Important theoretical and
inferential questions regarding DFMs
have been investigated in a series of papers by Forni and Lippi \cite
{fl01}, Forni et~al. \cite{fhlr00,fhlr04,fhlr05} and Stock and Watson
\cite{sw05}, to name
a few.
DFMs have also shown early promise for applications to other
interesting multivariate time series
problems such as the study of fMRI data.

A DFM can be described as follows. As in \cite{fl99}, let $Y_{jt}$ be
the response corresponding to the $j$th
individual/agent at time $t$, modeled as
%
%e3.2 #&#
\begin{equation}
\label{eqDFM} Y_{jt} = b_{j1}(L)U_{1t} + \cdots+
b_{jM}(L)U_{Mt} + \xi_{jt}, \qquad j=1,\ldots,p.
\end{equation}
The model specifies that $Y_{jt}$ is determined by a
small, fixed number of underlying common factors $U_{kt}$ and their lags,
determined by the polynomials $b_{jk}(L)$ in the lag-operator $L$,
plus an idiosyncratic component $\xi_{jt}$ assumed independent across
individuals. Typically, $(\xi_{jt}\dvtx t\in\mathbb{Z})$ is taken
to be a stationary linear processes, independent across $j$.

One of the key questions pertaining to DFM is the determination of the
number of
dynamic factors. This question has been investigated by Bai and Ng
\cite{bn07}, Stock and Watson \cite{sw05} and Hallin and Li\v{s}ka
\cite{hl07}. Unlike in PCA, here one has to deal with the additional
problem of detecting the lag orders of the dynamic factors. This can be
approached through the study of the behavior of the extreme eigenvalues
of the sample autocovariance matrices as in Jin et~al. \cite
{jwbnh14}. The issue becomes even more challenging
when the dimensionality of the problem increases. In such settings, one expects
that a form of phase transition phenomenon, well known in the context of
a high-dimensional static factor model (or spiked covariance model)
with i.i.d. observations
(see, e.g., Baik and Silverstein \cite{bs06}), will set in. In
particular, as Jin et~al. \cite{jwbnh14} argue,
a dynamic factor will be detectable from the data only if the
corresponding total
signal intensity, as measured, for example, by the sum of the variances
of the factor
loadings, is above a threshold. Moreover, the number of eigenvalues
that lie outside
the bulk of the eigenvalues of the symmetrized sample autocovariance of
a certain
lag order provide information about the lag order of the DFM.
Driven by the analogy with the static factor model with
i.i.d. observations, it is expected that the detection thresholds will
depend on the dimension-to-sample
size ratio, as well as the behavior of the bulk spectrum of the autocovariances
of the idiosyncratic terms at specific lag orders, including the
support of the
limiting ESD. Equation (\ref{eqidiot}) in Section~\ref{secind} constitutes
the ``null'' model for the DFM in which the dynamic common factors are absent.
Therefore follow-up studies on the different aspects of the ESD of the
symmetrized
sample autocovariances of such processes will be helpful in determining
the detection
thresholds and estimation characteristics of high-dimensional DFMs.

%s3.5 #&#
\subsection{An idealized production model}
\label{secexandappproduction}

Onatski \cite{o12} describes a model for production $Y_t$, at time
$t$, involving $p$ different industries
in an economy that is given by the equations
%
%e3.3 #&#
\begin{equation}
\label{eqproductionmodel} Y_t = \mathbf{W}^{-1} \Biggl(\sum
_{k=1}^M u_{kt} \mathbf{f}_k
+ \xi_t \Biggr).
\end{equation}
Model~(\ref{eqproductionmodel}) is a static factor model in which
$p\times1$ vectors $\mathbf{f}_k$ denote the (unobserved) common
static factors, $u_{kt}$ denote the (unobserved) factor scores
consisting of independent time series corresponding to different
factors $k$ and $\xi_t$ denote the $p\times1$ vectors of
idiosyncratic components. The entries of the matrix $\mathbf{W}$
indicate the interactions among the different industries.
In the following an enhanced version of the model is considered where
the economy is thought to be
divided into a finite number of distinct sectors for which the interaction
across the sectors is assumed ``weak'' in a suitable sense to be described.
In addition, the assumption of separable covariance structure of the
$\xi_t$ made in \cite{o12} is relaxed by requiring instead that the
temporal variation in $\xi_t$
for all the industrial units within a sector is the same and is
stationary in time.
This assumption means that the component of the vector $\xi_t$
corresponding to a particular
sector has a separable covariance structure with stationary time
variation, and the
components corresponding to different sectors are independent.
Specifically, if there
are $K$ sectors, we can divide $\mathbf{W}$ into $K\times K$ block matrices
\[
\mathbf{W} = \lleft[\matrix{ \mathbf{W}_{11} & \mathbf{W}_{12} & \cdots& \mathbf{W}_{1K}
\cr
\cdot& \cdot& \cdots& \cdot
\cr
\mathbf{W}_{K1} & \mathbf{W}_{K2} & \cdots&
\mathbf{W}_{KK}} \rright] = \widetilde{\mathbf{W}} + \Delta,
\]
where $\widetilde{\mathbf{W}} = \operatorname{diag}(\mathbf{W}_{11},\ldots
,\mathbf{W}_{KK})$.
If the sectors have no interaction at all, that is, if $\mathbf
{W}_{jk} = 0$ for all $j\neq k$,
then the corresponding data model is an instance of a blockwise
separable covariance model.
``Weak interaction'' means that the norms of the off-diagonal blocks
in the matrix $\mathbf{W}$ are small. More precisely,
if $p^{-1} \llVert\Delta\rrVert_F^2 \to0$ as $p \to\infty$, then the
limiting ESD of the
symmetrized autocovariances for the data matrix $\mathbf{Y}$ is the
same as that of
$\widetilde{\mathbf{Y}}$ obtained by replacing $\mathbf{W}$ by
$\widetilde{\mathbf{W}}$ in
(\ref{eqproductionmodel}). Under the assumption of a linear process
structure on the different\vspace*{1pt} components of $\xi_t$, and a natural
requirement on the stability
of the singular values of~$\widetilde{\mathbf{W}}$, the existence and
characterization of the
limiting ESDs of the symmetrized autocovariances of $\widetilde
{\mathbf{Y}}$ can be dealt within
the framework studied in Section~\ref{secmain}. These limiting
ESDs will help in determining the detection thresholds for the static factors,
or even dynamic factors, if the model were to be enhanced further.

%s4 #&#
\section{Extensions of the main results}\label{secex}

This section discusses three different extensions of the
main result. The arguments for the proof are similar to that of the
proof of
Theorem~\ref{thesdmainfty} and hence only a brief outline is provided.
Moreover, the results stated here apply to both real- and complex-valued
cases, the only difference being that, in the former case, the relevant matrices
are real symmetric while, in the latter case, they are Hermitian.

The first extension involves a rescaling of the process defined
in (\ref{eqmainfty}). Thus it is assumed that
%
%e4.1 #&#
\begin{equation}
\label{eqXtscaled} X_t = \mathbf{B}^{1/2}\sum
_{\ell=0}^\infty\mathbf{A}_\ell Z_{t-\ell},
\qquad t\in\mathbb{Z},
\end{equation}
where the processes $(Z_t\dvtx t\in\mathbb{Z})$ and matrices
$(A_\ell\dvtx
\ell\in\mathbb{N}_0)$ satisfy Assumption~\ref{asmainfty}, and the matrix
$\mathbf{B}^{1/2}$ is the square root of a $p\times p$ positive semidefinite
Hermitian (or symmetric) matrix $\mathbf{B}$ satisfying the following
assumption.

%as4.1 #&#
\begin{assumption}\label{asscaling}
Let\/ $\mathbf{U}$ be as in Assumption~\ref{asmainfty}. Then
\[
\mathbf{U}^*\mathbf{B}\mathbf{U} = \operatorname{diag}\bigl(g_{\mathbf{B}}(
\lambda_1),\ldots,g_{\mathbf
{B}}(\lambda_p)\bigr),
\]
where\vspace*{1pt} $g_{\mathbf{B}}\dvtx\mathbb{R}^m \to\mathbb{R}_+$ is
continuous and bounded
on $\mathbb{R}^m$, and $\{\lambda\dvtx g_{\mathbf{B}}(\lambda) >
0\}\cap
\operatorname{supp}(F^{\mathbf{A}})$ is nonempty.
\end{assumption}

As before, $\hat F_\tau$ is defined to be the ESD of the symmetrized
autocovariance $\mathbf{C}_\tau$ of lag order $\tau$. If the linear process
$(X_t\dvtx t\in\mathbb{Z})$ defined through (\ref{eqXtscaled}) satisfies
Assumptions~\ref{asmainfty} and~\ref{asscaling}, then the statement of
Theorem~\ref{thesdmainfty} holds with the function $h(\lambda,\nu
)$ replaced
by $g_{\mathbf{B}}(\lambda) h(\lambda,\nu)$.

%t\in\mathbb{Z})$ defined through (\ref{eqXtscaled}) satisfies
%Assumptions~\ref{asmainfty} and~\ref{asscaling}, then, with
%probability one
%and in the high-dimensional setting \eqref{eqpn}, $\hat F_\tau$
%converges
%weakly to a probability distribution $F_\tau$ with Stieltjes transform
%$s_\tau$
%determined by the equation
%s_\tau(z)=\int\bigg[\frac{1}{2\pi}\int_0^{2\pi}\frac{\cos(\tau\nu)g_{
%K_\tau(z,\nu)}\,d\nu-z\bigg]^{-1} \,dF^\bA(\lambda),
%where $K_\tau\dvtx\mathbb{C}^+\times[0,2\pi]\to\mathbb{C}^+$ is a
%Stieltjes
%kernel, that is, $K_\tau(\cdot,\nu)$ is the Stieltjes transform of a
%measure
%with total mass $\int g_{\mathbf{B}}(\lambda) h(\lambda,\nu)\,dF^\bA(
%for every fixed $\nu\in[0,2\pi]$. Moreover, $K_\tau$ is the unique
%solution to
%h(\lambda,\nu^\prime)}{1+c\cos(\tau\nu^\prime)K_\tau(z,\nu^\prime)}\,d
%subject to the restriction that $K_\tau$ is a Stieltjes kernel.

The second extension is about the existence and description of the
limiting ESD of the autocovariances of linear filters of the process
$(X_t\dvtx t\in\mathbb{Z})$ defined through~(\ref{eqXtscaled}).
A linear
filter of this process is of the form
%
%e4.2 #&#
\begin{equation}
Y_t = \sum_{k = 0}^\infty
b_k X_{t-k}, \qquad t \in\mathbb{Z},
\end{equation}
where $(b_k\dvtx k\in\mathbb{N}_0)$ is a sequence of real numbers
for which the following
summability condition is needed.

%as4.2 #&#
\begin{assumption}
The sequence $(b_k\dvtx k\geq0)$ satisfies $\sum_{k=0}^\infty
k \llvert b_k\rrvert< \infty$.
\end{assumption}

If the linear process $(X_t\dvtx t\in\mathbb{Z})$ defined through
(\ref{eqXtscaled}) satisfies Assumptions~\ref{asmainfty} and
\ref{asscaling}, then the statement of Theorem~\ref{thesdmainfty} holds
with the function $h(\lambda,\nu)$ replaced by $\zeta(\nu)
g_{\mathbf{B}}(\lambda) h(\lambda,\nu)$, where $\zeta(\nu) =
\llvert\sum_{k=0}^\infty
e^{ik\nu}b_k\rrvert^2$, $\nu\in[0,2\pi]$. This result follows using
the properties
of convolution and Fourier transform.

The third extension is about estimation of the
spectral density operator
%
%e4.3 #&#
\begin{equation}
\Gamma(\eta) = \sum_{\tau=-\infty}^\infty
e^{i\tau\eta} \mathbb{E}\bigl[X_t X_{t+\tau}^*\bigr],
\qquad\eta\in[0,2\pi].
\end{equation}
It is well known from classical multivariate time series analysis (see,
e.g., Chapter~10 of Hamilton~\cite{h94}) that the ``natural''
estimator that
replaces the population autocovariance $\mathbb{E}[X_t X_{t+\tau}^*]$
by the
corresponding sample autocovariance may not be positive definite. In
order to
obtain positive definite estimators and a better bias-variance
trade-off, it is
therefore standard in the literature to consider certain tapered estimators
with standard choices given, for example, by the Bartlett and Parzen
kernels as
described in \cite{h94}. In the following, the behavior of a class of tapered
estimators of $\Gamma(\eta)$, which are given by
%
%e4.4 #&#
\begin{equation}
\label{eqtaperedspectraldensity} \widehat\Gamma_n(\eta) = \sum
_{\tau= -\infty}^\infty T_n(\tau) e^{i\tau\eta}
\frac{1}{n} \sum_{t=1}^{n-\tau}
X_t X_{t+\tau}^*, \qquad\eta\in[0,2\pi],
\end{equation}
is studied, where $T_n(\cdot)$ is a sequence of even functions and the
quantities
$T_n(\tau)$ are known as tapering weights for which the following restriction
is imposed.

%as4.3 #&#
\begin{assumption}\label{astapering}
(i) The even functions $T_n(\cdot)$ are such that $T_n(\tau) = 0$ for
$\tau\geq
n$; (ii) there exists an even function $T(x)$ such that $T_n(x) \to
T(x)$ as $n
\to\infty$ and $\llvert T_n(x) \rrvert\leq C \llvert T(x) \rrvert
$ for some $C > 0$, for all
$x$; (iii)
$\sum_{\tau=0}^\infty(1+\llvert\tau\rrvert) \llvert T(\tau)
\rrvert< \infty$.
\end{assumption}

An implication of this assumption is that the function $\mathfrak
{f}_T(\eta)$
defined by
%
%e4.5 #&#
\begin{equation}
\mathfrak{f}_T(\eta) = 1 + 2 \sum_{\tau=1}^\infty
T(\tau) \cos(\tau\eta)
\end{equation}
is well defined and is uniformly Lipschitz, and %the sequence
$\mathfrak{f}_{T_n}(\eta) = 1+\sum_{\tau=1}^\infty T_n(\tau) \cos
(\tau\eta)$
converges to $\mathfrak{f}_T(\eta)$ uniformly in $\eta$. Examples\vspace*{2pt} of kernels
$T_n$ are $T_n(x) = (1+\llvert x\rrvert)^{-\alpha}\mathbb
{I}(\llvert x\rrvert< n)$ for $\alpha
> 2$ and
$T_n(x) = \beta^{\llvert x\rrvert}\mathbb{I}(\llvert x\rrvert<
n)$ for $\beta\in(0,1)$. It
can be seen
from Assumption~\ref{astapering} that in the high-dimensional setting under
consideration here, standard choices for tapering weights, such as those
given by the Bartlett and Parzen kernels, are ruled out.
Now the following generalization of
Theorem~\ref{thesdmainfty} is obtained.

%th4.1 #&#
\begin{theorem}
\label{thesdmainftyspectraldensity} Suppose that the linear process
$(X_t\dvtx t\in\mathbb{Z})$ defined through (\ref{eqXtscaled}) satisfies
Assumptions~\ref{asmainfty} and~\ref{asscaling}, and that the
estimated spectral
density operators $(\widehat\Gamma_n(\eta)\dvtx\eta\in[0,2\pi
])$ are
defined by (\ref{eqtaperedspectraldensity}) with tapering weights
satisfying Assumption~\ref{astapering}. Then, with probability one
and in the
high-\break dimensional setting~(\ref{eqpn}), for every $\eta\in[0,2\pi
]$, the ESD
of $\widehat\Gamma_n(\eta)$ converges weakly to a probability distribution
$F_{T,\eta}$ with Stieltjes transform $s_{T,\eta}$ determined by the equation
\[
s_{T,\eta}(z)=\int\biggl[\frac{1}{2\pi}\int_0^{2\pi}
\frac
{\mathfrak{f}_T(\eta-
\nu)g_{\mathbf{B}}(\lambda)h(\lambda,\nu)}{1+c\mathfrak{f}_T(\eta
- \nu)
K_{T,\eta}(z,\nu)}\,d\nu-z \biggr]^{-1} \,dF^\bA(\lambda),
\]
where $K_{T,\eta}\dvtx\mathbb{C}^+\times[0,2\pi]\to\mathbb
{C}^+$ is a
Stieltjes kernel; that is, $K_{T,\eta}(\cdot,\nu)$ is the Stieltjes transform
of a measure with total mass $m_\nu=\int g_{\mathbf{B}}(\lambda)
h(\lambda,\nu)\,dF^\bA(\lambda)$ for every fixed $\nu\in[0,2\pi]$,
whenever $m_\nu>0$. Moreover,
$K_{T,\eta}$ is the unique solution to
%
%e4.6 #&#
\begin{eqnarray*}%\label{eqkTeta}
&& K_{T,\eta}(z,\nu)\nonumber
\\
&&\qquad =\int\biggl[\frac{1}{2\pi}\int_0^{2\pi}
\frac{\mathfrak{f}_T(\eta- \nu^\prime)g_{\mathbf{B}}(\lambda)
h(\lambda,\nu^\prime)}{1+c\mathfrak{f}_T(\eta- \nu^\prime
)K_{T,\eta}(z,\nu^\prime)}\,d\nu^\prime-z \biggr]^{-1}
\\
&&\quad\qquad{}\hspace*{7pt}\times g_{\mathbf{B}}(\lambda)h(\lambda,\nu) \,dF^\bA(\lambda),\nonumber
\end{eqnarray*}
subject to the restriction that $K_{T,\eta}$ is a Stieltjes kernel.
Else, if
$m_\nu=0$, $K_{T,\eta}(z,\nu)$ is identically zero on $\mathbb
{C}^+$ and so
still satisfies the latter equation. %\eqref{eqkTeta}.
\end{theorem}

Section~\ref{secspecdens} outlines the main argument needed %in order
to prove Theorem~\ref{thesdmainftyspectraldensity}.

%s5 #&#
\section{Simulations}\label{secsimulations}

%f1 #&#
\begin{figure}%[t]

\includegraphics{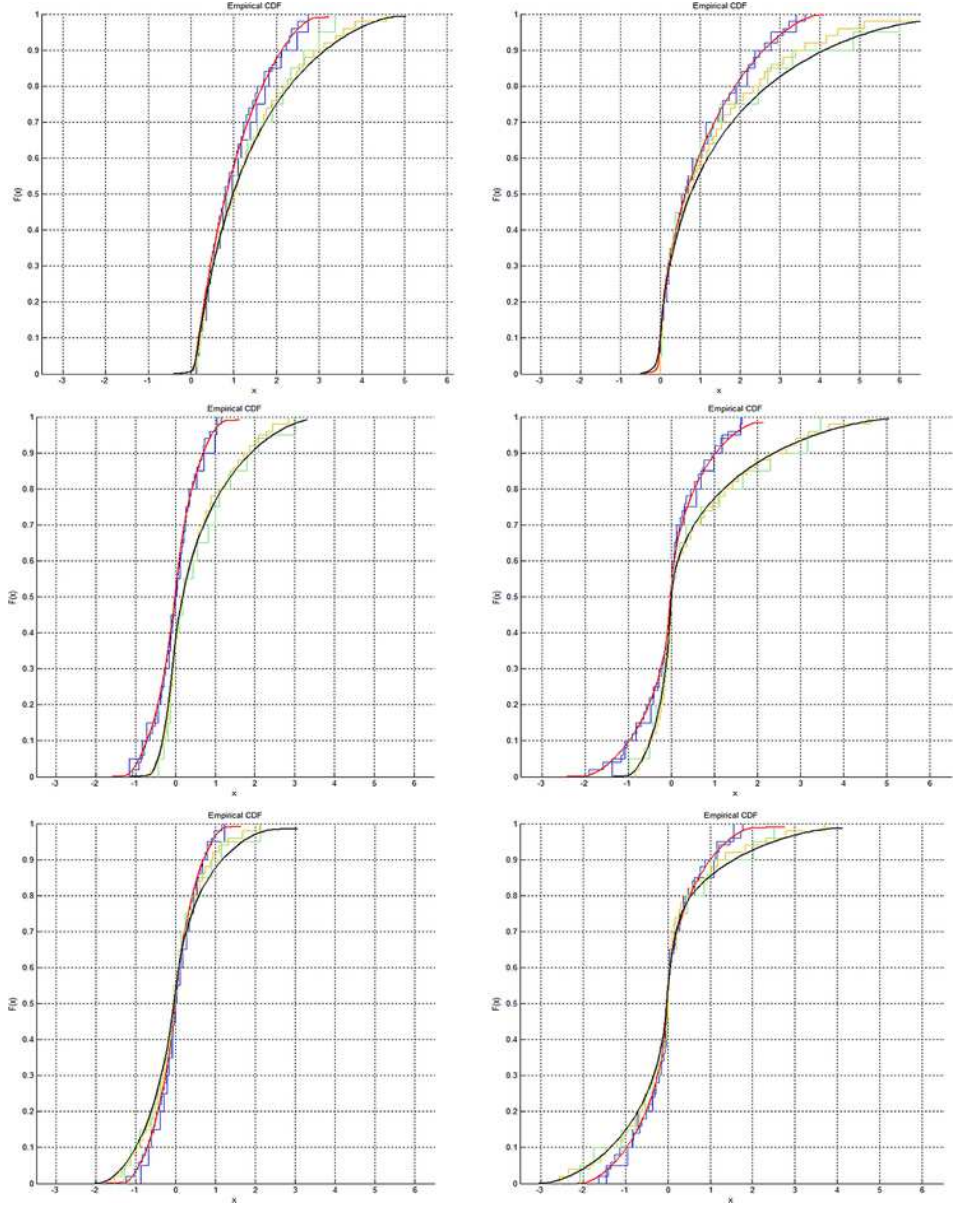}

\caption{ESDs and LSDs of the sample autocovariances for i.i.d. and $\mathrm{MA}(1)$
observations. {Left} panel: $p/n = 0.5$; {right} panel:
$p/n = 1$. {Top} panel: $\tau= 0$; {middle} panel: $\tau= 1$;
{bottom} panel: $\tau= 2$. {Red} curve: LSD for i.i.d.;
{black} curve: LSD for $\mathrm{MA}(1)$. ESDs are corresponding to
$p=20$ [deep blue for i.i.d., light green for $\mathrm{MA}(1)$] and $p=50$ [light blue for i.i.d.,
orange for $\mathrm{MA}(1)$].}\label{figsimul2050}
\end{figure}

In this section, a small simulation study is conducted to illustrate the
behavior of the LSD (limiting ESD) of the symmetrized sample autocovariances
$\bC_\tau$ for different lag orders $\tau$ when the observations are
i.i.d. ($X_t = Z_t$) versus when they come from an $\mathrm{MA}(1)$ process
($X_t = Z_t +
\mathbf{A}_1 Z_{t-1}$) with a symmetric coefficient matrix $\mathbf
{A}_1$. In
this study, two different sets of $(n,p)$ are chosen such that $p/n =
0.5$ and
$p/n = 1$, respectively, for $p=20$ and 50. Also, for the $\mathrm{MA}(1)$ case,
the ESD
of $\mathbf{A}_1$ is chosen to be $0.5 \delta_{0.2} + 0.5 \delta
_{0.8}$, where
$\delta_b$ denotes the degenerate distribution at $b$. The
distribution of the
$Z_{jt}$'s is chosen to be i.i.d. $N(0,1)$. In Figure~\ref{figsimul2050},
the ESD of $\bC_\tau$ is plotted for the lags $\tau=0,1,2$,
for one random realization corresponding to each setting. The theoretical
limits (c.d.f.) for all cases are plotted as solid smooth curves [red for
i.i.d. and black for $\mathrm{MA}(1)$]. These c.d.f.'s are obtained through numerically
inverting the Stieltjes transform $s_\tau$ of the LSD of $\bC_\tau$
using the
inversion formula for Stieltjes transforms; cf. \cite{bs10}, \cite
{pa13}. Here
$s_\tau$ is obtained from equations (\ref{eqstau}) and (\ref
{eqktau}) where
$c=p/n$ and $h(\lambda,\nu) = 1+2\cos(\nu) \lambda+ \lambda^2$.
Given $F^{\bA}
= 0.5 \delta_{0.2} + 0.5 \delta_{0.8}$, the Stieltjes kernel $K_\tau
(z,\nu)$
for the $\mathrm{MA}(1)$ case is solved numerically
%by writing it as $a_0(z) + 2\cos\nu a_1(z) + a_2(z)$ for appropriate
%functions $a_j(z)$'s and
by using calculus of residues. The computational algorithm can be found
in the
first author's Ph.D. thesis \cite{l13}.

The graphs clearly show the distinction in behavior of ESD of
symmetrized sample autocovariances between the i.i.d. and the $\mathrm{MA}(1)$
processes. The LSDs of $\bC_\tau$ for the i.i.d. case for
$\tau= 1$ and 2 are the same, which follows immediately from equations
(\ref{eqstau}) and (\ref{eqktau}) that reduce to a single equation
since the
i.i.d. case corresponds to $F^{\bA} = \delta_0$. The behavior of the
LSDs in
the $\mathrm{MA}(1)$ case is distinctly different. It can be shown that LSDs of
$\bC_\tau$ in the $\mathrm{MA}(1)$ case converge to the LSD for $\bC_1$ in
the i.i.d. case
as $\tau$ becomes larger. Owing to space constraints, the graphical
displays for
higher lags are omitted. Another important feature is that the LSDs approximate
the ESDs quite well even for $p$ as small as 20, indicating a fast convergence.

%s6 #&#
\section{The structure of the proof}\label{secstructure}

%s6.1 #&#
\subsection{Developing intuition for the Gaussian $\mathrm{MA}(1)$ case}\label{subsecintu}

\begin{Outline*}
In this section, the overall proof strategy is briefly outlined,
and the intuition behind the individual steps is developed for simpler
first-order moving average, $\mathrm{MA}(1)$, time series. The key ideas in the
proof of Theorem~\ref{thesdmainfty} consist of showing that:
\begin{itemize}
\item there is a unique Stieltjes kernel solution $K_\tau$ to equation
(\ref{eqktau});
\item almost surely, the Stieltjes transform of $\hat F_\tau$, say,
$s_{p,\tau}$ converges pointwise to a Stieltjes transform $s_\tau$ which will be identified
with $F_\tau$;
\item$\hat F_\tau$ is tight.
\end{itemize}
To achieve the second item, one can argue as follows. First, replace the
original linear process observations $X_1,\ldots,X_n$ with transformed vectors
$\tX_1,\ldots,\tX_n$ that are serially independent. Second,\vspace*{2pt} replace the
symmetrized lag-$\tau$ autocovariance matrix $\bC_\tau$ by a transformed
version $\tbC_\tau$ built from $\tX_1,\ldots,\tX_n$. A heuristic
formulation
for the simpler Gaussian $\mathrm{MA}(1)$ case is given below in some detail. Once these
two steps have been achieved, the proof proceeds by verifying some technical
conditions with the help of classical RMT results, available, for
example, in the
monograph Bai and Silverstein~\cite{bs10}. In the following, let
$(X_t\dvtx
t\in\mathbb{Z})$ denote the $p$-dimensional $\mathrm{MA}(1)$ process given by the
equations
%
%e6.1 #&#
\begin{equation}
\label{eqma1} X_t=Z_t+\bA_1Z_{t-1},
\qquad t\in\mathbb{Z},
\end{equation}
where $(Z_t\dvtx t\in\mathbb{Z})$ are assumed complex Gaussian in
addition to the requirements of Section~\ref{secmain}. For this time
series, the
conditions imposed through Assumption~\ref{asmainfty} simplify
considerably with
part (a) reducing to the condition that the eigenvalues of $\bA_1$ be
uniformly bounded and part (b) being\vspace*{1pt} satisfied by choosing $f_1$ as the
identity.
Moreover, $\bA=[\bI\dvtx\bA_1]$, $\psi(\lambda,\nu)=1+e^{i\nu
}\lambda$,
$\psi(\bA,\nu)=\bI+e^{i\nu}\bA_1$, $h(\lambda,\nu)=1+2\cos(\nu
)\lambda+\lambda^2$
and $\cH(\bA,\nu)=\bI+2\cos(\nu)\bA_1+\bA_1^2$, implying for
each $\lambda$,
$h(\lambda,\nu)$ is the spectral density (up to normalization) of a univariate
$\mathrm{MA}(1)$ process with parameter $\lambda$.
\end{Outline*}

\begin{TSI*}%{\sc Transformation to serial independence.}
The transformation to serial
independence requires two steps, the first consisting of an
approximation of
the lag operator by a circulant matrix and the second of a rotation
using the
complex Fourier basis to achieve independence. %The details are as
%follows.
Accordingly, let
\[
\bL=[o\dvtx e_1\dvtx\cdots\dvtx e_{n-1}] \quad\mbox{and}\quad\tbL
=[e_n\dvtx e_1\dvtx\cdots\dvtx e_{n-1}]
\]
be the $n \times n$ lag operator and its approximating circulant matrix,
respectively, where $o$ denotes the $n$-dimensional zero vector and
$e_j$ the
$j$th canonical unit vector in $\mathbb{R}^n$ taking the value 1 in
the $j$th
component and 0 elsewhere. Since $\tbL$ is a circulant matrix, its spectral
decomposition is $\tbL u_t= \eta_t u_t$, $t=1,\ldots,n$,
where $\eta_t=e^{i\nu_t}$, $\nu_t={2\pi t}/{n}$ and $u_t$ the vector whose
$j$th entry is $\eta_t^j$. It follows that $\tbL$ diagonalizes in the complex
Fourier basis with the usual Fourier frequencies. Let
$\bolds{\Lambda}_{\tbL}=\operatorname{diag}(\eta_1,\ldots,\eta_n)$ and
$\bU_{\tbL}=[u_1\dvtx\cdots\dvtx u_n]$ denote the corresponding
eigenvalue and
eigenvector matrices, so that
$\tbL=\bU_{\tbL}\bolds{\Lambda}_{\tbL}\bU_{\tbL}^*$. Using
$\bX=[X_1\dvtx\cdots\dvtx X_n]$ and $\bZ=[Z_1\dvtx\cdots\dvtx Z_n]$,
the $\mathrm{MA}(1)$ process
(\ref{eqma1}) can be transformed into
\[
\bX_1=\bZ+\bA_1\bZ\tbL,
\]
where $\bX_1$ constitutes a redefinition of $\bX$ such that the first
column is
changed to $Z_1+\bA_1Z_n$, while all other columns are as in the
original data
matrix $\bX$. Rotating in the complex Fourier basis, the observations are
transformed again into the vectors $\tX_1,\ldots,\tX_n$ given by
\[
\tbX=[\tX_1\dvtx\cdots\dvtx\tX_n]=\bX_1
\bU_{\tbL}. %(\bZ+\bA_1\bZ\tbL)
\]
Observe that $\tbX$ has independent columns. To see this, note first that
$\tbZ=\bZ\bU_{\tbL}$ possesses the same distribution as $\bZ$,
since $\bZ$ has
(complex) Gaussian entries and $\bU_{\tbL}$ is a unitary matrix.
Write then
\[
\tbX=\tbZ+\bA_1\tbZ\bolds{\Lambda}_{\tbL} = \bigl[(\bI+
\eta_1\bA_1)\tZ_1\dvtx\cdots\dvtx(\bI+
\eta_n\bA_1)\tZ_n\bigr],
\]
where $\tbZ=[\tZ_1\dvtx\cdots\dvtx\tZ_n]$, and independence of the columns
(and thus serial
independence) follows. Note also that $\bI+\eta_t\bA_1=\psi(\bA,\nu_t)$
and consequently
$\tbX=\psi(\bA,\nu_t)\tbZ$, using the transfer function $\psi$
defined in (\ref{eqtransferpsi}).
%introduced above Theorem~\ref{thesdmainfty}.
\end{TSI*}

\begin{TT*}%{\sc Transformation of $\bC_\tau$.}
The vectors $\tX_1,\ldots,\tX
_n$ give rise
to approximations $\tbC_\tau$ to the lag-$\tau$ symmetrized autocovariance
matrices $\bC_\tau$; in particular $\tbC_\tau$ and $\bC_\tau$
will be shown to
have the same large-sample spectral behavior, irrespective of
the distribution of the entries. So let
\[
\tbC_\tau\equiv\frac{1}{2n} \bX_1 \bigl(
\tbL^\tau+\bigl[\tbL^*\bigr]^\tau\bigr)\bX_1^* =
\frac{1}n \sum_{t=1}^n\cos(\tau
\nu_t)\tX_t\tX_t^*,
\]
where the latter equality follows from several small computations using the
quantities introduced in the preceding paragraph. Now,
\begin{eqnarray*}
\bC_\tau-\tbC_\tau&=& \biggl(\bC_\tau-
\frac{1}{2n}\bX_1 \bigl(\bL^\tau+\bigl[\bL^*
\bigr]^\tau\bigr)\bX_1^* \biggr)
\\
&&{} + \biggl(\frac
{1}{2n}\bX_1 \bigl(\bL^\tau+\bigl[
\bL^*\bigr]^\tau\bigr)\bX_1^*-\frac
{1}{2n}
\bX_1 \bigl(\tbL^\tau+\bigl[\tbL^*\bigr]^\tau
\bigr)\bX_1^* \biggr).
\end{eqnarray*}
The rank of the first difference on the right-hand side of the last
display is
at most~$2$, since $\bX$ and $\bX_1$ differ only in the first column.
The rank
of the second difference is at most $2\tau$. The rank of $\bC_\tau
-\tbC_\tau$
is therefore at most $2(1+\tau)$. Defining the resolvents
\[
\bR_\tau(z)=(\bC_\tau-z\bI)^{-1} \quad\mbox{and}\quad\tbR_\tau(z)=(\tbC_\tau-z\bI)^{-1},
\]
the Stieltjes transforms corresponding to the ESDs of $\bC_\tau$ and
$\tbC_\tau$ are, respectively, given by
\[
s_{\tau,p}(z)=\frac{1}p\operatorname{tr}\bigl[\bR_\tau(z)
\bigr] \quad\mbox{and}\quad\tilde s_{\tau,p}(z)=\frac{1}p
\operatorname{tr}\bigl[\tbR_\tau(z)\bigr].
\]
It follows from
%Lemma 2.6 of Silverstein and Bai~\cite{sb95} that,
Lemma~S.2 %\ref{lemrankinequality}
that with probability one, the ESDs of
$\bC_\tau$ and $\tbC_\tau$ converge to the same limit, provided the limits
exist. We conclude from Lemma~S.3 %\ref{lemStieltjesconvergence}
that the ESD
of $\tbC_\tau$ converges a.s. to a nonrandom distribution by showing that
$\tilde s_{\tau,p}$ converge pointwise a.s. to the Stieltjes
transform of a
probability measure,
%and in that case the Stieltjes transform $s_{\tau,p}$
%converges pointwise a.s. to the same limit,
thus establishing a.s. convergence of the ESD of $\bC_\tau$.
\end{TT*}

\begin{AT*}%{\sc Approximating equations for the Stieltjes kernel $K_\tau$.}
In order to
derive limiting equation (\ref{eqktau}), a finite sample counterpart is
needed. This can\vspace*{1pt} be derived as follows. The transformed data $\tbX$
gives rise
to the transformed Stieltjes kernel $\tK_\tau\dvtx\mathbb
{C}^+\times
[0,2\pi]\to\mathbb{C}^+$ given by
\[
\tK_\tau(z,\nu)=\frac{1}p\operatorname{tr} \bigl[
\tbR_\tau(z)\cH(\bA,\nu) \bigr].
\]
Following arguments typically used to establish the deterministic
equivalent of
a resolvent matrix, $p\times p$ matrix-valued function solutions
$(\bH_{\tau,p}(z)\dvtx p\in\mathbb{N})$ are needed such that, for
sufficiently
large $p$,
%
%e6.2 #&#
\begin{equation}
\label{eqapprox} \frac{1}p\operatorname{tr} \bigl[ \bigl( \bigl(\bI+
\bH_{\tau,p}(z) \bigr)^{-1}+z\tbR_\tau(z) \bigr)
\bD_p \bigr]\approx0
\end{equation}
for all $z\in\mathbb{C}^+$ and all $p\times p$ Hermitian matrix sequences
$(\bD_p\dvtx p\in\mathbb{N})$ with uniformly bounded norms $\llVert
\bD
_p\rrVert$. If one
uses $\bD_p=\cH(\bA,\nu)$ and the definition of $\tK_\tau$, the latter
approximate equation becomes
\[
\tK_\tau(z,\nu)\approx-\frac{1}{pz}\operatorname{tr} \bigl[ \bigl(
\bI+\bH_{\tau,p}(z) \bigr)^{-1}\cH(\bA,\nu) \bigr].
\]
Section~\ref{secgauss} below is devoted to making precise the use of
$\approx$
in the above equations and to showing that choosing
%
%e6.3 #&#
\begin{equation}
\label{eqhtaup} \bH_{\tau,p}(z) =-\frac{1}{zn}\sum
_{t=1}^n\frac{\cos(\tau\nu_t)\cH(\bA,\nu
_t)}{1+c_n\cos(\tau\nu_t)\tK_\tau(z,\nu_t)},
\end{equation}
with $\nu_t=2\pi t/n$ and $c_n=p/n$, is appropriate.
\end{AT*}
%Finally it should be noted that the approach taken here is similar to
%finding the
%deterministic equivalent of a resolvent matrix. However, since there
%is still randomness
%involved in $\bH_{\tau,p}$, the two approaches are not identical.

%s6.2 #&#
\subsection{Extension to the non-Gaussian case}\label{secintextnon-gauss}

In order to verify the statement of Theorem~\ref{thesdmainfty} for
non-Gaussian innovations $(Z_t\dvtx t\in\mathbb{Z})$, two key ideas are
invoked, namely showing that:
\begin{itemize}
\itemsep0ex
\item for any $z\in\mathbb{C}^+$, the Stieltjes transform $s_p(z)$
concentrates around $\mathbb{E}[s_p(z)]$ regardless of the
underlying distributional assumption;
\item the difference between the expectations $\mathbb{E}[s_p(z)]$
under the
Gaussian model and the non-Gaussian model is asymptotically
negligible.
\end{itemize}
To establish the concentration property of the first item, McDiarmid's
inequality is used to bound probabilities of the type
$\mathbb{P}(\llvert s_p(z)-\mathbb{E}[s_p(z)]\rrvert\geq\varepsilon)$
for arbitrary
\mbox{$\varepsilon>0$}. These probabilities are then shown to converge to zero
exponentially fast under (\ref{eqpn}). To establish the second item, the
generalized Lindeberg principle of Chatterjee~\cite{c06} is applied.
To this
end, the argument $z$ is viewed as a parameter and $s_p(z)$ as a
function of
the real\vspace*{1pt} parts $Z_{jt}^R$ and the imaginary parts $Z_{jt}^I$ of the innovation
entries $Z_{jt}$, for $j=1,\ldots,p$ and $t=1,\ldots,n$. The
difference between
Gaussian $Z_{t}$ and non-Gaussian $Z_{t}$ can then be analyzed by consecutively
changing one pair $(Z_{jt}^R,Z_{jt}^I)$ from Gaussian to non-Gaussian, thereby
expressing the respective differences in the expected Stieltjes
transforms as a sum of these entrywise changes. These differences will
be evaluated through a Taylor series expansion, bounding certain third-order
partial derivatives of $s_p$. Details are given in Section~\ref{secpnon-gauss}.

%s6.3 #&#
\subsection{Extension to the linear process case}\label{secintlin}

While the arguments established so far work in the same fashion also for
MA($q$) processes, certain difficulties arise when making the
transition to the
MA($\infty$) or linear process case. First, if one constructs the data matrix
$\bX$ not from $\mathrm{MA}(1)$ observations as above but from the linear process
$X_t=\sum_{\ell=0}^\infty\bA_\ell Z_{t-\ell}$, then every column
of $\bX$ is
different from the corresponding column in the transformed matrix
$\bX_\infty=\sum_{\ell=0}^\infty\bA_\ell\bZ\tbL^\ell$ and not
only the first
column [or the first $q$ columns for the MA($q$) case]. Second, for the $\mathrm{MA}(1)$
case one can write the Stieltjes transform $s_p(z)$ as a function of $2p(n+1)$
variables\vadjust{\goodbreak} $Z_{jt}^R$ and $Z_{jt}^I$ [or $2p(n+q)$ variables for the MA($q$)
case], but for linear processes, even for finite $p$, $s_p(z)$ is a
function of
infinitely many $Z_{jt}^R$ and $Z_{jt}^I$. This makes their study substantially
harder.

Linear processes are thus, for the purposes of this paper, approached through
truncation, that is, by approximation through finite-order MA processes
$X_t=\sum_{\ell=0}^{q(p)}\bA_\ell Z_{t-\ell}$ whose order $q(p)$ is
a function
of the dimension $p$ and therefore grows with the sample size under
(\ref{eqpn}). Obviously $q(p)\to\infty$ is a necessary condition to
make this
approximation work. However, $q(p)$ cannot grow too fast (leading to
the same difficulties in transitioning from the Gaussian to the non-Gaussian
case as for the linear process itself) or too slow (showing that the
LSDs of
the linear process and its truncated version are identical becomes an issue).
It turns out in Section~\ref{secplin} that $q(p)=\lceil
p^{1/3}\rceil$, with
$\lceil\cdot\rceil$ denoting the ceiling function, is an appropriate choice.

%s6.4 #&#
\subsection{Including the real-valued case}\label{secintreal}

To address the statements of Theorem~\ref{coesdmainfty}, the arguments
presented thus far have to be adjusted for real-valued innovations
$(Z_t\dvtx
t\in\mathbb{Z})$. This is done using the eigen-decomposition of the
coefficient
matrices in the real Fourier basis, after which arguments already
developed for
the complex case apply. Detailed steps are given in
Section~\ref{secpreal}.

%For the $\mathrm{MA}(1)$ case, the necessary changes in the proofs are
%based on the eigen-decomposition of the coefficient matrix, that is,
%on the
%decomposition $\bA_1=\bU_{\bA_1}\bolds{\Lambda}_{\bA_1}\bU_{\bA_1}^*$
%with an
%orthogonal eigenvector matrix $\bU_{\bA_1}$ and the real Fourier
%basis, say,
%$\bU$. With these, one can transform the original data matrix $X$ into
%the
%matrix $\tbX=\bU_{\bA_1}(\bZ+\bA_1\bZ\tbL)\bU$, which will be shown to
%have
%independent entries. Then arguments already developed for the complex
%case
%apply. The detailed steps are carried out in Section~{sec:p:real}.

%s6.5 #&#
\subsection{Dealing with the spectral density operator}\label{secspecdens}

The key step toward proving Theorem~\ref{thesdmainftyspectraldensity} is to express $\widehat\Gamma
_n(\eta)$
as
\begin{eqnarray*}
\widehat\Gamma_n(\eta) &=& \sum_{\tau=-\infty}^\infty
T_n(\tau) e^{i\tau\eta} \frac{1}{n} \mathbf{X}
\mathbf{L}^{-\tau} \mathbf{X}^*
\\
&=& \frac{1}{n} \mathbf{X} \Biggl(\bI_p +\sum
_{\tau=1}^\infty T_n(\tau)
\bigl(e^{i\tau\eta} \mathbf{L}^{-\tau} + e^{-i\tau\eta}
\mathbf{L}^{\tau}\bigr) \Biggr) \mathbf{X}^*,
\end{eqnarray*}
and then noticing that the matrix $\bI_p +\sum_{\tau=1}^\infty
T_n(\tau)
(e^{i\tau\eta} \tilde{\mathbf{L}}^\tau+ e^{-i\tau\eta}
\tilde{\mathbf{L}}^{-\tau})$ diagonalizes in the (real or complex) Fourier
basis with eigenvalues $\mathfrak{f}_{T_n}(\eta- \nu_t)$, for
$t=1,\ldots,n$,
so that the ESD of $\widehat\Gamma_n(\eta)$ can be\vspace*{1pt} approximated by
the ESD of
the matrix $\widetilde\Gamma_n(\eta) = n^{-1} \sum_{t=1}^n
\mathfrak{f}_{T_n}(\eta- \nu_t) \tilde X_t \tilde X_t^*$, where
$\tilde X_t^*$ is the $t$th column of the matrix $\mathbf
{X}\mathbf{F}_n^*$,
and $\mathbf{F}_n$ denotes the $n\times n$ Fourier rotation matrix. We
give the main
steps of the arguments leading to this result. First, suppose that $m_n
\to\infty$
such that $m_n/n \to0$ as $n\to\infty$. Then we can write
%
%e6.4 #&#
\begin{equation}
\label{eqhatGammanetarepresentation} \widehat\Gamma_n(\eta) = \widetilde
\Gamma_{n,m_n}^{(1)}(\eta) + \bigl(\widehat\Gamma_{n,m_n}^{(1)}(
\eta) - \widetilde\Gamma_{n,m_n}^{(1)}(\eta)\bigr) + \widehat
\Gamma_{n,m_n}^{(2)}(\eta),
\end{equation}
where
\begin{eqnarray*}
\widetilde\Gamma_{n,m_n}^{(1)}(\eta) &=& \frac{1}{n}
\mathbf{X} \Biggl(\bI_p +\sum_{\tau=1}^{m_n}
T_n(\tau) \bigl(e^{i\tau\eta} \tilde{\mathbf{L}}^{-\tau} +
e^{-i\tau\eta} \tilde{\mathbf{L}}^{\tau}\bigr) \Biggr) \mathbf{X}^*
\\
&=& \frac{1}{n} \sum_{t=1}^n
\mathfrak{f}_{T_n,m_n}(\eta-\nu_t) \tilde X_t
\tilde X_t^*,
\end{eqnarray*}
with $\mathfrak{f}_{T_n,m_n}(\theta) = 1 + 2\sum_{\tau=1}^{m_n}
\cos(\tau\theta)$,
\[
\widehat\Gamma_{n,m_n}^{(1)}(\eta) = \frac{1}{n}
\mathbf{X} \Biggl(\bI_p +\sum_{\tau=1}^{m_n}
T_n(\tau) \bigl(e^{i\tau\eta} \mathbf{L}^{-\tau} +
e^{-i\tau\eta} \mathbf{L}^{\tau}\bigr) \Biggr) \mathbf{X}^*
\]
and
\[
\widehat\Gamma_{n,m_n}^{(2)}(\eta) = \frac{1}{n}
\mathbf{X} \Biggl(\sum_{\tau=m_n+1}^n
T_n(\tau) \bigl(e^{i\tau\eta} \mathbf{L}^{-\tau} +
e^{-i\tau\eta} \mathbf{L}^{\tau}\bigr) \Biggr) \mathbf{X}^*.
\]
Now, the following facts together with representation (\ref
{eqhatGammanetarepresentation})
and Theorem~A.43 of Bai and Silverstein \cite{bs10} (rank inequality)
and Lemma~S.1 %\ref{lemnorminequality}
(norm inequality) prove the assertion:
\begin{longlist}[(iii)]
\item[(i)]
$\sup_{\theta\in[0,2\pi]}\max\{\llvert\mathfrak{f}_{T_n,m_n}(\theta) -
\mathfrak{f}_{T_n}(\theta)\rrvert,\llvert\mathfrak{f}_{T_n}(\theta
)-\mathfrak
{f}_T(\theta)\rrvert\}
\to0$, as $n\to\infty$;\vspace*{3pt}

\item[(ii)]
rank$(\widehat\Gamma_{n,m_n}^{(1)}(\eta) - \widetilde\Gamma
_{n,m_n}^{(1)}(\eta)) \leq2m_n = o(p)$;\vspace*{3pt}

\item[(iii)]
\[
p^{-1}\bigl\|\widehat\Gamma_{n,m_n}^{(2)}(\eta)\bigr\|_F^2
\leq\bigl\|\widehat\Gamma_{n,m_n}^{(2)}(\eta)\bigr\|^2
\]
and
\[
\big\|\widehat\Gamma_{n,m_n}^{(2)}(\eta)\big\|\leq
\big\| n^{-1}\mathbf{X}\mathbf{X}^*\big\| \sum_{\tau
=m_n+1}^\infty
\big\llvert T_n(\tau)\big\rrvert\to0\qquad\mbox{a.s.}
\]
\end{longlist}

%s7 #&#
\section{Proof for the complex Gaussian $\mathrm{MA}(q)$ case}\label{secgauss}

Throughout, $(\mathbf{A}_\ell\dvtx\ell\in\mathbb{N})$ is treated
as a
sequence of nonrandom matrices, and all the arguments are valid
conditionally on
this sequence. In this section, the result of Theorem~\ref
{thesdmainfty} is
first verified for the MA($q$) process $X_t=\sum_{\ell=0}^q\bA_\ell
Z_{t-\ell}$
when $Z_{jt}$'s are i.i.d. standard complex Gaussian; that is, real
and imaginary parts of $Z_{jt}$ are independent normals with mean zero
and variance one half. Following the outline in
Section~\ref{subsecintu}, the data matrix $\bX=[X_1\dvtx\cdots\dvtx
X_n]$ is
transformed into the matrix $\tbX=[\tX_1\dvtx\cdots\dvtx\tX_n]$, with
each column
satisfying $\tX_t=\psi(\bA,\nu_t)\tZ_t$. Then, since the rank
of $\mathbf{X} - (\mathbf{Z}+\sum_{\ell=1}^q \mathbf{A}_\ell
\mathbf{Z}\tilde
{\mathbf{L}}^\ell)$ is at most $q$, by Lemma~S.2, %
it follows that the ESDs of $\mathbf{C}_\tau$ and $\tilde{\mathbf
{C}}_\tau$ have
the same limit, if the latter exists. For simplicity, let
$\bpsi_t=\psi(\bA,\nu_t)$. To keep notation more compact, the extra
subscripts
$p$ and~$\tau$ (indicating the lag of the autocovariance matrix under
consideration) are often suppressed when no confusion can arise. For example,
in (\ref{eqapprox}) the notation $\bH(z)$ will be preferred over the more
complex $\bH_{\tau,p}(z)$. The proof is given in several steps.
First, a bound
on the approximation error is derived if the Stieltjes kernel $K$ in~(\ref{eqktau}) is replaced with its finite sample counterpart $\tK$. Second,
existence, convergence and continuity of the solution to (\ref
{eqapprox}) are
verified. Third, tightness of the ESDs $\hat F_p$ and convergence of the
corresponding Stieltjes transforms $s_p$ is shown.

%s7.1 #&#
\subsection{Bound on the approximation error}\label{secgaussbound}

The goal is to provide a rigorous formulation of (\ref{eqapprox}) and a
bound on the resulting approximation error. The first step consists of
giving a
heuristic argument for the definition of $\bH(z)$ in (\ref
{eqhtaup}). To this\vadjust{\goodbreak}
end, note that
%
%e7.1 #&#
\begin{eqnarray}\label{eq*}
&& \bigl(\bI+\bH(z) \bigr)^{-1} +z\tbR(z)\nonumber
\\
&&\qquad = \tbR(z)\tbR(z)^{-1} \bigl(\bI+\bH(z) \bigr)^{-1}-\tbR(z)
\bigl(\bI+\bH(z) \bigr) \bigl(\bI+\bH(z) \bigr)^{-1}
\\
&&\qquad = \tbR(z) \bigl(\tbC_\tau+z\bH(z) \bigr) \bigl(\bI+\bH(z)
\bigr)^{-1}.
\nonumber
\end{eqnarray}
It follows that, to achieve (\ref{eqapprox}), it is sufficient to solve
$\tbR(z)\tbC_\tau\approx-z\tbR(z)\bH(z)$. (The use of $\approx$
will be clarified below.) Let
\[
\gamma_{\tau,t}=\cos(\tau\nu_t),
\]
and define the rank-one perturbation and its corresponding resolvent,
respectively, given by
\[
\tbC_{\tau,t}=\tbC_\tau-\frac{1}n\gamma_{\tau,t}
\tX_t\tX_t^* \quad\mbox{and}\quad
\tbR_t(z)=(\tbC_{\tau,t}-z\bI)^{-1}.
\]
Using
$\tbR(z)=([\tbC_{\tau,t}-z\bI]+\frac{1}n\gamma_{\tau,t}X_tX_t^*)^{-1}$
and defining
\[
\cH_t=\cH(\bA,\nu_t),
\]
the Sherman--Morrison formula and some matrix algebra lead to
%
%e7.2 #&#
\begin{eqnarray}
\tbR(z)\tbC_\tau&=&\frac{1}n\sum
_{t=1}^n \biggl(\tbR_t(z)-
\frac{(1/n)\gamma_{\tau,t}\tbR_t(z)\tX_t\tX_t^*\tbR
_t(z)}{1+(1/n)\gamma_{\tau,t}\tX_t^*\tbR_t(z)\tX_t} \biggr) \gamma
_{\tau,t}\tX_t
\tX_t^*
\nonumber
\\
&=&\frac{1}n\sum_{t=1}^n
\tbR_t(z)\frac{\gamma_{\tau,t}\tX_t\tX
_t^*}{1+(1/n)\gamma_{\tau,t}\tX_t^*\tbR_t(z)\tX_t}
\nonumber
\\
&=&\frac{1}n\sum_{t=1}^n
\tbR_t(z)\frac{\gamma_{\tau,t}\bpsi_t\tZ
_t\tZ_t^*\bpsi_t^*}{
1+(1/n)\gamma_{\tau,t}\tZ_t^*\bpsi_t^*\tbR_t(z)\bpsi_t\tZ_t} \label{eq**}
\\
&\approx&\tbR(z)\frac{1}n\sum_{t=1}^n
\frac{\gamma_{\tau,t}\mathcal
{H}_t}{1+(1/n)\gamma_{\tau,t}\operatorname{tr}[\tbR(z)\mathcal{H}_t]}
\nonumber
\\
&=&\tbR(z)\frac{1}n\sum_{t=1}^n
\frac{\gamma_{\tau,t}\cH_t}{1+c_n\gamma_{\tau
,t}\tK(z,\nu_t)},
\nonumber
\end{eqnarray}
thus\vspace*{1pt} validating the choice of $\bH(z)$ as given in (\ref{eqhtaup}). The
$\approx$ sign is due to substituting $\tbR_t(z)$, $\tZ_t\tZ_t^*$,
$\tZ_t^*\bpsi_t^*\tbR_t(z)\bpsi_t\tZ_t$ with $\tbR(z)$, $\bI$,
$\operatorname{tr}[\bpsi_t^*\tbR(z)\bpsi_t]$, respectively.
%the normalized trace of the
%difference between the sides multiplied by a nonrandom Hermitian
%matrix with
%bounded norm converges to zero in an appropriate sense.}
The arguments are made precise in the following theorem.

%th7.1 #&#
\begin{theorem}
\label{thapperr}
Let\/ $\bD$ be a Hermitian matrix such that $\llVert\bD\rrVert\leq
\bar
{\lambda}_{\bD}$
and $z\in\mathbb{C}^+$. If the assumptions of Theorem~\ref
{thesdmainfty} are satisfied, then
%
%e7.3 #&#
\begin{equation}
\label{eq***} \frac{1}{zp}\operatorname{tr} \bigl[ \bigl( \bigl(\bI+\bH(z)
\bigr)^{-1}+z\tbR(z) \bigr)\bD\bigr] \to0 \qquad\mbox{a.s.}
\end{equation}
under (\ref{eqpn}), where $\bH$ is defined in (\ref{eqhtaup}), and
$\tbR(z)=(\tbC_\tau-z\bI)^{-1}$ is the resolvent of the %transformed
symmetrized autocvariance matrix $\tbC_\tau=\frac{1}n\sum_{t=1}^n\gamma
_{\tau,t}\tX_t\tX_t^*$ with
$\gamma_{\tau,t}=\cos(\tau\nu_t)$.
\end{theorem}

\begin{pf}
Observe that using (\ref{eq*}), (\ref{eq**}) and the definition of
$\bH$ in
(\ref{eqhtaup}), the a.s. convergence in (\ref{eq***}) is shown to be
equivalent to
\[
\bar{d}_{\tau}^{(n)}=\frac{1}n\sum
_{t=1}^n d_{\tau,t}\to0 \qquad\mbox{a.s.},
\]
where
\[
d_{\tau,t}= \frac{1}{zp}\operatorname{tr} \bigl[\gamma_{\tau,t}
\bigl( \tilde\beta_{\tau,t}\tbR_t(z)\tcH_t
-\beta_{\tau,t}\tbR(z)\cH_t %\bpsi_t\bpsi_t^*
\bigr)
\bigl(\bI+\bH(z) \bigr)^{-1}\bD\bigr],
\]
with $\tcH_t=\bpsi_t\tZ_t\tZ_t^*\bpsi_t^* = \tilde X_t \tilde X_t^*$,
%
%e7.4 #&#
\begin{eqnarray}\label{eqbeta}
\tilde\beta_{\tau,t} &=& \biggl(1+\frac{1}n
\gamma_{\tau,t}\tX_t^*\tbR_t(z)
\tX_t \biggr)^{-1},
\nonumber\\[-8pt]\\[-8pt]\nonumber %\quad\mbox{and}\quad
\beta_{\tau,t} &=& \biggl(1+\frac{1}n\gamma_{\tau,t}\operatorname{tr}\bigl[\tbR(z)
\cH_t\bigr] \biggr)^{-1}.
\end{eqnarray}
Decomposing further, write next $d_{\tau,t}=d_{\tau,t}^{(1)}+\cdots
+d_{\tau,t}^{(5)}$, where
\begin{eqnarray*}
d_{\tau,t}^{(1)} &=&\frac{1}{zp}\operatorname{tr} \bigl[
\gamma_{\tau,t}\tilde{\beta}_{\tau
,t}\tbR_t(z)
\tcH_t \bigl( \bigl(\bI+\bH(z) \bigr)^{-1}- \bigl(\bI+
\bH_t(z) \bigr)^{-1} \bigr)\bD\bigr],
\\
d_{\tau,t}^{(2)} &=&\frac{1}{zp}\operatorname{tr} \bigl[
\gamma_{\tau,t}\tilde{\beta}_{\tau,t}\tbR_t(z) (
\tcH_t-\cH_t ) \bigl(\bI+\bH_t(z)
\bigr)^{-1}\bD\bigr],
\\
d_{\tau,t}^{(3)} &=&\frac{1}{zp}\operatorname{tr} \bigl[
\gamma_{\tau,t}\tilde{\beta}_{\tau,t} \tbR_t(z)
\cH_t \bigl( \bigl(\bI+\bH_t(z) \bigr)^{-1}-
\bigl(\bI+\bH(z) \bigr)^{-1} \bigr)\bD\bigr],
\\
d_{\tau,t}^{(4)} &=&\frac{1}{zp}\operatorname{tr} \bigl[
\gamma_{\tau,t}\tilde\beta_{\tau,t} \bigl(\tbR_t(z)-
\tbR(z) \bigr)\cH_t \bigl(\bI+\bH(z) \bigr)^{-1}\bD\bigr],
\\
d_{\tau,t}^{(5)} &=&\frac{1}{zp}\operatorname{tr} \bigl[
\gamma_{\tau,t} (\tilde\beta_{\tau,t}-\beta_{\tau,t} )
\tbR(z)\cH_t \bigl(\bI+\bH(z) \bigr)^{-1}\bD\bigr],
\end{eqnarray*}
with
%
%e7.5 #&#
\begin{equation}
\label{eqbetabart} \qquad \bH_t=-\frac{1}{zn}\sum
_{t=1}^n\bar{\beta}_{\tau,t}
\gamma_{\tau,t}\mathcal{H}_t, \qquad%\mbox{and}\qquad
\bar{
\beta}_{\tau,t}= \biggl(1+\frac{1}n\gamma_{\tau,t}\operatorname{tr}\bigl
[\tbR_t(z)\mathcal{H}_t\bigr]
\biggr)^{-1},
\end{equation}
thus exhibiting the various approximations being made in the proof.\vspace*{1pt}

The Borel--Cantelli lemma provides that $\bar{d}_{\tau}^{(n)} \to0$
almost surely is
implied if $\mathbb{P}(\max_{t\leq n}\llvert{d}_{\tau,t}\rrvert
>\varepsilon)\to0$
%%at a rate
faster than $\frac{1}n$ for all $\varepsilon>0$. Since $\mathbb{P}(\max
_{t\leq
n}\llvert{d}_{\tau,t}\rrvert>\varepsilon)\leq\sum_{\ell=1}^5\mathbb
{P}(\max_{t\leq
n}\llvert{d}_{\tau,t}^{(\ell)}\rrvert>\frac{\varepsilon}{5})$, in
order to verify
(\ref{eq***}), it is sufficient to show that $\mathbb{P}(\max_{t\leq
n}\llvert{d}_{\tau,t}^{(\ell)}\rrvert>\varepsilon)$ goes to zero
faster than $\frac{1}n$ for
all $\varepsilon>0$ and $\ell=1,\ldots,5$. The corresponding arguments are
detailed in Section S1.1 %\ref{secappapperr}
of the online SM.
\end{pf}

%s7.2 #&#
\subsection{Existence and uniqueness of the solution}\label{secgaussexconcon}

In this section, the proof of Theorem~\ref{thesdmainfty} is
completed for
the complex Gaussian innovation model. In what follows, $\bA$ is
without loss
of generality assumed to be nonrandom, thereby restricting randomness
to the
innovations $\bZ=[Z_1\dvtx\cdots\dvtx Z_n]$. For a fixed $\omega$ in
the underlying
sample space $\Omega$, notation such as $\bZ(\omega)$ will be
utilized to
indicate realizations of the respective random quantities.

Noticing first that Theorem~\ref{thesdmainfty} makes an almost sure
convergence statement, a~suitable subset $\Omega_0\subset\Omega$ with
$\mathbb{P}(\Omega_0)=1$ is determined. This subset is used for all subsequent
arguments. To this end, observe that since
the matrix $\tbZ=[\tZ_1\dvtx\cdots\dvtx\tZ_n]$ has i.i.d. entries with
zero mean and unit variance
the norm of $n^{-1}\tbZ\tbZ^*$ converges almost surely to a number
not exceeding $(1+\sqrt{c})^2+1$. Let
${F}_{\tbZ}$ denote the ESD of $n^{-1}\tbZ\tbZ^*$, and define
\[
\Omega_1= \bigl\{\omega\in\Omega\dvtx F_{\tbZ(\omega)}< (1+\sqrt
{c})^2 + 1\mbox{ for } p\geq p_0(\omega)\bigr\},
\]
with suitably chosen $p_0(\omega)$. Then $\mathbb{P}(\Omega_1)=1$.
Define next
\begin{eqnarray*}
\bar{d}(z,\nu,\omega) &=&\frac{1}{zp}\operatorname{tr} \bigl[ \bigl(
\bigl(\bI+
\bH(z) \bigr)^{-1}+z\tbR(z,\omega) \bigr)\cH(\bA,\nu) \bigr],
\\
\bar{d}_\bI(z,\omega) &=&\frac{1}{zp}\operatorname{tr} \bigl[ \bigl(
\bI+\bH(z) \bigr)^{-1}+z\tbR(z,\omega) \bigr].
\end{eqnarray*}
Let\vspace*{1pt} $\mathbb{C}_\mathbb{Q}^+$ denote the set of complex numbers with rational
real part and positive rational imaginary part and
$[0,2\pi]_\mathbb{Q}=[0,2\pi]\cap\mathbb{Q}$. Define the set
\[
\Omega_2= \bigl\{ \omega\in\Omega\dvtx\bar{d}(z,\nu,\omega)\to0,
\bar{d}_\bI(z,\omega)\to0, z\in\mathbb{C}_\mathbb{Q}^+, \nu
\in[0,2\pi]_\mathbb{Q} \bigr\}.
\]
In view of Theorem~\ref{thapperr}, it follows that $\bar{d}(z,\nu,\cdot
)\to
0$ a.s. and $\bar{d}_\bI(z,\cdot)\to0$ a.s. for all fixed
$z\in\mathbb{C}_\mathbb{Q}^+$ and $\nu\in[0,2\pi]_\mathbb{Q}$. Thus\vspace*{1pt}
$\mathbb{P}(\Omega_2)=1$. Henceforth only
$\omega\in\Omega_0=\Omega_1\cap\Omega_2$, so that $\mathbb
{P}(\Omega_0)=1$, are
considered.

Recall that
$\tK_\tau(z,\nu,\omega)=p^{-1}\operatorname{tr}[\tbR(z,\omega)\cH(\bA
,\nu)]$. The
following theorem establishes existence of a Stieltjes kernel solution
to the
equations in (\ref{eqktau}) along a subsequence.

%th7.2 #&#
\begin{theorem}[(Existence)]\label{thexist}
Suppose that the assumptions of Theorem~\ref{thesdmainfty} are satisfied:
\begin{longlist}[(a)]
\item[(a)] For\vspace*{1pt} all $\omega\in\Omega_0$ and for all subsequences of $\{p\}$, there
exists another subsequence $\{p_\ell\}$ along which $\tilde
K_\tau(z,\nu,\omega)$ converges pointwise in $z\in\mathbb{C}^+$
and uniformly
in $\nu\in[0,2\pi]$ to a limit $\kK_\tau(z,\nu,\omega)$ analytic
in $z$ and
continuous in $\nu$.

\item[(b)] For every subsequence $\{p_\ell\}$ satisfying \textup{(a)}, $\kK_\tau
(z,\nu,\omega)$
satisfies (\ref{eqktau}) for any $z\in\mathbb{C}^+$. Moreover,
$\kK_\tau(z,\nu,\omega)$ is the Stieltjes transform of a measure
with mass
$m_\nu=\int h(\lambda,\nu)\,dF^\bA(\lambda)$, provided that $m_\nu>0$.
\end{longlist}
\end{theorem}

\begin{pf}
(a) Let $z=w+iv$. Then for a compact subset $S\subset\mathbb{C}^+$,
$\tK_\tau(z,\nu,\break \omega)\leq\bar{\lambda}_\bA^2(\min_{z\in
S}v)^{-1}$ for all
$\nu\in[0,2\pi]$ and $z\in S$. Enumerate $[0,2\pi]_\mathbb{Q}=\{
\nu_\ell\dvtx
\ell\in\mathbb{N}\}$. Let $\{p_\ell^j(\omega)\}\subset\{p_\ell
^{j-1}(\omega)\}$
mean that $\{p_\ell^j(\omega)\}$ is a further subsequence of
$\{p_\ell^{j-1}(\omega)\}$, and let $\{p\}$ denote the original
sequence. An
application of Lemma~3 in Geronimo and Hill \cite{gh03} yields that,
for any
fixed $\omega\in\Omega_0$, there is a sequence of subsequences
\[
\bigl\{p_\ell^j(\omega)\bigr\}\subset\bigl
\{p_\ell^{j-1}(\omega)\bigr\}\subset\cdots\subset\bigl
\{p_\ell^1(\omega)\bigr\}\subset\{p\},
\]
so that $K_\tau^{(p_\ell^j)}(z,\nu,\omega)$ converges to an
analytic function
of $z$ on $\{\nu_\ell\dvtx\ell=1,\ldots,\ell\}$. A standard
diagonal argument
implies that $K_\tau^{(p_\ell^\ell)}(z,\nu,\omega)$ converges to
an analytic
function of $z$ on $[0,2\pi]_\mathbb{Q}$. To simplify notation, write
$K_\tau^{(p_\ell)}(z,\break \nu,\omega)$ in place of
$K_\tau^{(p_\ell^\ell)}(z,\nu,\omega)$. Observe that the thus
obtained limit,
which will be denoted by $\kK_\tau(z,\nu,\omega)$, is so far
defined only on
$\mathbb{C}^+\times[0,2\pi]_\mathbb{Q}$. It remains to obtain the
extension of
the limit to $\mathbb{C}^+\times[0,2\pi]$. Note that Lemma~S.10 %
implies that, for any $z\in\mathbb{C}^+$, $K_\tau^{(p_\ell)}(z,\nu
,\omega)$ are
equicontinuous in $\nu$ and converge pointwise to $\kK_\tau(z,\nu,\omega
)$ on
the dense subset $[0,2\pi]_\mathbb{Q}$ of $[0,2\pi]$. By the Arzela--Ascoli
theorem, $K_\tau^{(p_\ell)}(z,\nu,\omega)$ converges therefore
uniformly to a
continuous function of $\nu$ on $\mathbb{C}^+\times[0,2\pi]$. This limit,
denoted again by $\kK_\tau(z,\nu,\omega)$, is also analytic on
$\mathbb{C}^+$.
To see this, pick $\nu_0\in[0,2\pi]\setminus\mathbb{Q}$ and a sequence
$\{\nu_p\}\subset[0,2\pi]_\mathbb{Q}$ such that $\nu_p\to\nu_0$. Then
$\{\kK_\tau(z,\nu_p,\omega)\}$ satisfies the conditions of Lemma~3 in
\cite{gh03}, and consequently there is a subsequence of
$\{\kK_\tau(z,\nu_p,\omega)\}$ that converges to an analytic
function. The
limit of this subsequence has to coincide with $\kK_\tau(z,\nu
_0,\omega)$ by
continuity on $[0,2\pi]$. It follows that $\kK_\tau(z,\nu,\omega)$
analytic.

(b) By Lemma~S.11 %\ref{lemexist2}
and the definition of $\Omega_0$, for all
$\omega\in\Omega_0$, $\kK_\tau(z,\nu,\omega)$ satisfies~(\ref{eqktau}) for all
$z\in\mathbb{C}^+_\mathbb{Q}$. So, by analyticity in $z$, it holds
for all $\mathbb{C}^+$.

Suppose first that $m_\nu=0$. Then $\llvert K_\tau(z,\nu,\omega
)\rrvert\leq
%(vp)^{-1}\operatorname{tr}(h(\bA,\nu))=
v^{-1}\int h(\lambda,\nu)\,dF_p^\bA(\lambda)\to
v^{-1}\int h(\lambda,\nu)\,dF^\bA(\lambda)=0$. Thus $\kK(z,\nu,\omega)=0$
for all
$z\in\mathbb{C}^+$, and the claim is verified.

For the remainder, suppose that $m_\nu>0$. Showing that
$\kK_\tau(z,\nu,\omega)$ is a Stieltjes transform of a measure with mass
$m_\nu=\int h(\lambda,\nu)\,dF^\bA(\lambda)$ is equivalent to
showing that
$m_\nu^{-1}\kK_\tau(z,\nu,\omega)$ is Stieltjes transform of a Borel
probability measure. Let $\bLambda_{\tbC}$ and $\bU_{\tbC}$ denote the
eigenvalue and eigenvector matrices of $\tbC_\tau$. Then
$\tK_\tau^{(p_\ell)}(z,\nu,\omega)=p_\ell^{-1}\operatorname
{tr}[(\bLambda_{\tbC}-z\bI)^{-1}\bU_{\tbC}^*h(\bA,\nu)\bU_{\tbC}]$
is the Stieltjes transform of a measure with mass
$m_\nu^{(p_\ell)}=p_\ell^{-1}\operatorname{tr}[h(\bA,\nu)]$. By the
weak convergence
of $F^\bA_{p_\ell}$ to $F^\bA$, $m_\nu^{(p_\ell)}\to m_\nu>0$ as
$p_\ell\to\infty$. This shows that $m_\nu^{(p_\ell)}\geq\bar
{m}_\nu>0$ for all
$p_\ell$. Since the diagonal entries of $\bU_{\tbC}^*h(\bA,\nu)\bU
_{\tbC}$ are
bounded from above by $\bar{\lambda}_\bA^2$, it follows that
$(m_\nu^{(p_\ell)})^{-1}\tK_\tau^{(p_\ell)}(z,\nu,\omega)$ is
the Stieltjes
transform of a measure $\mu_{p_\ell}$, say, such that, for all real $x$,
$\mu_{p_\ell}((x,\infty))\leq
\bar{\lambda}_\bA^2\bar{m}_\nu^{-1}F^{\tbC_\tau}((x,\infty))$,
where $F^{\tbC_\tau}$ denotes the ESD of $\tbC_\tau$. It follows
from Lemma~S.12 %\ref{lemexist3}
that $\{F^{\tbC_\tau}\}$ is a tight sequence. Therefore
$(m_\nu^{(p_\ell)})^{-1}\times\break \tK_\tau^{(p_\ell)}(z,\nu,\omega)$ are
the Stieltjes
transforms of a tight sequence of Borel measures. An application of
Lemma~S.13 %\ref{lemunique1}
yields that $K_\tau(z,\nu,\omega)$ is the Stieltjes
transform of a measure with mass $m_\nu$, completing the proof.
\end{pf}

%th7.3 #&#
\begin{theorem}[(Uniqueness)]\label{thunique}
If the assumptions of Theorem~\ref{thesdmainfty} are satisfied,
then there
is a unique solution $K_\tau(z,\nu)$ to (\ref{eqktau}) that is
analytic in
$z\in\mathbb{C}^+$ and continuous in $\nu\in[0,2\pi]$ with $K_\tau
(z,\nu)$
being a Stieltjes transform of a measure with mass $\int
h(\lambda,\nu)\,dF^\bA(\lambda)$.
\end{theorem}

\begin{pf}
Suppose there are two solutions $K_{\tau,1}(z,\nu)$ and $K_{\tau
,2}(z,\nu)$ to
(\ref{eqktau}). Let $\gamma_\tau(\nu)=\cos(\tau\nu)$ and
$\beta_{\tau,j}(z,\nu)=(1+c\gamma_\tau(\nu)K_{\tau,j}(z,\nu
))^{-1}$. Define
then
$U_{\tau,j}(z,\lambda)=(2\pi)^{-1}\int_0^{2\pi}\beta_{\tau,j}(z,\lambda
)\gamma_{\tau}(\nu)h(\lambda,\nu)\,d\nu$,
$j=1,2$. Note that\break $U_{\tau,1}(z,\lambda)$ and $U_{\tau,2}(z,\lambda
)$ have
nonpositive imaginary parts. Now
\begin{eqnarray*}
&& U_{\tau,1}(z,\lambda)-U_{\tau,2}(z,\lambda)
\\
&&\qquad =\frac{c}{2\pi}\int_0^{2\pi}
\beta_{\tau,1}(z,\nu)\beta_{\tau,2}(z,\nu)\gamma_{\tau}^2(
\nu) \bigl(K_{\tau,1}(z,\nu)-K_{\tau,2}(z,\nu) \bigr)h(\lambda,\nu)\,d
\nu
\end{eqnarray*}
and thus
\begin{eqnarray*}
\hspace*{-4pt}&&K_{\tau,1}(z,\nu)-K_{\tau,2}(z,\nu)
\\
\hspace*{-4pt}&&\qquad =-\int\frac{(U_{\tau,1}(z,\lambda)-U_{\tau,2}(z,\lambda
))h(\lambda,\nu)}{(U_{\tau,1}(z,\lambda)-z)(U_{\tau,2}(z,\lambda
)-z)}\,dF^\bA(\lambda)
\\
\hspace*{-4pt}&&\qquad =-\int\frac{c}{2\pi}\int_0^{2\pi}\prod
_{j=1}^2\frac{\beta
_{\tau,j}(z,\nu^\prime)}{U_{\tau,j}(z,\lambda)-z}
\gamma_{\tau
}^2\bigl(\nu^\prime\bigr)
\\
\hspace*{-4pt}&&\quad\qquad{}\hspace*{75pt} \times\bigl(K_{\tau,1}\bigl(z,\nu^\prime\bigr)-K_{\tau,2}
\bigl(z,\nu^\prime\bigr) \bigr)h\bigl(\lambda,\nu^\prime\bigr)\,d
\nu^\prime h(\lambda,\nu)\,dF^\bA(\lambda).
\end{eqnarray*}
Using the fact that $K_{\tau,j}(z,\nu)$ is a Stieltjes transform with mass
bounded from above by $\bar{\lambda}_\bA^2$, it follows that
$\prod_{j=1}^2\beta_{\tau,j}(z,\nu^\prime)(U_{\tau,j}(z,\lambda
)-z)^{-1}$ is
bounded by $C(v)=\max\{64c^2\bar{\lambda}_\bA^4v^{-4},4v^{-2}\}$. Thus
\begin{eqnarray*}
&& \int_0^{2\pi}\bigl\llvert K_{\tau,1}(z,
\nu)-K_{\tau,2}(z,\nu)\bigr\rrvert \,d\nu
\\
&&\qquad \leq4C(v)c\bar{
\lambda}_\bA^4\int_0^{2\pi}
\bigl\llvert K_{\tau,1}\bigl(z,\nu^\prime\bigr)-K_{\tau,2}
\bigl(z,\nu^\prime\bigr)\bigr\rrvert \,d\nu^\prime.
\end{eqnarray*}
If $v>4\bar{\lambda}^2_\bA\max\{c^{3/4},\sqrt{c}\}$, then
$4C(v)c\bar{\lambda}_\bA^4<1$ and thus
\[
\int_0^{2\pi}\bigl\llvert K_{\tau,1}(z,\nu)-K_{\tau,2}(z,\nu)\bigr\rrvert \,d\nu=0,
\]
which by continuity
in $\nu$ implies that $K_{\tau,1}(z,\nu)=K_{\tau,2}(z,\nu)$ for
any fixed
$z\in\mathbb{C}^+$. Since both solutions are analytic, the equality holds
indeed for all $z\in\mathbb{C}^+$. This proves uniqueness.
\end{pf}

In the remainder of this section, the proof of Theorem~\ref
{thesdmainfty} is
completed for the Gaussian MA($q$) case. This is done by establishing
that (a)
the convergence along subsequences as stated in Theorem~\ref{thexist} holds
indeed for the whole sequence and (b) the relevant ESDs converge.

Toward (a), it is necessary to prove that, for every $\omega\in\Omega_0$,
$\tK_\tau(z,\nu,\omega)$ converges to $K_\tau(z,\nu)$ pointwise in
$z\in\mathbb{C}^+$ and uniformly for $\nu\in[0,2\pi]$ under~(\ref{eqpn}).
Assume the contrary,\vspace*{1pt} and suppose that there are $z_0\in\mathbb{C}^+$,
$\nu_0\in[0,2\pi]$ and $\omega_0\in\Omega_0$ such that
$\tK_\tau(z_0,\nu_0,\omega_0)$ does not converge to $K_\tau
(z_0,\nu_0)$. By
boundedness of $\tK_\tau(z_0,\nu_0,\omega_0)$, there is a subsequence
$\{p_\ell\}$ along which $\tK_\tau(z_0,\nu_0,\omega_0)$ converges
to a limit
different from $K_\tau(z_0,\nu_0)$. Invoking Theorems~\ref{thexist} and
\ref{thunique}, there is a further subsequence $\{p_{\ell^\prime}\}
$ of
$\{p_\ell\}$ along which $\tK_\tau(z_0,\nu,\omega_0)$ converges to
$K_\tau(z_0,\nu)$ uniformly in $\nu\in[0,2\pi]$. This is a
contradiction. It
follows that for every $\omega\in\Omega_0$, $\tK_\tau(z,\nu,\omega)$ converges
to $K_\tau(z,\nu)$ pointwise in $z\in\mathbb{C}^+$ and $\nu\in
[0,2\pi]$. An
application of Theorem~\ref{thexist} and the Arzela--Ascoli theorem
shows that
the convergence is uniform on $[0,2\pi]$. Note that, for any
$z\in\mathbb{C}^+$,
$K_{\tau,p}(z,\nu,\omega)=p^{-1}\operatorname{tr}[\bR(z)\cH(\bA,\nu
)]$ converges to
$K_\tau(z,\nu)$ uniformly on $[0,2\pi]$. Since we have
$\llvert\tK_\tau(z,\nu,\omega)-K_{\tau,p}(z,\nu,\omega)\rrvert<\bar
{\lambda
}_\bA^2(vp)^{-1}$,
assertion (a) follows.

Toward (b), let $\omega\in\Omega_0$. It needs to be shown that
$s_{\tau,p}(z,\omega)\to s_\tau(z)$ on~$\mathbb{C}^+$. By arguments
as in the
proof of Lemma~S.11, %\ref{lemexist2},
it is already established that
$\ts_{\tau,p}(z,\omega)\to s_\tau(z)$ on $\mathbb{C}_\mathbb
{Q}^+$. Now, for
any compact $S\subset\mathbb{C}^+$ and $z_1,z_2\in S$,
\begin{eqnarray*}
\bigl\llvert\ts_{\tau,p}(z_1,\omega)-\ts_{\tau,p}(z_2,
\omega)\bigr\rrvert&=& \frac{1}p\bigl\llvert\operatorname{tr} \bigl[
\tbR(z_1,\omega)-\tbR(z_2,\omega) \bigr]\bigr\rrvert
\\
&=&\frac{1}{p}\bigl\llvert\operatorname{tr} \bigl[(z_1-z_2)
\tbR(z_1,\omega)\tbR(z_2,\omega) \bigr]\bigr\rrvert
\\
&\leq& \llvert z_1-z_2\rrvert\Bigl(\min
_{z\in S}v \Bigr)^2.
\end{eqnarray*}
Thus $\{\ts_{\tau,p}(z,\omega)\}$ are equicontinuous in $z$ (with
$p$ and
$\omega$ as parameters) on $S$. By Arzela--Ascoli, $\ts_{\tau
,p}(z,\omega)$
thus converges uniformly to $s_\tau(z)$ on $S$. Consequently
$\ts_{\tau,p}(z,\omega)\to s_\tau(z)$ on $\mathbb{C}^+$. Since
$\tilde F_\tau$,
the ESD of $\tbC_\tau$, is tight (by Lemma~S.12), %\ref{lemexist3}),
it follows from Lemmas S.13 %\ref{lemunique1}
and S.3 %\ref{lemStieltjesconvergence}
that $s_\tau(z)$ is a Stieltjes transform of a (nonrandom) probability
measure, and $\tilde
F_\tau$ converges a.s. to the distribution whose Stieltjes transform
is given
by $s_\tau$. Since, by Lemma S.2, %\ref{lemrankinequality},
$\sup_\sigma\llvert\tilde
F_\tau(\sigma) - \hat F_\tau(\sigma)\rrvert\to0$ a.s., it follows that
$\hat
F_\tau$ converges a.s. to the same limit, and hence $s_{\tau,p}(z)$ converges
a.s. to $s_\tau(z)$. The proof for the Gaussian MA($q$) case is
complete. It
can be checked that all the statements remain valid even if $q=q(p)\to
\infty$
sufficiently slowly under~(\ref{eqpn}), for example, if $q(p)=o(p^{1/2})$.

%s8 #&#
\section{Truncation, centering and rescaling}\label{sectrunc}

The extension of the result to non-Gaussian innovations requires in its first
step, a truncation argument, followed by a centering and rescaling of the
innovations. This section justifies that the symmetrized autocovariance
matrices obtained from Gaussian innovations and from their truncated, centered
and rescaled counterparts have the same LSD. The extension to the non-Gaussian
case is then completed in Section~\ref{secpnon-gauss}.

Since the underlying process is an MA($q$) series, the observations
$X_1,\ldots,X_n$ are functions of the innovations $Z_{1-q},\ldots,Z_n$. For
$j=1,\ldots,p$ and $t=1-q,\ldots,n$, define the quantities
\begin{eqnarray*}
\bar{Z}_{jt}^R&=&Z_{jt}^R\mathbb{I}
\bigl(\bigl\llvert Z_{jt}^R\bigr\rrvert<\sqrt{p}\bigr),
\qquad\check{Z}_{jt}^R=\bar Z_{jt}^R-
\mathbb{E}\bigl[\bar Z_{jt}^R\bigr],
\\
\bar{Z}_{jt}^I&=&Z_{jt}^I\mathbb{I}
\bigl(\bigl\llvert Z_{jt}^I\bigr\rrvert<\sqrt{p}\bigr),
\qquad\check{Z}_{jt}^I=\bar Z_{jt}^R-
\mathbb{E}\bigl[\bar Z_{jt}^I\bigr],
\end{eqnarray*}
where\vspace*{-1pt} $Z_{jt}=Z_{jt}^R+iZ_{jt}^I$ and $\mathbb{I}$ the indicator function.
Correspondingly define $\bar{\bC}_\tau$ and $\check{\bC}_\tau$ to
be the
autocovariance matrices obtained from $\bar{Z}_{jt}^R,\bar{Z}_{jt}^I$ and
$\check{Z}_{jt}^R$, $\check{Z}_{jt}^I$, respectively.

%pr8.1 #&#
\begin{proposition}
If the assumptions of Theorem~\ref{thesdmainfty} are satisfied,
then a.s. under (\ref{eqpn}),
\[
\sup_\lambda\bigl\llvert F^{\bC_\tau}-F^{\bar{\bC}_\tau}\bigr
\rrvert\to0 \quad\mbox{and}\quad\sup_\lambda\bigl\llvert
F^{\bar{\bC}_\tau}-F^{\check{\bC}_\tau}\bigr\rrvert\to0,
\]
where $F^{\bC_\tau}$, $F^{\bar{\bC}_\tau}$ and $F^{\check{\bC
}_\tau}$ denote
the ESDs of $\bC_\tau$, $\bar{\bC}_\tau$ and $\check{\bC}_\tau
$, respectively.
\end{proposition}

\begin{pf}
Let $\eta_{jt}^R=1-\mathbb{I}(\llvert Z_{jt}^R\rrvert>\sqrt{p})$ and
$\eta_{jt}^I=1-\mathbb{I}(\llvert Z_{jt}^I\rrvert>\sqrt{p})$. Let further
$\delta^R=\mathbb{P}(\llvert Z_{jt}^R\rrvert>\sqrt{p})$ and
$\delta^I=\mathbb{P}(\llvert Z_{jt}^I\rrvert>\sqrt{p})$, and note\vspace*{1pt}
that these
quantities are
independent due to the assumed i.i.d. structure on $Z_{jt}$. Since the fourth
moments of the latter random variables are assumed finite, it also
follows that
$\delta^R<p^{-2}\mathbb{E}[Z_{11}^R| ^4]<\infty$ and
$\delta^I<p^{-2}\mathbb{E}[Z_{11}^I| ^4]<\infty$.

Observe next that the rank of a matrix does not exceed the number of its
nonzero columns and that each nonzero $Z_{jt}$ causes at most $2(q+1)$ nonzero
columns. Recalling that $\bC_\tau=(2n)^{-1}\bX(\bL+\bL^*)\bX^*$,
Theorem~A.44
of \cite{bs10} (using $\mathbf{F}=\mathbf{0}$ and $\mathbf{D}=\bL
+\bL^*$)
implies that, for any $\varepsilon>0$,
\begin{eqnarray*}
&& \mathbb{P} \Bigl(\sup_{\lambda}\bigl\llvert F^{\bC_\tau}(
\lambda)-F^{\bar{\bC}_\tau}(\lambda)\bigr\rrvert>\varepsilon\Bigr)
\\
&&\qquad \leq \mathbb{P}
\biggl(\frac{2(q+1)}{p}\sum_{jt}\bigl(
\eta_{jt}^R+\eta_{jt}^I\bigr)>
\varepsilon\biggr)
\\
&&\qquad \leq \mathbb{P} \biggl(\sum_{jt}
\eta_{jt}^R>\frac{\varepsilon
p}{4(q+1)} \biggr) +\mathbb{P} \biggl(
\sum_{jt}\eta_{jt}^I>
\frac{\varepsilon
p}{4(q+1)} \biggr)
\\
&&\qquad =\mathbb{P} \biggl(\sum_{jt}\eta_{jt}^R-p(q+n)
\delta^R>\frac
{\varepsilon p}{4(q+1)}-p(q+n)\delta^R \biggr)
\\
&&\quad\qquad{} +\mathbb{P} \biggl(\sum_{jt}
\eta_{jt}^I-p(q+n)\delta^R>
\frac
{\varepsilon p}{4(q+1)}-p(q+n)\delta^I \biggr).
\end{eqnarray*}
Let $\varepsilon(p,q)=\varepsilon p[4(q+1)]^{-1}$. For $p$ large enough
so that $2p(q+n)\min\{\delta^R,\delta^I\}<\varepsilon(p,q)$,
Hoeffding's inequality
yields
\begin{eqnarray*}
&& \mathbb{P} \biggl(\sum_{jt}\eta_{jt}^R-p(q+n)
\delta^R>\varepsilon(p,q)-p(q+n)\delta^R \biggr)
\\
&&\qquad \leq2\exp\bigl(- \bigl(\varepsilon(p,q)-p(q+n)\delta^R
\bigr)^2 \bigl(p(q+n)\delta^R+\varepsilon(p,q)
\bigr)^{-1} \bigr)
\\
&&\qquad \leq2\exp\biggl(- \biggl(\frac{\varepsilon(p,q)}2 \biggr)^2 \biggl(
\frac{3\varepsilon(p,q)}{2} \biggr)^{-1} \biggr)
\\
&&\qquad =2\exp\biggl(-\frac{\varepsilon p}{24(q+1)} \biggr)
\end{eqnarray*}
as well as
\[
\mathbb{P} \biggl(\sum_{jt}\eta_{jt}^I-p(q+n)
\delta^I>\varepsilon(p,q)-p(q+n)\delta^R \biggr) \leq2\exp
\biggl(-\frac{\varepsilon p}{24(q+1)} \biggr).
\]
The Borel--Cantelli lemma now implies that
$\sup_\sigma\llvert F^{\bC_\tau}(\sigma)-F^{\bar{\bC}_\tau}(\sigma
)\rrvert\to0$ a.s., which
is the first claim of the proposition.

To verify the second, note that the equality
\[
\sum_{\ell=0}^q\bA_\ell
\bar{Z}_{t-\ell}-\sum_{\ell=0}^q\bA
_\ell\check{Z}_{t-\ell} =\sum_{\ell=0}^q
\bA_\ell\bigl(\mathbb{E}\bigl[\bar{Z}_{11}^R
\bigr]+i\mathbb{E}\bigl[\bar{Z}_{11}^I\bigr] \bigr)
\mathbf{1},
\]
with $\mathbf{1}$ being the vector with all entries equal to 1, shows that
$\sum_{\ell=0}^q\bA_\ell\bar{Z}_{t-\ell}-\sum_{\ell=0}^q\bA
_\ell\check{Z}_{t-\ell}$
is independent of $t$. Thus an application of Lemma S.1 %
leads to
\[
\sup_\sigma\bigl\llvert F^{\bar{\bC}_\tau}(\sigma)-F^{\check{\bC}_\tau
}(
\sigma)\bigr\rrvert\leq\frac{1}{p},
\]
which converges to 0 a.s. under (\ref{eqpn}). This is the second assertion.
\end{pf}

After truncation and centering, it does not necessarily follow that\break 
$\mathbb{E}[\llvert\check{Z}_{11}\rrvert^2]=\mathbb{E}[\llvert
\check
{Z}_{11}^R+i\check{Z}_{11}^I\rrvert^2]$
is equal to 1. However, rescaling $\check{Z}_{11}^R+i\check
{Z}_{11}^I$ by
dividing with $(\mathbb{E}[\llvert\check{Z}_{11}^R+i\check
{Z}_{11}^I\rrvert^2])^{1/2}$ (in
order to obtain unit variance) does not affect the LSD because
$\mathbb{E}[\llvert\check{Z}_{11}\rrvert^2]\to1$ under (\ref
{eqpn}). The detailed
arguments follow as in Section~3.1.1 of Bai and Silverstein~\cite
{bs10}. This
shows that the symmetrized autocovariances from~$Z_{jt}$ and their truncated,
centered and rescaled counterparts have the same LSD. Thus to simplify the
argument, it can be assumed that the recentered process has variance one.

%s9 #&#
\section{Extension to the non-Gaussian case}\label{secpnon-gauss}

In this section, the results for the Gaussian MA($q$) case are extended to
general innovation sequences $(W_t\dvtx t\in\mathbb{Z})$ satisfying
the same
moment conditions as their Gaussian counterparts $(Z_t\dvtx t\in
\mathbb{Z})$.
The processes of interest are then the two MA($q$) processes
\[
X_t=\sum_{\ell=0}^q
\bA_\ell Z_{t-\ell}\quad\mbox{and}\quad\oX_t=\sum
_{\ell=0}^q\bA_\ell
W_{t-\ell}.
\]
Define the symmetrized autocovariance matrix
$\obC_\tau=(2n)^{-1}\sum_{t=1}^{n-\tau}(\oX_t\oX_{t+\tau}^*+\oX
_{t+\tau}\oX_t^*)$,
the resolvent $\obR(z)=(\obC_\tau-z\bI)^{-1}$ and the Stieltjes transform
$\bar{s}_p(z)=p^{-1}\operatorname{tr}[\obR(z)]$. In the following, it will
be shown
that the ESDs of $\bC_\tau$ and $\obC_\tau$ converge to the same
limit. This is
done via verifying that, for all $z\in\mathbb{C}^+$, $\bar{s}_p(z)$ and
$s_p(z)$ converge to the same limit $s(z)$ under (\ref{eqpn}), which
in turn
requires us to show that:
\begin{longlist}[(a)]
\item[(a)] $\mathbb{E}[\bar{s}_p(z)-s_p(z)]\to0$ under (\ref
{eqpn}) for all $z\in\mathbb{C}^+$;
\item[(b)] $\mathbb{P}[\llvert\bar{s}_p(z)-\mathbb{E}[\bar
{s}_p(z)]\rrvert\geq
\varepsilon]\to0$
under (\ref{eqpn}) for all $z\in\mathbb{C}^+$ and $\varepsilon>0$.
\end{longlist}
Part (a) requires the use of the Lindeberg principle, and part (b) is
achieved via
an application of McDiarmid's inequality.

%s9.1 #&#
\subsection{Showing that \texorpdfstring{$\mathbb{E}[\bar{s}_p(z)-s_p(z)]\to0$}{E[barsp(z)-sp(z)]to0}}\label{secpnon-gausslindeberg}

For the use in this section, redefine $\bZ=[Z_{1-q}\dvtx\cdots\dvtx
Z_n]$, define
$\bW=[W_{1-q}\dvtx\cdots\dvtx W_n]$ and let $\bZ^R$, $\bW^R$ and $\bZ^I$,
$\bW^I$ be the
corresponding matrices of real and imaginary parts. Claim (a) will be verified
via the Lindeberg principle developed in Chatterjee \cite{c06}. This involves
successive replacements of Gaussian variables with non-Gaussian
counterparts in
a telescoping sum. To this end, define an order on the index set $\{
(j,t)\dvtx
j=1,\ldots,p, t=1-q,\ldots,n\}$ by letting $(j,t)<(j^\prime,t^\prime)$ if
either (1) $t<t^\prime$ or (2) $t=t^\prime$ and $j<j^\prime$, so
that one
successively moves columnwise through the entries of a matrix.
Moreover, let
\[
(j,t)-1= \cases{ (j-1,t), &\quad if $j>1$ and $t\geq1-q$,
\cr
(p,t-1), &\quad if
$j=1$ and $t>1-q$,
\cr
(0,-q), &\quad if $j=1$ and $t=1-q$.}
\]
Let $\bV_{j,t}^R$ denote the $p\times(q+n)$ matrix given by the entries
\[
\bigl(\bV_{jt}^R\bigr)_{j^\prime t^\prime}= \cases{
Z_{j^\prime t^\prime}^R, &\quad if $\bigl(j^\prime,t^\prime
\bigr)\leq(j,t)$,
\vspace*{3pt}\cr
W_{j^\prime t^\prime}^R, &\quad if $
\bigl(j^\prime,t^\prime\bigr)>(j,t)$.}
\]
Further let the $p\times(q+n)$ matrix $\obV_{jt}^R$ be equal to $\bV_{jt}^R$
for all entries but the $(j,t)$th one, which is set to equal 0, and define
analogously the matrices $\bV_{jt}^I$ and $\obV_{jt}^I$. These matrices
determine how many of the original Gaussian $Z_{jt}$'s have been
replaced by the
non-Gaussian $W_{jt}$. In the following, $s_p(z)$ will be viewed as a function
of $\bZ^R$ and $\bZ^I$, fixing $z$ and $p$ as parameters, that is,
$s_p^z(\bZ^R,\bZ^I)
%&
=s_p(z)$. %\\
%&=&\frac1p\operatorname{tr}\bigg[\bigg(\frac1{2n}\bigg\{\sum_{\ell=0}^q
%where $\bF=[\mathbf{0}_{n\times q}:\bI_n]^\prime$ is an $(n+q)\times
%n$ matrix.
Similarly, let $s_p^z(\bW^R,\bW^I)
=\bar{s}_p(z)$. Utilizing this notation, the quantity to be bounded in
expectation can be written as
\begin{eqnarray*}
&& s_p^z\bigl(\bZ^R,\bZ^I
\bigr)-s_p^z\bigl(\bW^R,\bW^I
\bigr)
\\
&&\qquad = s_p^z\bigl(\bV_{pn}^R,
\bV_{pn}^I\bigr)-s_p^z\bigl(
\bV_{0,-q}^R,\bV_{0,-q}^I\bigr)
\\
&&\qquad = \sum_{(j,t)=(1,1-q)}^{(p,n)} \bigl[s_p^z
\bigl(\bV_{jt}^R,\bV_{pn}^I
\bigr)-s_p^z\bigl(\bV_{(j,t)-1}^R,
\bV_{pn}^I\bigr) \bigr]
\\
&&\quad\qquad{} +\sum_{(j,t)=(1,1-q)}^{(p,n)}
\bigl[s_p^z\bigl(\bV_{0,-q}^R,\bV
_{jt}^I\bigr)-s_p^z\bigl(
\bV_{0,-q}^R,\bV_{(j,t)-1}^I\bigr) \bigr]
\\
&&\qquad =\Delta_1+\Delta_2,
\end{eqnarray*}
where $\Delta_1=\sum_{(j,t)=(1,1-q)}^{(p,n)}\Delta_{j,t}^{(1)}$ and
$\Delta_2=\sum_{(j,t)=(1,1-q)}^{(p,n)}\Delta_{j,t}^{(2)}$ are the
real and
imaginary part of the difference with
$\Delta_{j,t}^{(1)}=s_p^z(\bV_{jt}^R,\bV_{pn}^I)-s_p^z(\bV
_{(j,t)-1}^R,\bV_{pn}^I)$
and
$\Delta_{j,t}^{(2)}=s_p^z(\bV_{0,-q}^R,\bV_{jt}^I)-s_p^z(\bV
_{0,-q}^R,\bV_{(j,t)-1}^I)$.
In the\vspace*{2pt} following only the telescoping real parts will be discussed, as the
imaginary parts can be estimated along the same lines. Inserting $\obV_{jt}^R$,
one obtains
\begin{eqnarray*}
\Delta_{j,t}^{(1)} &=& \bigl[s_p^z
\bigl(\bV_{j,t}^R,\bV_{p,n}^I
\bigr)-s_p^z\bigl(\obV_{j,t}^R,\bV
_{p,n}^I\bigr) \bigr] +\bigl[s_p^z
\bigl(\obV_{j,t}^R,\bV_{p,n}^I
\bigr)-s_p^z\bigl(\bV_{(j,t)-1}^R,\bV
_{p,n}^I\bigr) \bigr]
\\
&=&\Delta_{j,t}^{(1,1)}+\Delta_{j,t}^{(1,2)}.
\end{eqnarray*}
Let $\partial_{j,t,1}^{(k)}s_p^z$ denote the $k$th-order partial
derivative of
$s_p^z$ with respect to $Z_{jt}^R$. A~Taylor series expansion gives
\begin{eqnarray*}
&& \biggl\llvert\Delta_{j,t}^{(1,1)}-Z_{jt}^R
\partial_{j,t,1}^{(1)}s_p^z\bigl(\obV
^R_{jt},\bV_{pn}^I\bigr) -
\frac{1}2\bigl(Z_{jt}^R\bigr)^2
\partial_{j,t,1}^{(2)}s_p^z\bigl(
\obV^R_{jt},\bV_{pn}^I\bigr)\biggr
\rrvert
\\
&&\qquad \leq\frac{1}6\bigl\llvert Z_{jt}^R\bigr\rrvert
^3 \max_{\alpha\in[0,1]}\bigl\llvert\partial_{j,t,1}^{(3)}s_p^z
\bigl(\alpha\bV_{jt}^R+(1-\alpha)\obV_{jt}^R,
\bV_{pn}^I\bigr)\bigr\rrvert
\end{eqnarray*}
and
\begin{eqnarray*}
&& \biggl\llvert\Delta_{j,t}^{(1,2)}+W_{jt}^R
\partial_{j,t,1}^{(1)}s_p^z\bigl(\obV
^R_{jt},\bV_{pn}^I\bigr) +
\frac{1}2\bigl(W_{jt}^R\bigr)^2
\partial_{j,t,1}^{(2)}s_p^z\bigl(
\obV^R_{jt},\bV_{pn}^I\bigr)\biggr
\rrvert
\\
&&\qquad \leq\frac{1}6\bigl\llvert W_{jt}^R\bigr\rrvert
^3 \max_{\alpha\in[0,1]}\bigl\llvert\partial_{j,t,1}^{(3)}s_p^z
\bigl(\alpha\bV_{(j,t)-1}^R+(1-\alpha)\obV_{jt}^R,
\bV_{pn}^I\bigr)\bigr\rrvert.
\end{eqnarray*}
The entries of the matrices $\bZ^R$, $\bZ^I$, $\bW^R$ and $\bW^I$
are all
independent of each other and the first and second moments of the
various real
parts (and imaginary parts) coincide, so that the bound in the last two
inequalities also hold for higher-order terms (HOT). This leads to
%
%e9.1 #&#
\begin{eqnarray}\label{eqtaylor}
\bigl\llvert\mathbb{E} \bigl[\Delta_{j,t}^{(1)} \bigr]\bigr
\rrvert&=&\biggl\llvert\mathbb{E} \biggl[ \bigl(Z_{jt}^R-W_{jt}^R
\bigr)\partial_{j,t,1}^{(1)}s_p^z
\bigl(\obV^R_{jt},\bV_{pn}^I\bigr) \nonumber
\\
&&\hspace*{15pt}{}+
\frac{1}2 \bigl[ \bigl(Z_{jt}^R
\bigr)^2-\bigl(W_{jt}^R \bigr)^2
\bigr]\partial_{j,t,1}^{(2)}s_p^z\bigl(
\obV^R_{jt},\bV_{pn}^I\bigr)\pm
\mathrm{HOT} \biggr]\biggr\rrvert
\nonumber\\[-8pt]\\[-8pt]\nonumber
&\leq&\frac{1}6\mathbb{E} \Bigl[\bigl\llvert Z_{jt}^R
\bigr\rrvert^3 \max_{\alpha\in[0,1]}\bigl\llvert
\partial_{j,t,1}^{(3)}s_p^z\bigl(\alpha
\bV_{jt}^R+(1-\alpha)\obV_{jt}^R,
\bV_{pn}^I\bigr)\bigr\rrvert\Bigr]
\\
&&{} + \frac{1}6\mathbb{E} \Bigl[\bigl\llvert
W_{jt}^R\bigr\rrvert^3 \max
_{\alpha\in[0,1]}\bigl\llvert\partial_{j,t,1}^{(3)}s_p^z
\bigl(\alpha\bV_{(j,t)-1}^R+(1-\alpha)\obV_{jt}^R,
\bV_{pn}^I\bigr)\bigr\rrvert\Bigr].\nonumber\hspace*{-25pt}
\end{eqnarray}
Dealing\vspace*{1pt} with the right-hand side of (\ref{eqtaylor}) requires the computation
and estimation of the third-order derivatives
$\partial_{j,t,1}^{(3)}s_p^z(\bZ^R,\bZ^I)=p^{-1}\operatorname
{tr}[\partial_{j,t,1}^{(3)}\bR(z)]$.
Focusing only on the first term of the right-hand side of (\ref{eqtaylor})
(the second can be handled similarly), Lemma S.14 %
shows that this term converges to zero under (\ref{eqpn}) if, almost
surely under
(\ref{eqpn}),
%
%e9.2 #&#
%e9.3 #&#
\begin{eqnarray}
&& \frac{q+1}{n^3p}\sum_{(j,t)=(1,1-q)}^{(p,n)}
\mathbb{E} \biggl[\bigl\llvert Z_{t,j}^R\bigr\rrvert
^3\max_{\alpha\in[0,1]} \biggl(\sum
_{\ell\in\mathcal{I}_+(t)}\llVert\bA_\ell\rrVert\bigl\llVert
\tX_{t+\tau
+\ell}^\alpha\bigr\rrVert
\nonumber\\[-8pt] \label{eqremainder1} \\[-8pt] \nonumber
&&\hspace*{165pt}{} +\sum_{\ell\in\mathcal{I}_-(t)}\bigl\llVert\tX_{t-\tau+\ell}^\alpha
\bigr\rrVert\llVert\bA_\ell\rrVert\biggr)^3 \biggr]\to 0,\hspace*{-20pt}
\nonumber
\\
&& \frac{q+1}{n^2p}\sum_{(j,t)=(1,1-q)}^{(p,n)}
\mathbb{E} \biggl[\bigl\llvert Z_{t,j}^R\bigr\rrvert
^3\max_{\alpha\in[0,1]} \biggl(\sum
_{\ell\in\mathcal{I}_+(t)}\llVert\bA_\ell\rrVert\bigl\llVert
\tX_{t+\tau+\ell}^\alpha\bigr\rrVert
\nonumber\\[-8pt] \label{eqremainder2} \\[-8pt]\nonumber
&&\hspace*{165pt}{} +\sum_{\ell\in\mathcal{I}_-(t)}\bigl\llVert\tX_{t-\tau+\ell}^\alpha
\bigr\rrVert\llVert\bA_\ell\rrVert\biggr) \biggr]\to0,\hspace*{-20pt}
\nonumber
\end{eqnarray}
where
$\mathcal{I}_+(t)=\{\ell\dvtx\max\{0,1-t-\tau\}\leq\ell\leq
\min\{q,n-t-\tau\}\}$,
$\mathcal{I}_-(t)=\{\ell\dvtx\break  \max\{0,1-t+\tau\}\leq\ell\leq
\min\{q,n-t\}\}$,
$\tbX^\alpha=[\tX_1^\alpha\dvtx\cdots\dvtx\tX_n^\alpha]=\sum_{\ell
=0}^q\bA_\ell\bV^\alpha\bL^\ell\bF$,
$\bV^\alpha=\alpha\bV_{jt}^R+(1-\alpha)\obV_{jt}^R+i\bZ^I$ and
$\bF=[\mathbf{0}_{n\times q}\dvtx\bI_n]^\prime$ being\vspace*{1pt} an $(n+q)\times
n$ matrix.
The choice of $\alpha\in[0,1]$ only affects the value of the $j$th
entry of
$V_t^\alpha$. For $t\neq t^\prime$, the notation
$V_{t^\prime}=V_{t^\prime}^\alpha$ is therefore preferred. By definition,
$V_t^0$ is the vector whose $j$th entry has a real part of zero. Let
$\mathcal{J}_{\pm}(t,\ell)=\{\ell^\prime\in\{0,\ldots,q\}\dvtx
\ell^\prime\neq\tau\pm\ell,t\pm\tau+\ell-\ell^\prime\geq
1\}$. Then
\begin{eqnarray*}
&& \max_{\alpha\in[0,1]}\bigl\llVert\tX^\alpha_{t\pm\tau+\ell}
\bigr\rrVert
\\
&&\qquad \leq\sum_{\ell^\prime\in\mathcal{J}_\pm(t,\ell)}\bar
{\lambda
}_{\bA_{\ell^\prime}}\llVert Z_{t\pm\tau+\ell-\ell^\prime}\rrVert+\bar
{\lambda}_{\bA_{\ell\pm\tau}}
\bigl(\bigl\llVert V_t^0\bigr\rrVert+\bigl\llvert
Z_{j,t}^R\bigr\rrvert\bigr).
\end{eqnarray*}
Hence, setting $\bar{\lambda}_\ell^{\mathrm{rem}}=(\sum_{\ell^\prime
=\llvert\ell\rrvert}^q\bar{\lambda}_{\bA
_{\ell
^\prime}}^2)^{1/2}$,
\begin{eqnarray*}
&& \max_{\alpha\in[0,1]}\sum_{\ell\in\mathcal{J}_\pm(t)} \llVert
\bA_\ell\rrVert\bigl\llVert\tX^\alpha_{t\pm\tau+\ell}\bigr\rrVert
\\
&&\qquad \leq\bar{\lambda}_\bA\sum_{\ell=-q,\ell\neq\mp\tau}^q
\bar{\lambda}_\ell^{\mathrm{rem}}\llVert Z_{t\pm\tau+\ell}\rrVert+
\bar{\lambda}_\bA\bar{\lambda}_{\pm\tau}^{\mathrm{rem}} \bigl(
\bigl\llVert V_t^0\bigr\rrVert+\bigl\llvert
Z_{j,t}^R\bigr\rrvert\bigr).
\end{eqnarray*}
Using this, it follows from Lemma S.15 %\ref{lemnon-gaussremainder}
that the left-hand sides of (\ref{eqremainder1}) and (\ref
{eqremainder2}) converge
to zero a.s., thus establishing that $\mathbb{E}[\bar
{s}_p(z)-s_p(z)]\to0$.

%s9.2 #&#
\subsection{Showing that \texorpdfstring{$\mathbb{P}[|\bar{s}_p(z)-\mathbb{E}[\bar{s}_p(z)]|\geq\varepsilon]\to0$}
{P[|barsp(z)-E[barsp(z)]|geqvarepsilon]to0}}\label{secpnon-gauss-mcdiarmid}

For a fixed $z$, $\bar{s}_p(z)$ is a function of the $n+q$ vectors
$W_{1-q},\ldots,W_n$. Letting $m=\lceil q^{-1}(n+q)\rceil$, these are now
segmented into the groups $W_{(k-1)q+1},\ldots,W_{kq}$, $k=0,1,\ldots,m$,
possibly adding additional vectors to the last group to ensure all
groups have
the same length [even though the value $\bar{s}_p(z)$ does not depend
on the
additions]. To satisfy the conditions needed in order to apply McDiarmid's
inequality, note that a change of the values $W_{(k-1)q+1},\ldots
,W_{kq}$ in
one group to, say, $W^\prime_{(k-1)q+1},\ldots,W^\prime_{kq}$,
causes the
values of $\bar{X}_{(k-1)q+1},\ldots,\bar{X}_{(k+1)q}$ to change to, say,
$\bar{X}^\prime_{(k-1)q+1},\ldots,\bar{X}^\prime_{(k+1)q}$. In the
following,
the focus is on changes applied to the first group of innovations
$W_1,\ldots,W_q$. Consider the case $\tau\leq2q$ and let
\[
\bar{\bC}_\tau^\prime=\frac{1}{2n}\sum
_{t=1}^{2q} \bigl(\bar{X}_t^\prime
\bar{X}_{t+\tau
}^*+\bar{X}_{t+\tau} {\bar{X}}_t^{\prime,*}
\bigr) +\frac{1}{2n}\sum_{t=2q+1}^{n-\tau}
\bigl(\bar{X}_t\bar{X}_{t+\tau
}^*+\bar{X}_{t+\tau}
\bar{X}_t^* \bigr),
\]
$\bar{\bR}_\tau^\prime(z)=(\bar{\bC}_\tau^\prime-z\bI)^{-1}$ and
$\bar{s}_p^\prime(z)=p^{-1}\operatorname{tr}[\bar{\bR}_\tau^\prime
(z)]$. The goal is
now to represent $\bar{\bC}_\tau^\prime$ as a finite rank
perturbation of
$\bar{\bC}_\tau$ in the form $\sum_{j=1}^Jr_jr_j^*$ with
appropriate $r_j$ and~$J$. Write
\begin{eqnarray*}
\bar{\bC}_\tau^\prime-\bar{\bC}_\tau%%
&=&
\frac{1}{2n}\sum_{t=1}^{2q} \bigl(
\bar{X}_t^\prime\bar{X}_{t+\tau
}^*+
\bar{X}_{t+\tau} \bar{X}_t^{\prime,*} -
\bar{X}_t\bar{X}_{t+\tau}^*-\bar{X}_{t+\tau}
\bar{X}_t^* \bigr)
\\
&=&\frac{1}{2n}\sum_{t=1}^{2q} \bigl[
\bigl(\bar{X}_t^\prime+\bar{X}_{t+\tau}\bigr) \bigl(
\bar{X}_t^\prime+\bar{X}_{t+\tau}\bigr)^*
\\
&&\hspace*{31pt}{} -(\bar{X}_t+\bar{X}_{t+\tau}) (
\bar{X}_t+\bar{X}_{t+\tau})^*-\bar{X}_t^\prime
\bar{X}_t^{\prime,*}+\bar{X}_t \bar{X}_t^*
\bigr].
\end{eqnarray*}
Choosing $J=8q$ and repeatedly utilizing (S.3) %\eqref{eqdiffresolvent}
with $\mathcal{H}_t$ replaced by $\bI_p$, it follows that
$\llvert\bar{s}_p(z)-\bar{s}_p^\prime(z)\rrvert\leq C_1q(vp)^{-1}$
for some appropriately
chosen constant $C_1>0$. This bound holds for any of the $m$ groups of
innovations. McDiarmid's inequality consequently implies, for any
$\varepsilon>0$
and a suitable constant $C_2>0$,
\[
\mathbb{P} \bigl(\bigl\llvert\bar{s}_p(z)-\mathbb{E}\bigl[
\bar{s}_p(z)\bigr]\bigr\rrvert\geq\varepsilon\bigr) \leq4\exp
\biggl(-C_2\frac{\varepsilon^2v^2p^2}{mq^2} \biggr) \sim4\exp\biggl(-C_2
\frac{\varepsilon^2v^2c}{q}p \biggr).
\]
The right-hand side converges to zero at a rate faster than $qp^{-1}$ and
concentration of the Stieltjes transform around its mean is
established, since
the case $\tau>2q$ can be handled in a similar fashion. Note that the last
argument remains valid if $q=q(p)\to\infty$ as $p\to\infty$ at a
sufficiently
slow rate, for example, if $q=o(p^{1/2})$.

%s10 #&#
\section{Extension to the linear process case}\label{secplin}

Let $(X_t\dvtx t\in\mathbb{Z})$ now denote a linear process. To
complete the
proof of Theorem~\ref{thesdmainfty}, a truncation argument is
invoked. Let
$(X_t^{\mathrm{tr}}\dvtx\mathbb{Z})$ denote the truncated process
given by
$X_t^{\mathrm{tr}}=\sum_{\ell=0}^{q(p)}\bA_\ell Z_{t-\ell}$, $t\in
\mathbb{Z}$, where
$q(p)$ depends on the dimension. Let $(\check{X}_t\dvtx
t\in\mathbb{Z})$ be the process given by $\check{X}_t=X_t-X_t^{\mathrm{tr}}$,
$t\in\mathbb{Z}$, and denote by $L(F,G)$ the L\'evy distance between
distribution functions $F$ and $G$.

%le10.1 #&#
\begin{lemma}\label{lemlinproclevy}
If the assumptions of Theorem~\ref{thesdmainfty} are satisfied and if
$q(p)=\lceil p\rceil^{1/3}$, then
\[
L\bigl(F^{\bC_\tau},F^{\bC_\tau^{\mathrm{tr}}}\bigr)\to0 \qquad\mbox{a.s.}
\]
under (\ref{eqpn}), where $\bC_\tau^{\mathrm{tr}}=(2n)^{-1}\sum_{\tau
=1}^{n-\tau} (X_t^{\mathrm{tr}}(X_{t+\tau
}^{\mathrm{tr}})^*+X_{t+\tau}^{\mathrm{tr}}(X_t^{\mathrm{tr}})^* )$.
\end{lemma}

\begin{pf}
By Lemma S.1, %\ref{lemnorminequality},
it suffices to show that
$p^{-1} \operatorname{tr}[(\tbC_\tau- \tbC_\tau^{\mathrm{tr}})^2] \to0$ a.s.
Write $\bX=
\bX^{\mathrm{tr}} + \check\bX$ and $\bD_\tau= (\bL_\tau+ \bL_\tau
^*)/2$. Then
\[
\tbC_\tau- \tbC_\tau^{\mathrm{tr}} = \frac{1}{n}
\bX^{\mathrm{tr}} \bD_\tau\check\bX^* + \frac{1}{n} \check\bX
\bD_\tau\bigl(\bX^{\mathrm{tr}}\bigr)^* + \frac{1}{n} \check\bX
\bD_\tau\check\bX^* = \bP+ \bP^* + \bQ,
\]
say. From repeated applications of the Cauchy--Schwarz inequality, we have
\begin{eqnarray*}
\operatorname{tr}\bigl[\bigl(\tbC_\tau- \tbC_\tau^{\mathrm{tr}}
\bigr)^2\bigr] &=& 2\operatorname{tr}\bigl(\bP^2\bigr) + \operatorname{tr}
\bigl(\bQ^2\bigr) + 2 \operatorname{tr}\bigl(\bP\bP^*\bigr) + 4 \operatorname{tr}(\bP
\bQ)
\\
& \leq & 4 \operatorname{tr}\bigl(\bP\bP^*\bigr) + \operatorname{tr}\bigl(\bQ^2\bigr) +
4 \sqrt{\operatorname{tr}\bigl(\bP\bP^*\bigr)} \sqrt{\operatorname{tr}\bigl(\bQ^2
\bigr)}.
\end{eqnarray*}
Since $\llVert\bD_\tau\rrVert\leq1$, $\operatorname{tr}(\bQ^2) \leq
n^{-2}\operatorname{tr}[(\check
\bX\check\bX^*)^2]$, and by another application of Cauchy--Schwarz,
\[
\operatorname{tr}\bigl(\bP\bP^*\bigr) \leq\biggl(\frac{1}{n^2} \operatorname{tr}\bigl[
\bigl(\bX^{\mathrm{tr}}\bigl(\bX^{\mathrm{tr}}\bigr)^*\bigr)^2\bigr]
\biggr)^{1/2} \biggl(\frac{1}{n^2}\operatorname{tr}\bigl[\bigl(\check\bX
\check\bX^*\bigr)^2\bigr] \biggr)^{1/2}.
\]
Since it is easy to see that $(pn^2)^{-1}\operatorname{tr}[(\bX^{\mathrm
{tr}}(\bX
^{\mathrm{tr}})^*)^2]$ is stochastically bounded (e.g., by showing that the
expectation is
finite), it is enough to show that $(pn^2)^{-1} \operatorname{tr}[(\check\bX
\check
\bX^*)^2] \to0$ a.s. This is established by showing that the sum
$\sum_{p=1}^\infty
(pn^2)^{-1}\mathbb{E} [\operatorname{tr}[(\check\bX\check\bX^*)^2]
] <
\infty$, and then applying the Borel--Cantelli lem\-ma. To this end,
note that
\begin{eqnarray*}
&& \mathbb{E} \bigl[\operatorname{tr}\bigl[\bigl(\check\bX\check\bX^*\bigr)^2
\bigr] \bigr]
\\
&&\qquad = \sum_{t=1}^n \sum
_{s=1}^n \mathbb{E} \bigl[\bigl\llvert\check
X_t^* \check X_s\bigr\rrvert^2 \bigr]
\\
&&\qquad  = \sum_{t=1}^n \sum
_{s=1}^n \sum_{\ell=q+1}^\infty
\sum_{\ell'=q+1}^\infty\sum
_{m=q+1}^\infty\sum_{m'=q+1}^\infty
\mathbb{E} \bigl[\operatorname{tr}\bigl(Z_{t-\ell'}Z_{t-\ell}^* \bA
_\ell\bA_m
\\
&&\hspace*{204pt}{}\times Z_{s-m}Z_{s-m'}
\bA_{\ell'} \bA_{m'}\bigr) \bigr].
\end{eqnarray*}
It is clear from the independence of $Z_t$'s that the summands are nonzero
only if the indices of $Z$'s pair up. Direct calculation shows that the total
contribution of all four types of pairings: (i) $t-\ell= t-\ell' \neq
s - m =
s - m'$; (ii) $t-\ell= s-m \neq s - m' = t - \ell'$; (iii) $t-\ell=
s-m' \neq
s - m = t - \ell'$ and (iv) $t-\ell= t-\ell' = s - m = s-m'$ can be
bounded by
$C(p^2n^2)(\sum_{\ell=q}^\infty\bar\lambda_{\bA_\ell}^2)^2$ for
some $C > 0$,
using the fact that $p$ and $n$ are of the same order.
Thus since $q=q(p) = \lceil p\rceil^{1/3}$,
\begin{eqnarray*}
\sum_{p=1}^\infty p \Biggl(\sum
_{\ell=q}^\infty\bar\lambda_{\bA_\ell}^2
\Biggr)^2 &\leq& 2 \sum_{p=1}^\infty
p \sum_{\ell=q}^\infty\biggl(\bar
\lambda_{\bA
_\ell}^2 \sum_{j \geq\ell} \bar
\lambda_{\bA_j}^2\biggr)
\\
&=& 2\sum_{\ell=1}^\infty\biggl(\bar
\lambda_{\bA_\ell}^2 \sum_{j\geq
\ell} \bar
\lambda_{\bA_j}^2 \biggr) \sum_{p=1}^\infty
p \mathbf{1}\bigl\{p\dvtx \lceil p\rceil^{1/3} \leq\ell\bigr\}
\\
&\leq& 2\sum
_{\ell=1}^\infty\biggl(\bar\lambda_{\bA_\ell}^2
\sum_{j\geq\ell
} \bar\lambda_{\bA_j}^2
\biggr) \sum_{p=1}^{\ell^3} p
\\
&\leq& 2 \sum_{\ell=1}^\infty\ell^4
\bar\lambda_{\bA_\ell}^2 \sum_{j\geq\ell}
\bar\lambda_{\bA_j}^2 \leq2 \sum
_{\ell
=1}^\infty\ell^2 \bar
\lambda_{\bA_\ell}^2 \sum_{j\geq\ell}
j^2 \bar\lambda_{\bA_j}^2
\\
&\leq& 2 \Biggl(\sum
_{\ell=1}^\infty\ell\bar\lambda_{\bA_\ell}
\Biggr)^4.
\end{eqnarray*}
This proves the result.
\end{pf}

Using Gaussian innovations $(Z_t\dvtx t\in\mathbb{Z})$, let $\tbX
^{\mathrm{tr}}=[X_1^{\mathrm{tr}}\dvtx\cdots\dvtx X_n^{\mathrm
{tr}}]=\break  \sum_{\ell=0}^q\bA_\ell\bZ\tbL^\ell\bU_{\tbL}$ be the
corresponding
transformed data matrix. Then define %the truncated quantities
$\psi^{\mathrm{tr}}(\lambda,\nu)=\sum_{\ell=0}^{q(p)}e^{i\ell\nu}f_\ell
(\lambda
)$, $\psi^{\mathrm{tr}}(\bA,\nu)=\sum_{\ell=0}^{q(p)}e^{i\ell\nu}\bA
_\ell$, $\tbC
_\tau^{\mathrm{tr}}=\gamma_{\tau,t}\tX_t^{\mathrm{tr}}(\tX_t^{\mathrm
{tr}})^*$,\break  $\tbR_\tau
^{\mathrm{tr}}=(\tbC_\tau^{\mathrm{tr}}-z\bI)^{-1}$, $\tilde{s}^{\mathrm
{tr}}(z)=p^{-1}\operatorname{tr}[\tbR_\tau^{\mathrm{tr}}(z)]$, $\tK
_p^{\mathrm{tr}}(z,\nu)=p^{-1}\operatorname{tr}[\tbR_\tau^{\mathrm
{tr}}(z)\psi^{\mathrm{tr}}(\bA,\nu)]$,
$h^{\mathrm{tr}}(\lambda,\nu)=\llvert\psi^{\mathrm{tr}}(\lambda,\nu
)\rrvert
^2$ and $\cH
^{\mathrm{tr}}(\bA,\nu)=\psi^{\mathrm{tr}}(\bA,\nu)\psi^{\mathrm
{tr}}(\bA,\nu)^*$.
Then one verifies in a similar vein,
as in the proofs of Theorems~\ref{thexist} and~\ref{thunique}, that
for all
$\omega\in\Omega_0$ with $\Omega_0$ defined in the beginning of
Section~\ref{secgaussexconcon}, under (\ref{eqpn}),
\[
\tK_{p_\ell}^{\mathrm{tr}}(z,\nu,\omega) -\int\bigl[U_{p_\ell}^{\mathrm{tr}}(z,
\lambda)-z\bigr]^{-1}h^{\mathrm{tr}}(\lambda,\nu)\,dF_{p_\ell}^\bA(
\lambda) \to0,
\]
where $U_{p_\ell}^{\mathrm{tr}}(z,\lambda)=(n(p_\ell))^{-1}
\sum_{t=1}^{n(p_\ell)}[1+c_{p_\ell}\gamma_{\tau,t}\tK_{p_\ell
}^{\mathrm{tr}}(z,\nu_{tj},\omega)]^{-1}h^{\mathrm{tr}}(\lambda,\nu
_{tj})$. This is
done by
exploiting the convergence $h^{\mathrm{tr}}(\lambda,\nu)\to h(\lambda
,\nu
)$, which
is uniform in $\nu$ and $\lambda$, and $U_{p_\ell}^{\mathrm
{tr}}(z,\lambda)\to
U(z,\lambda,\omega)=(2\pi)^{-1}\int_0^{2\pi}\beta_{\tau
}(z,\lambda)\gamma_\tau(\nu)h(\lambda,\nu)\,d\nu$,
which is uniform in $\lambda$. Therein,
$\beta_{\tau}(z,\nu)=[1+c\gamma_{\tau}(\nu)K_{\tau}(z,\nu
)]^{-1}$. From these
facts it follows that the limiting version of the truncated version satisfies
the defining equations for the Stieltjes kernel (\ref{eqktau}).
Therefore the
results for the complex Gaussian innovation model with fixed order $q$ are,
subject to minor modifications, still applicable when orders $q(p)$
grow at a suitable rate.

%s11 #&#
\section{The real-valued case}\label{secpreal}

The idea of the proof of Theorem~\ref{coesdmainfty} is motivated by
focusing on the $\mathrm{MA}(1)$ case. The
derivation for the $\mathrm{MA}(q)$ and finally $\mathrm{MA}(\infty)$ cases follows the
from the
corresponding transformations and subsequent constructions analogous to the
complex case. So, let $\bX=\bZ+\bA_1\bZ\bL$ be the data matrix
obtained from a
Gaussian $\mathrm{MA}(1)$ time series, and suppose that $\bA_1$ possesses an
eigendecomposition $\bA_1=\bU_{\bA_1}\bLambda_{\bA_1}\bU^*_{\bA
_1}$ with
$\bU_{\bA_1}$ orthogonal.\vspace*{1pt} Let $\bU$ denote the real Fourier basis
(see,\vspace*{1pt} e.g., Chapter~10 of \cite{bd92}), and let for the real case
$\tbX=\bU_{\bA_1}(\bZ+\bA_1\bZ\tbL)\bU$. Since $\bZ\bU$ and
$\bZ\tbL\bU$ have
independent columns, it follows that $\tbX$ has independent columns. Moreover,
$\tbX$ has also independent rows. To see this, note that the transpose
of the
$j$th row of $\tbX$ is
\[
\bigl(\bU_{\bA_1}(\bZ+\bA_1\bZ\tbL)\bU
\bigr)^Te_j =\bU^T(\bI+\lambda_j
\tbL)Z_j,
\]
where $Z_j$ is the $j$th column of $\bU_{\bA_1}\bZ\bU$ and $\lambda
_j$ the $j$th
eigenvalue of $\bA_1$. The covariance of the $j$th column is
\[
\mathbb{E} \bigl[\bU^T(\bI+\lambda_j
\tbL)Z_jZ_j^T(\bI+\lambda_j
\tbL)^T\bU\bigr] =\bU^T \bigl(\bI+\lambda_j
\bigl(\tbL+\tbL^T\bigr)+\lambda_j^2\tbL\tbL
^T \bigr)\bU.
\]
Since $\bI+\lambda_j(\tbL+\tbL^T)+\lambda_j^2\tbL\tbL^T =
(1+\lambda_j^2)\bI+
\lambda_j(\tbL+\tbL^T)$ is a symmetric circulant matrix, it
diagonalizes in the
real Fourier basis, and hence the covariance matrix of the last display is
diagonal. From the same display it follows also that the variance of the
$(j,t)$th entry is $h(\lambda_j,\nu_t)$, so that the rest of the
proof follows
as in the complex case (Theorem~\ref{thesdmainfty}).


\begin{supplement}[id=suppA]
\stitle{Supplement to ``On the Mar\v{c}enko--Pastur law for linear time series''}
\slink[doi]{10.1214/14-AOS1294SUPP} %[doi,text={...}] - jei reikia
%suskaldyti doi
\sdatatype{.pdf}
\sfilename{aos1294\_supp.pdf}
\sdescription{The supplementary material provides additional technical lemmas and their proofs.}
\end{supplement}

% imsref loaded by linak, 2015-01-29 10:49:50
%

\printaddresses

\begin{thebibliography}{37}
% pybtex-1.21. Style name=ims, version=2.91, label_style=nolabel,
%sorting_style=complex, cfg=None, language=None.
%b1 ###
%b1 #&#
\bibitem{aw10}
%
\begin{barticle}[auto:parserefs-M02]
\bauthor{\bsnm{Aguilar},~\bfnm{O.}\binits{O.}} \AND
\bauthor{\bsnm{West},~\bfnm{M.}\binits{M.}}
(\byear{2000}).
\btitle{Bayesian dynamic factor models and portfolio allocation}.
\bjournal{J. Bus. Econom. Statist.}
\bvolume{18}
\bpages{338--357}.
\end{barticle}
%
%\OrigBibText
%Aguilar, O. and West, M. (2000). Bayesian dynamic factor models and
%portfolio allocation.
%{\it Journal of Business \& Economic Statistics} {\bf18}, 338--357.
%\endOrigBibText
\bptok{imsref}%
\endbibitem

%b2 ###
%b2 #&#
\bibitem{agz09}
%
\begin{bbook}[mr]
\bauthor{\bsnm{Anderson},~\bfnm{Greg~W.}\binits{G.~W.}},
\bauthor{\bsnm{Guionnet},~\bfnm{Alice}\binits{A.}} \AND
\bauthor{\bsnm{Zeitouni},~\bfnm{Ofer}\binits{O.}}
(\byear{2010}).
\btitle{An Introduction to Random Matrices}.
\bpublisher{Cambridge Univ. Press},
\blocation{Cambridge}.
\bid{mr={2760897}}
\bptnote{check year}%
\end{bbook}
%%
%\OrigBibText
%Anderson, G. Guionnet, A. and Zeitouni, O. (2009). {\it An
%Introduction to
%Random Matrices}. Cambridge University Press.
%\endOrigBibText
\bptok{imsref}%
% NOT OUTPUTTED:
% isbn = 978-0-521-19452-5
% fpage = xiv+492
\endbibitem

%b3 ###
%b3 #&#
\bibitem{abr10}
%
\begin{barticle}[auto:parserefs-M02]
\bauthor{\bsnm{Angelini},~\bfnm{M.}\binits{M.}},
\bauthor{\bsnm{Ba{\'{n}}bura},~\bfnm{M.}\binits{M.}} \AND
\bauthor{\bsnm{R{\"{u}}nstler},~\bfnm{G.}\binits{G.}}
(\byear{2010}).
\btitle{Estimating and forecasting the Euro area monthly national
accounts from a dynamic factor model}.
\bjournal{J. Bus. Cycle Meas. Anal.}
\bvolume{2010}
\bpages{1--22}.
\end{barticle}
%
%\OrigBibText
%Angelini, M., Ba\'{n}bura, M. and R\"{u}nstler, G. (2010).
%Estimating and forecasting the Euro
%area monthly national accounts from a dynamic factor model.
%{\it Journal of Business Cycle Measurement and Analysis} {\bf2010}, 1--22.
%\endOrigBibText
\bptok{imsref}%
\endbibitem

%b4 ###
%b4 #&#
\bibitem{bn07}
%
\begin{barticle}[mr]
\bauthor{\bsnm{Bai},~\bfnm{Jushan}\binits{J.}} \AND
\bauthor{\bsnm{Ng},~\bfnm{Serena}\binits{S.}}
(\byear{2007}).
\btitle{Determining the number of primitive shocks in factor models}.
\bjournal{J.~Bus. Econom. Statist.}
\bvolume{25}
\bpages{52--60}.
\bid{doi={10.1198/073500106000000413}, issn={0735-0015}, mr={2338870}}
\end{barticle}
%
%\OrigBibText
%Bai, J. and Ng, S. (2007). Determining the number of primitive shocks
%in factor models.
%{\it Journal of Business \& Economic Statistics} {\bf25}, 52--60.
%\endOrigBibText
\bptok{imsref}%
% NOT OUTPUTTED:
% number = 1
% doi = http://dx.doi.org/10.1198/073500106000000413
% fjournal = Journal of Business \& Economic Statistics
\endbibitem

%b5 ###
%b5 #&#
\bibitem{bs10}
%
\begin{bbook}[mr]
\bauthor{\bsnm{Bai},~\bfnm{Zhidong}\binits{Z.}} \AND
\bauthor{\bsnm{Silverstein},~\bfnm{Jack~W.}\binits{J.~W.}}
(\byear{2010}).
\btitle{Spectral Analysis of Large Dimensional Random Matrices},
\bedition{2nd} ed.
\bpublisher{Springer},
\blocation{New York}.
\bid{doi={10.1007/978-1-4419-0661-8}, mr={2567175}}
\end{bbook}
%
%\OrigBibText
%Bai, Z.D. and Silverstein, J.W. (2010). {\it Spectral Analysis of Large
%Dimensional Random Matrices}. Springer-Verlag, New York.
%\endOrigBibText
\bptok{imsref}%
% NOT OUTPUTTED:
% doi = http://dx.doi.org/10.1007/978-1-4419-0661-8
% isbn = 978-1-4419-0660-1
% fpage = xvi+551
\endbibitem



%b7 ###
%b7 #&#
\bibitem{bs06}
%
\begin{barticle}[mr]
\bauthor{\bsnm{Baik},~\bfnm{Jinho}\binits{J.}} \AND
\bauthor{\bsnm{Silverstein},~\bfnm{Jack~W.}\binits{J.~W.}}
(\byear{2006}).
\btitle{Eigenvalues of large sample covariance matrices of spiked
population models}.
\bjournal{J. Multivariate Anal.}
\bvolume{97}
\bpages{1382--1408}.
\bid{doi={10.1016/j.jmva.2005.08.003}, issn={0047-259X}, mr={2279680}}
\end{barticle}
%
%\OrigBibText
%Baik, J. and Silverstein, J.W. (2006). Eigenvalues of large sample covariance
%matrices of spiked population models. {\it Journal of Multivariate Analysis}
%{\bf97}, 1382--1408.
%\endOrigBibText
\bptok{imsref}%
% NOT OUTPUTTED:
% number = 6
% doi = http://dx.doi.org/10.1016/j.jmva.2005.08.003
% fjournal = Journal of Multivariate Analysis
\endbibitem

%b8 ###
%b8 #&#
\bibitem{bd92}
%
\begin{bbook}[mr]
\bauthor{\bsnm{Brockwell},~\bfnm{Peter~J.}\binits{P.~J.}} \AND
\bauthor{\bsnm{Davis},~\bfnm{Richard~A.}\binits{R.~A.}}
(\byear{1991}).
\btitle{Time Series: Theory and Methods},
\bedition{2nd} ed.
\bpublisher{Springer},
\blocation{New York}.
\bid{doi={10.1007/978-1-4419-0320-4}, mr={1093459}}
\end{bbook}
%
%\OrigBibText
%Brockwell, P.D. and Davis, R.A. (1991).
%{\it Time Series Analysis: Theory and Methods (2nd ed.).}
%Springer, New York.
%\endOrigBibText
\bptok{imsref}%
% NOT OUTPUTTED:
% doi = http://dx.doi.org/10.1007/978-1-4419-0320-4
% isbn = 0-387-97429-6
% fpage = xvi+577
\endbibitem

%b9 ###
%b9 #&#
\bibitem{c06}
%
\begin{barticle}[mr]
\bauthor{\bsnm{Chatterjee},~\bfnm{Sourav}\binits{S.}}
(\byear{2006}).
\btitle{A generalization of the {L}indeberg principle}.
\bjournal{Ann. Probab.}
\bvolume{34}
\bpages{2061--2076}.
\bid{doi={10.1214/009117906000000575}, issn={0091-1798}, mr={2294976}}
\end{barticle}
%
%\OrigBibText
%Chatterjee, S. (2006). A generalization of the Lindeberg principle.
%{\it The
%Annals of Probability} {\bf34}, 2061--2076.
%\endOrigBibText
\bptok{imsref}%
% NOT OUTPUTTED:
% number = 6
% doi = http://dx.doi.org/10.1214/009117906000000575
% coden = APBYAE
% fjournal = The Annals of Probability
\endbibitem

%b10 ###
%b10 #&#
\bibitem{e09}
%
\begin{barticle}[mr]
\bauthor{\bsnm{El Karoui},~\bfnm{Noureddine}\binits{N.}}
(\byear{2009}).
\btitle{Concentration of measure and spectra of random matrices:
Applications to correlation matrices, elliptical distributions and beyond}.
\bjournal{Ann. Appl. Probab.}
\bvolume{19}
\bpages{2362--2405}.
\bid{doi={10.1214/08-AAP548}, issn={1050-5164}, mr={2588248}}
\end{barticle}
%
%\OrigBibText
%El Karoui, N. (2009). Concentration of measure and spectra of random matrices:
%applications to correlation matrices, elliptical distributions and
%beyond. {\it
%The Annals of Applied Probability} {\bf19}, 2362--2405.
%\endOrigBibText
\bptok{imsref}%
% NOT OUTPUTTED:
% number = 6
% doi = http://dx.doi.org/10.1214/08-AAP548
% fjournal = The Annals of Applied Probability
\endbibitem

%b11 ###
%b11 #&#
\bibitem{fhlr00}
%
\begin{barticle}[auto:parserefs-M02]
\bauthor{\bsnm{Forni},~\bfnm{M.}\binits{M.}},
\bauthor{\bsnm{Hallin},~\bfnm{M.}\binits{M.}},
\bauthor{\bsnm{Lippi},~\bfnm{M.}\binits{M.}} \AND
\bauthor{\bsnm{Reichlin},~\bfnm{L.}\binits{L.}}
(\byear{2000}).
\btitle{The generalized factor model: Identification and estimation}.
\bjournal{Rev. Econ. Stat.}
\bvolume{82}
\bpages{540--554}.
\end{barticle}
%
%\OrigBibText
%Forni, M., Hallin, M., Lippi, M. and Reichlin, L. (2000).
%The generalized factor model: identification and estimation. {\it The
%Review of Economics
%and Statistics} {\bf82}, 540--554.
%\endOrigBibText
\bptok{imsref}%
\endbibitem

%b12 ###
%b12 #&#
\bibitem{fhlr04}
%
\begin{barticle}[mr]
\bauthor{\bsnm{Forni},~\bfnm{Mario}\binits{M.}},
\bauthor{\bsnm{Hallin},~\bfnm{Marc}\binits{M.}},
\bauthor{\bsnm{Lippi},~\bfnm{Marco}\binits{M.}} \AND
\bauthor{\bsnm{Reichlin},~\bfnm{Lucrezia}\binits{L.}}
(\byear{2004}).
\btitle{The generalized dynamic factor model: Consistency and rates}.
\bjournal{J. Econometrics}
\bvolume{119}
\bpages{231--255}.
\bid{doi={10.1016/S0304-4076(03)00196-9}, issn={0304-4076}, mr={2057100}}
\end{barticle}
%
%\OrigBibText
%Forni, M., Hallin, M., Lippi, M. and Reichlin, L. (2004). The generalized
%dynamic factor model: consistency and rates. {\it Journal of Econometrics}
%{\bf119}, 231--255.
%\endOrigBibText
\bptok{imsref}%
% NOT OUTPUTTED:
% number = 2
% doi = http://dx.doi.org/10.1016/S0304-4076(03)00196-9
% coden = JECMB6
% fjournal = Journal of Econometrics
\endbibitem

%b13 ###
%b13 #&#
\bibitem{fhlr05}
%
\begin{barticle}[mr]
\bauthor{\bsnm{Forni},~\bfnm{Mario}\binits{M.}},
\bauthor{\bsnm{Hallin},~\bfnm{Marc}\binits{M.}},
\bauthor{\bsnm{Lippi},~\bfnm{Marco}\binits{M.}} \AND
\bauthor{\bsnm{Reichlin},~\bfnm{Lucrezia}\binits{L.}}
(\byear{2005}).
\btitle{The generalized dynamic factor model: One-sided estimation and
forecasting}.
\bjournal{J. Amer. Statist. Assoc.}
\bvolume{100}
\bpages{830--840}.
\bid{doi={10.1198/016214504000002050}, issn={0162-1459}, mr={2201012}}
\end{barticle}
%
%\OrigBibText
%Forni, M., Hallin, M., Lippi, M. and Reichlin, L. (2005). The generalized
%dynamic factor model: one-sided estimation and forecasting.
%{\it Journal of the American Statistical Association} {\bf100}, 830--840.
%\endOrigBibText
\bptok{imsref}%
% NOT OUTPUTTED:
% number = 471
% doi = http://dx.doi.org/10.1198/016214504000002050
% coden = JSTNAL
% fjournal = Journal of the American Statistical Association
\endbibitem

%b14 ###
%b14 #&#
\bibitem{fl99}
%
\begin{barticle}[mr]
\bauthor{\bsnm{Forni},~\bfnm{Mario}\binits{M.}} \AND
\bauthor{\bsnm{Lippi},~\bfnm{Marco}\binits{M.}}
(\byear{1999}).
\btitle{Aggregation of linear dynamic microeconomic models}.
\bjournal{J.~Math. Econom.}
\bvolume{31}
\bpages{131--158}.
\bid{doi={10.1016/S0304-4068(98)00060-3}, issn={0304-4068}, mr={1676417}}
\end{barticle}
%
%\OrigBibText
%Forni, M. and Lippi, M. (1999). Aggregation of linear dynamic microeconomic
%models. {\it Journal of Mathematical Economics} {\bf31}, 131--158.
%\endOrigBibText
\bptok{imsref}%
% NOT OUTPUTTED:
% number = 1
% doi = http://dx.doi.org/10.1016/S0304-4068(98)00060-3
% coden = JMECDA
% fjournal = Journal of Mathematical Economics
\endbibitem

%b15 ###
%b15 #&#
\bibitem{fl01}
%
\begin{barticle}[mr]
\bauthor{\bsnm{Forni},~\bfnm{Mario}\binits{M.}} \AND
\bauthor{\bsnm{Lippi},~\bfnm{Marco}\binits{M.}}
(\byear{2001}).
\btitle{The generalized dynamic factor model: Representation theory}.
\bjournal{Econometric Theory}
\bvolume{17}
\bpages{1113--1141}.
\bid{issn={0266-4666}, mr={1867540}}
\end{barticle}
%
%\OrigBibText
%Forni, M. and Lippi, M. (2001). The generalized dynamic factor model:
%representation theory.
%{\it Econometric Theory} {\bf17}, 1113--1141.
%\endOrigBibText
\bptok{imsref}%
% NOT OUTPUTTED:
% number = 6
% fjournal = Econometric Theory
\endbibitem

%b16 ###
%b16 #&#
\bibitem{gh03}
%
\begin{barticle}[mr]
\bauthor{\bsnm{Geronimo},~\bfnm{Jeffrey~S.}\binits{J.~S.}} \AND
\bauthor{\bsnm{Hill},~\bfnm{Theodore~P.}\binits{T.~P.}}
(\byear{2003}).
\btitle{Necessary and sufficient condition that the limit of
{S}tieltjes transforms is a {S}tieltjes transform}.
\bjournal{J. Approx. Theory}
\bvolume{121}
\bpages{54--60}.
\bid{doi={10.1016/S0021-9045(02)00042-4}, issn={0021-9045}, mr={1962995}}
\end{barticle}
%
%\OrigBibText
%Geronimo, J.S. and Hill, T.P. (2003). Necessary and sufficient
%condition that
%the limit of Stieltjes transforms is a Stieltjes transform. {\it The
%Annals of
%Probability} {\bf31}, 54--60.
%\endOrigBibText
\bptok{imsref}%
% NOT OUTPUTTED:
% number = 1
% doi = http://dx.doi.org/10.1016/S0021-9045(02)00042-4
% coden = JAXTAZ
% fjournal = Journal of Approximation Theory
\endbibitem

%b17 ###
%b17 #&#
\bibitem{hln05}
%
\begin{barticle}[mr]
\bauthor{\bsnm{Hachem},~\bfnm{W.}\binits{W.}},
\bauthor{\bsnm{Loubaton},~\bfnm{P.}\binits{P.}} \AND
\bauthor{\bsnm{Najim},~\bfnm{J.}\binits{J.}}
(\byear{2005}).
\btitle{The empirical eigenvalue distribution of a Gram matrix: From
independence to stationarity}.
\bjournal{Markov Process. Related Fields}
\bvolume{11}
\bpages{629--648}.
\bid{issn={1024-2953}, mr={2191967}}
\end{barticle}
%
%\OrigBibText
%Hachem, W., Loubaton, P. and Najim, J. (2005). The empirical eigenvalue
%distribution of a Gram matrix: from independence to stationarity. {\it Markov
%Processes and Related Fields} {\bf11}, 629--648.
%\endOrigBibText
\bptok{imsref}%
% NOT OUTPUTTED:
% number = 4
% fjournal = Markov Processes and Related Fields
\endbibitem

%b18 ###
%b18 #&#
\bibitem{hln07}
%
\begin{barticle}[mr]
\bauthor{\bsnm{Hachem},~\bfnm{Walid}\binits{W.}},
\bauthor{\bsnm{Loubaton},~\bfnm{Philippe}\binits{P.}} \AND
\bauthor{\bsnm{Najim},~\bfnm{Jamal}\binits{J.}}
(\byear{2007}).
\btitle{Deterministic equivalents for certain functionals of large
random matrices}.
\bjournal{Ann. Appl. Probab.}
\bvolume{17}
\bpages{875--930}.
\bid{doi={10.1214/105051606000000925}, issn={1050-5164}, mr={2326235}}
\end{barticle}
%
%\OrigBibText
%Hachem, W., Loubaton, P. and Najim, J. (2007). Deterministic
%equivalents for certain
%functionals of large random matrices. {\it Annals of Applied
%Probability} {\bf17}, 875--930.
%\endOrigBibText
\bptok{imsref}%
% NOT OUTPUTTED:
% number = 3
% doi = http://dx.doi.org/10.1214/105051606000000925
% fjournal = The Annals of Applied Probability
\endbibitem

%b19 ###
%b19 #&#
\bibitem{hl07}
%
\begin{barticle}[mr]
\bauthor{\bsnm{Hallin},~\bfnm{Marc}\binits{M.}} \AND
\bauthor{\bsnm{Li{\u{s}}ka},~\bfnm{Roman}\binits{R.}}
(\byear{2007}).
\btitle{Determining the number of factors in the general dynamic
factor model}.
\bjournal{J. Amer. Statist. Assoc.}
\bvolume{102}
\bpages{603--617}.
\bid{doi={10.1198/016214506000001275}, issn={0162-1459}, mr={2325115}}
\end{barticle}
%%
%\OrigBibText
%Hallin, M. and Li\v{s}ka, R. (2007). Determining the number of
%factors in the general
%dynamic factor model. {\it Journal of the American Statistical Association}
%{\bf102}, 603--617.
%\endOrigBibText
\bptok{imsref}%
% NOT OUTPUTTED:
% number = 478
% doi = http://dx.doi.org/10.1198/016214506000001275
% coden = JSTNAL
% fjournal = Journal of the American Statistical Association
\endbibitem

%b20 ###
%b20 #&#
\bibitem{h94}
%
\begin{bbook}[mr]
\bauthor{\bsnm{Hamilton},~\bfnm{James~D.}\binits{J.~D.}}
(\byear{1994}).
\btitle{Time Series Analysis}.
\bpublisher{Princeton Univ. Press},
\blocation{Princeton, NJ}.
\bid{mr={1278033}}
\end{bbook}
%
%\OrigBibText
%Hamilton, J.D. (1994).
%{\it Time Series Analysis.}
%Princeton University Press, Princeton, NJ.
%\endOrigBibText
\bptok{imsref}%
% NOT OUTPUTTED:
% isbn = 0-691-04289-6
% fpage = xvi+799
\endbibitem

%b21 ###
%b21 #&#
\bibitem{jwbnh14}
%
\begin{barticle}[mr]
\bauthor{\bsnm{Jin},~\bfnm{Baisuo}\binits{B.}},
\bauthor{\bsnm{Wang},~\bfnm{Chen}\binits{C.}},
\bauthor{\bsnm{Bai},~\bfnm{Z.~D.}\binits{Z.~D.}},
\bauthor{\bsnm{Nair},~\bfnm{K.~Krishnan}\binits{K.~K.}} \AND
\bauthor{\bsnm{Harding},~\bfnm{Matthew}\binits{M.}}
(\byear{2014}).
\btitle{Limiting spectral distribution of a symmetrized auto-cross
covariance matrix}.
\bjournal{Ann. Appl. Probab.}
\bvolume{24}
\bpages{1199--1225}.
\bid{doi={10.1214/13-AAP945}, issn={1050-5164}, mr={3199984}}
\end{barticle}
%
%\OrigBibText
%Jin, B., Wang, C. Bai, Z.D. Nair, K.K. and Harding, M.C. (2014). Limiting
%spectral distribution of a symmetrized auto-cross covariance matrix.
%{\it Annals of Applied Probability} (to appear).
%\endOrigBibText
\bptok{imsref}%
% NOT OUTPUTTED:
% number = 3
% doi = http://dx.doi.org/10.1214/13-AAP945
% fjournal = The Annals of Applied Probability
\endbibitem

%b22 ###
%b22 #&#
\bibitem{jwml09}
%
\begin{barticle}[mr]
\bauthor{\bsnm{Jin},~\bfnm{Baisuo}\binits{B.}},
\bauthor{\bsnm{Wang},~\bfnm{Cheng}\binits{C.}},
\bauthor{\bsnm{Miao},~\bfnm{Baiqi}\binits{B.}} \AND
\bauthor{\bsnm{Lo Huang},~\bfnm{Mong-Na}\binits{M.-N.}}
(\byear{2009}).
\btitle{Limiting spectral distribution of large-dimensional sample
covariance matrices generated by {VARMA}}.
\bjournal{J. Multivariate Anal.}
\bvolume{100}
\bpages{2112--2125}.
\bid{doi={10.1016/j.jmva.2009.06.011}, issn={0047-259X}, mr={2543090}}
\end{barticle}
%
%\OrigBibText
%Jin, B., Wang, C., Miao, B. and Lo Huang, M.-N. (2009). Limiting spectral
%distribution of large-dimensional sample covariance matrices generated by
%VARMA. {\it Journal of Multivariate Analysis} {\bf100}, 2112--2125.
%\endOrigBibText
\bptok{imsref}%
% NOT OUTPUTTED:
% number = 9
% doi = http://dx.doi.org/10.1016/j.jmva.2009.06.011
% fjournal = Journal of Multivariate Analysis
\endbibitem

%b23 ###
%b23 #&#
\bibitem{j07}
%
\begin{bincollection}[mr]
\bauthor{\bsnm{Johnstone},~\bfnm{Iain~M.}\binits{I.~M.}}
(\byear{2007}).
\btitle{High dimensional statistical inference and random matrices}.
In \bbooktitle{International {C}ongress of {M}athematicians. {V}ol. {I}}
\bpages{307--333}.
\bpublisher{Eur. Math. Soc.},
\blocation{Z\"urich}.
\bid{doi={10.4171/022-1/13}, mr={2334195}}
\end{bincollection}
%%
%\OrigBibText
%Johnstone, I. M. (2007). High dimensional statistical inference and random
%matrices. In {\it Proceedings of the International Congress of Mathematicians
%I}, 307--333. European Mathematical Society, Zurich.
%\endOrigBibText
\bptok{imsref}%
% NOT OUTPUTTED:
% doi = http://dx.doi.org/10.4171/022-1/13
\endbibitem

%b24 ###
%b24 #&#
\bibitem{l13}
%
\begin{bmisc}[mr]
\bauthor{\bsnm{Liu},~\bfnm{Haoyang}\binits{H.}}
(\byear{2013}).
\bhowpublished{Spectral analysis of high dimensional time series.
Ph.D. thesis, Univ. California, Davis.}
\bid{mr={3232212}}
\end{bmisc}
%
%\OrigBibText
%Liu, H. (2013). {\it Spectral Analysis of High Dimensional Time
%Series}. Ph.D.
%Thesis. University of California, Davis.
%\endOrigBibText
\bptok{imsref}%
% NOT OUTPUTTED:
% isbn = 978-1303-79226-7
% url =
%%%http://gateway.proquest.com/openurl?url_ver=Z39.88-2004&rft_val_fmt=info:ofi/fmt:kev:mtx:dissertation&res_dat=xri:pqm&rft_dat=xri:pqdiss:3614237
% fpage = 173
\endbibitem

%b25 ###
%b25 #&#
\bibitem{lap-sm}
%
\begin{bmisc}[auto:parserefs-M02]
\bauthor{\bsnm{Liu},~\bfnm{H.}\binits{H.}},
\bauthor{\bsnm{Aue},~\bfnm{A.}\binits{A.}} \AND
\bauthor{\bsnm{Paul},~\bfnm{D.}\binits{D.}}
(\byear{2014}).
\bhowpublished{Supplement to ``On the Mar\v{c}enko--Pastur law for
linear time series.''
DOI:\doiurl{10.1214/14-AOS1294SUPP}}.
\bptok{imsref}%
\end{bmisc}
%
%\OrigBibText
%Liu, H., Aue, A. and Paul, D. (2014).
%Online supplementary material to ``On the Mar\v{c}enko--Pastur law for
%linear time series.
%Available online at...
%\endOrigBibText
\bptok{imsref}%
\endbibitem

%b26 ###
%b26 #&#
\bibitem{l06}
%
\begin{bbook}[mr]
\bauthor{\bsnm{L{\"u}tkepohl},~\bfnm{Helmut}\binits{H.}}
(\byear{2005}).
\btitle{New Introduction to Multiple Time Series Analysis}.
\bpublisher{Springer},
\blocation{Berlin}.
\bid{doi={10.1007/978-3-540-27752-1}, mr={2172368}}
\bptnote{check year}%
\end{bbook}
%
%\OrigBibText
%L\"utkepohl, H. (2006). {\it New Introduction to Multiple Time Series
%Analysis}. Springer-Verlag, New York.
%\endOrigBibText
\bptok{imsref}%
% NOT OUTPUTTED:
% doi = http://dx.doi.org/10.1007/978-3-540-27752-1
% isbn = 3-540-40172-5
% fpage = xxii+764
\endbibitem

%b27 ###
%b27 #&#
\bibitem{mp67}
%
\begin{barticle}[auto:parserefs-M02]
\bauthor{\bsnm{Mar{\v{c}}enko},~\bfnm{V.}\binits{V.}} \AND
\bauthor{\bsnm{Pastur},~\bfnm{L.}\binits{L.}}
(\byear{1967}).
\btitle{Distribution of eigenvalues for some sets of random matrices}.
\bjournal{Mathematics of the USSR. Sbornik}
\bvolume{1}
\bpages{457--483}.
\end{barticle}
%
%\OrigBibText
%Mar\v{c}enko, V. and Pastur, L. (1967). Distribution of eigenvalues
%for some
%sets of random matrices. {\it Mathematics of the USSR-Sbornik} {\bf1},
%457--483.
%\endOrigBibText
\bptok{imsref}%
\endbibitem

%b28 ###
%b28 #&#
\bibitem{m94}
%
\begin{barticle}[auto:parserefs-M02]
\bauthor{\bsnm{Molenaar},~\bfnm{P.~C.~M.}\binits{P.~C.~M.}}
(\byear{1994}).
\btitle{A manifesto on psycholog as idiographic science: Bringing the
person back into scientific psychology, this time forever}.
\bjournal{Measurement}
\bvolume{2}
\bpages{201--218}.
\end{barticle}
%%
%\OrigBibText
%Molenaar, P.C.M. (1994). A manifesto on psycholog as idiographic
%science: bringing the
%person back into scientific psychology, this time forever. {\it
%Measurement} {\bf2},
%201--218.
%\endOrigBibText
\bptok{imsref}%
\endbibitem

%b29 ###
%b29 #&#
\bibitem{ner92}
%
\begin{barticle}[auto:parserefs-M02]
\bauthor{\bsnm{Ng},~\bfnm{V.}\binits{V.}},
\bauthor{\bsnm{Engel},~\bfnm{R.~F.}\binits{R.~F.}} \AND
\bauthor{\bsnm{Rothschild},~\bfnm{M.}\binits{M.}}
(\byear{1992}).
\btitle{A multi-dynamic-factor model for stock returns}.
\bjournal{J. Econometrics}
\bvolume{52}
\bpages{245--266}.
\end{barticle}
%
%\OrigBibText
%Ng, V., Engel, R.F. and Rothschild, M. (1992). A
%multi-dynamic-factor model for stock returns.
%{\it Journal of Econometrics} {\bf52}, 245--266.
%\endOrigBibText
\bptok{imsref}%
\endbibitem

%b30 ###
%b30 #&#
\bibitem{o12}
%
\begin{barticle}[mr]
\bauthor{\bsnm{Onatski},~\bfnm{Alexei}\binits{A.}}
(\byear{2012}).
\btitle{Asymptotics of the principal components estimator of large
factor models with weakly influential factors}.
\bjournal{J. Econometrics}
\bvolume{168}
\bpages{244--258}.
\bid{doi={10.1016/j.jeconom.2012.01.034}, issn={0304-4076}, mr={2923766}}
\end{barticle}
%
%\OrigBibText
%Onatski, A. (2012). Asymptotics of the principal components estimator
%of large
%factor models with weakly influential factors. {\it Journal of Econometrics}
%{\bf168}, 244--258.
%\endOrigBibText
\bptok{imsref}%
% NOT OUTPUTTED:
% number = 2
% doi = http://dx.doi.org/10.1016/j.jeconom.2012.01.034
% coden = JECMB6
% fjournal = Journal of Econometrics
\endbibitem

%b31 ###
%b31 #&#
\bibitem{pa13}
%
\begin{barticle}[mr]
\bauthor{\bsnm{Paul},~\bfnm{Debashis}\binits{D.}} \AND
\bauthor{\bsnm{Aue},~\bfnm{Alexander}\binits{A.}}
(\byear{2014}).
\btitle{Random matrix theory in statistics: A review}.
\bjournal{J. Statist. Plann. Inference}
\bvolume{150}
\bpages{1--29}.
\bid{doi={10.1016/j.jspi.2013.09.005}, issn={0378-3758}, mr={3206718}}
\end{barticle}
%
%\OrigBibText
%Paul, D. and Aue, A. (2014). Random matrix theory in statistics: a review.
%{\it Journal of Statistical Planning and Inference} {\bf150}, 1--29.
%\endOrigBibText
\bptok{imsref}%
% NOT OUTPUTTED:
% doi = http://dx.doi.org/10.1016/j.jspi.2013.09.005
% fjournal = Journal of Statistical Planning and Inference
\endbibitem

%b32 ###
%b32 #&#
\bibitem{ps09}
%
\begin{barticle}[mr]
\bauthor{\bsnm{Paul},~\bfnm{Debashis}\binits{D.}} \AND
\bauthor{\bsnm{Silverstein},~\bfnm{Jack~W.}\binits{J.~W.}}
(\byear{2009}).
\btitle{No eigenvalues outside the support of the limiting empirical
spectral distribution of a separable covariance matrix}.
\bjournal{J. Multivariate Anal.}
\bvolume{100}
\bpages{37--57}.
\bid{doi={10.1016/j.jmva.2008.03.010}, issn={0047-259X}, mr={2460475}}
\end{barticle}
%%
%\OrigBibText
%Paul, D. and Silverstein, J.W. (2009). No eigenvalues outside the
%support of
%the limiting empirical spectral distribution of a separable covariance matrix.
%{\it Journal of Multivariate Analysis} {\bf100}, 37--57.
%\endOrigBibText
\bptok{imsref}%
% NOT OUTPUTTED:
% number = 1
% doi = http://dx.doi.org/10.1016/j.jmva.2008.03.010
% fjournal = Journal of Multivariate Analysis
\endbibitem

%b33 ###
%b33 #&#
\bibitem{ps12}
%
\begin{barticle}[mr]
\bauthor{\bsnm{Pfaffel},~\bfnm{Oliver}\binits{O.}} \AND
\bauthor{\bsnm{Schlemm},~\bfnm{Eckhard}\binits{E.}}
(\byear{2011}).
\btitle{Eigenvalue distribution of large sample covariance matrices of
linear processes}.
\bjournal{Probab. Math. Statist.}
\bvolume{31}
\bpages{313--329}.
\bid{issn={0208-4147}, mr={2853681}}
\bptnote{check year}%
\end{barticle}
%
%\OrigBibText
%Pfaffel, O. and Schlemm, E. (2012). Eigenvalue distribution of large sample
%covariance matrices of linear processes. {\it arXiv:1201.3828}.
%\endOrigBibText
\bptok{imsref}%
% NOT OUTPUTTED:
% number = 2
% fjournal = Probability and Mathematical Statistics
\endbibitem

%b34 ###
%b34 #&#
\bibitem{sw05}
%
\begin{bmisc}[auto:parserefs-M02]
\bauthor{\bsnm{Stock},~\bfnm{J.~H.}\binits{J.~H.}} \AND
\bauthor{\bsnm{Watson},~\bfnm{M.~W.}\binits{M.~W.}}
(\byear{2005}).
\bhowpublished{Implications of dynamic factor models for VAR analysis.
NBER Working Paper No. 11467.}
\end{bmisc}
%
%\OrigBibText
%Stock, J.H. and Watson, M.W. (2005). Implications of dynamic factor
%models for VAR analysis.
%{\it NBER Working Paper No. 11467}.
%\endOrigBibText
\bptok{imsref}%
\endbibitem

%b35 ###
%b35 #&#
\bibitem{t12}
%
\begin{bbook}[mr]
\bauthor{\bsnm{Tao},~\bfnm{Terence}\binits{T.}}
(\byear{2012}).
\btitle{Topics in Random Matrix Theory}.
\bseries{Graduate Studies in Mathematics}
\bvolume{132}.
\bpublisher{Amer. Math. Soc.},
\blocation{Providence, RI}.
\bid{mr={2906465}}
\end{bbook}
%
%\OrigBibText
%Tao, T. (2012). {\it Topics in Random Matrix Theory}. American Mathematical
%Society.
%\endOrigBibText
\bptok{imsref}%
% NOT OUTPUTTED:
% isbn = 978-0-8218-7430-1
% fpage = x+282
\endbibitem

%b36 ###
%b36 #&#
\bibitem{y12}
%
\begin{barticle}[mr]
\bauthor{\bsnm{Yao},~\bfnm{Jianfeng}\binits{J.}}
(\byear{2012}).
\btitle{A note on a {M}ar\v{c}enko--{P}astur type theorem for time series}.
\bjournal{Statist. Probab. Lett.}
\bvolume{82}
\bpages{22--28}.
\bid{doi={10.1016/j.spl.2011.08.011}, issn={0167-7152}, mr={2863018}}
\end{barticle}
%
%\OrigBibText
%Yao, J.-F. (2012). A note on a {M}ar\v{c}enko-{P}astur type theorem
%for time
%series. {\it Statistics \& Probability Letters} {\bf82}, 22--28.
%\endOrigBibText
\bptok{imsref}%
% NOT OUTPUTTED:
% number = 1
% doi = http://dx.doi.org/10.1016/j.spl.2011.08.011
% coden = SPLTDC
% fjournal = Statistics \& Probability Letters
\endbibitem

%b37 ###
%b37 #&#
\bibitem{z06}
%
\begin{bmisc}[auto:parserefs-M02]
\bauthor{\bsnm{Zhang},~\bfnm{L.}\binits{L.}}
(\byear{2006}).
\bhowpublished{Spectral analysis of large dimensional random matrices.
Ph.D. thesis, National Univ. Singapore}.
\end{bmisc}
%
%\OrigBibText
%Zhang, L. (2006). {\it Spectral Analysis of Large Dimensional Random
%Matrices}. Ph.D. Thesis. National University of Singapore.
%\endOrigBibText
\bptok{imsref}%
\endbibitem
\end{thebibliography}
\end{document}